\documentclass[12pt,a4paper]{amsart}
\usepackage{amsmath,amscd,amssymb,amsthm,array}
\usepackage{stix}

\usepackage{cancel}
\usepackage{mathtools}

\usepackage{stmaryrd}

\usepackage{hyperref}

\usepackage{tikz-cd}

\newcommand {\Lra} {\Longrightarrow}

\usepackage{hyperref}

\let\Oproof\proof
\let\Oendproof\endproof
\let\proof\relzx
\let\endproof\relax

\usepackage{DLdef1}

\makeatletter
\theoremstyle{thm}
\@ifundefined{Theorem}{}{}
\@ifundefined{Lemma}{}{}
\@ifundefined{Proposition}{}{}
\@ifundefined{Statement}{}{}
\@ifundefined{Corollary}{}{}
\@ifundefined{Problem}{}{}
\@ifundefined{Question}{\newtheorem*{Question}{\secno Question}}{}
\@ifundefined{Conjecture}{}{}
\theoremstyle{definition}
\@ifundefined{Remark}{}{}
\@ifundefined{Example}{}{}
\@ifundefined{Convention}{}{}
\@ifundefined{rank}{\DeclareMathOperator{\rank}{rank}}{}
\@ifundefined{qtr}{\DeclareMathOperator{\qtr}{qtr}}{}
\@ifundefined{Cent}{\DeclareMathOperator{\Cent}{Cent}}{}
\makeatother

\usepackage{tabularx,mdwtab}
\let\proof\Oproof
\let\endproof\Oendproof

\renewcommand {\ssbegin}[2][*]
 {\refstepcounter{subsection}%
\if#1*
\addcontentsline{toc}{subsection}{\thesubsection.\hskip 1pc #2}%
\else
\addcontentsline{toc}{subsection}{\thesubsection.\hskip 1pc #2. #1}%
\fi
 \def \secno {\gdef \secno {}{\ssecfont
\thesubsection.\hskip 2ex}%
 }%
 \begin{#2}}

\renewcommand {\sssbegin}[2][*]
 {\refstepcounter{subsubsection}
\if#1*
\addcontentsline{toc}{subsubsection}{\thesubsubsection.\hskip 1pc #2}%
\else
\addcontentsline{toc}{subsubsection}{\thesubsubsection.\hskip 1pc #2. #1}
\fi
 \def \secno {\gdef \secno {}{\ssecfont \thesubsubsection.\hskip 2ex}%
 }%
 \begin{#2}}

\renewcommand {\parbegin}[2][*]
 {\refstepcounter{paragraph}
\if#1*
\addcontentsline{toc}{paragraph}{\theparagraph.\hskip 1pc #2}%
\else
\addcontentsline{toc}{paragraph}{\theparagraph.\hskip 1pc #2. #1}
\fi
 \def \secno {\gdef \secno {}{\ssecfont \theparagraph.\hskip 2ex}%
 }%
 \begin{#2}}

\renewcommand {\thesubsection} {\arabic{section}.\arabic{subsection}}
\renewcommand {\thesubsubsection} {\thesubsection.\arabic{subsubsection}}

\renewcommand\thetable{\arabic{subsection}.\arabic{table}}

\newcommand {\Llra} {\Longleftrightarrow}

\makeatletter

\makeatletter
\let\JS@makepapertitle\relax
\def\JS@makepapertitle{
  \JS@doublepagebreak
  \global\JS@emptyheadlinetrue
  \global\JS@papertitlepagetrue
  \lhead{\ifJS@emptyheadline
      \ifJS@papertitlepage
        \JS@headlinesmallfont
        \rlap{\smash{\hskip-50pt\lower4pt\hbox{\logoDG}}}%
      \hbox{\vbox to 0pt{\vskip3.5pt\hbox{%
        \theshort@@journalname\ \the@VOLUME\ (\the@YEAR)
        \JS@firstpage--\JS@lastpage}
\JS@DOIinheadline
            \hbox{\vbox{\textcopyright \the@YEAR\ \theshort@author\\This is an open access article licensed under the CC BY-SA 4.0}}
          \vss}}%
      \fi
    \else
      \ifodd\pageno\else\the\pageno\fi
    \fi}
  \chead{\ifJS@emptyheadline
    \else
      \JS@headlinesmallfont
      \ifodd\pageno\theshort@title\else\theshort@author\fi
    \fi}
  \rhead{\ifJS@emptyheadline
      \ifJS@papertitlepage
        \the\pageno
      \fi
      \global\JS@emptyheadlinefalse
      \global\JS@papertitlepagefalse
    \else
      \ifodd\pageno\the\pageno\fi
    \fi}
  \vspace*{1cm}
  \JS@fnornothing\the@MscName\the@msc
  \JS@fnornothing\the@KeywordsName\the@keywords
 \JS@fnornothing{\hskip-0.5em}\the@acknowledgment
  \JS@fnornothing\the@AffiliationName\the@affiliation
  \begin{flushleft}
  {\baselineskip18pt \lineskiplimit6pt \lineskip6pt
  \noindent
  \vrule height1cm width0pt\relax
  \sffamily\bfseries\mathversion{bold}\Large\the@title\par}
  \JS@ifempty\the@author\else
    \baselineskip30pt
    \noindent
    \em\fbox{\the@author}
    \vskip24pt
  \fi
    \JS@ifempty\the@dedication\else
    \baselineskip12pt
    \noindent
    \sl\the@dedication
    \vskip24pt
  \fi
  \end{flushleft}
  \JS@narrowparornothing\the@AbstractName\the@abstract%
  }

\makeatother

\title{How to superize $\mathfrak{gl}(\infty)$}

\author{George Egorov} 

\address{c/o Dimitry Leites, Department of Mathematics, Stockholm University, Albanov\"agen 28, SE-114 19 Stockholm, Sweden; dimleites@gmail.com}

\makeatletter
\@namedef{subjclassname@2020}{\textup{2020} Mathematics Subject Classification}
\makeatother


\begin{document}

\begin{abstract} Penkov with co-authors studied several versions of the Lie algebra $\mathfrak{gl}(n)$, where $n=\infty$. Completely different new types of super versions of $\mathfrak{gl}(\infty)$ are introduced in this paper: these are Lie superalgebras of supermatrices infinite in all directions, the non-zero entries of each matrix being confined to a region around the main diagonal. For three of the four superalgebras introduced here, that region is a~band whose boundaries are curves, not straight lines: the band is of a~fixed width in one case, and in the other two it flares out at infinity, slower than any power in one of them and slower than the distance to the origin in the other. The fourth and largest superalgebra is not confined inside a~band: its definition describes only the quadrants the support may enter, the region itself is cut out by no inequality between $|i-j|$ and $|i+j|$. The largest of the four turns out to be precisely the space of matrices on which the cocycle defining the central extension is given by a~finite sum, and this largest algebra lies properly inside the Mackey Lie algebra of Penkov and Serganova, on which the same sum is infinite and the cocycle therefore undefined. Several open questions related with further study of the Lie superalgebras introduced here and with their possible applications to dynamical systems, such as KdV, are formulated.
\end{abstract}

\maketitle

\markboth{\itshape George Egorov}{{\itshape How to superize $\fgl(\infty )$}}



\section{Introduction}

\textit{Editor's note}. After G.~Egorov died, I composed an article from scrap papers he left to me; he intended to improve them as a~ contribution to proceedings of my ``Seminar on Supersymmetries\rq\rq\  a.k.a. SoS, but his sudden death intervened. Egorov was almost blind, so it was sometimes difficult to decipher intersecting lines of the notes he wrote to himself, or recover signs omitted. Being guided by what I recalled of our discussions, I concocted the abstract and introduction, and supplied with a~bibliography, but was unable to recover some of the proofs. Since that article was published in not easily accessible collection, see \cite{Eg}, for the current version of Egorov's paper I updated the references trying to provide with the latest sources containing definitions of the basics (e.g., Cartan matrix, Dynkin-Kac diagram, double extension) more reasonable than the first hasty and clumsy definitions and more convenient notation. Questions from Irina Shchepochkina, Alexey Lebedev,  Ivan Penkov and Alex Petukhov considerably improved (corrected)  the draft,  to them and O.~Lozhechnyk at the last stages I am very thankful for help. I heartily thank Osmo Pekonen who sent me the pdf of \cite{Eg}, so one could work on it trying to make Egorov's findings more available and apply them to dynamical systems, such as KdV. 
~\textit{D.Leites.} 
\smallskip

In this paper, the ground field is $\Cee$.

\subsection{Structures of a Lie algebra and of an associative algebra on the same space}\label{ssLiealg} We write $\Mat$ speaking about the associative structure with matrix product and $\fgl$ when the matrix product is replaced with a~bracket. The two structures live on one and the same space, and each of the superalgebras of this paper carries both, see Lemma \ref{Lideal}\textup{;} which of the two is meant matters at exactly one place, Convention \ref{2.3.1}, whose price is settled in the Appendix, see \ref{Sapp}. 

There are many contenders for the role of $\fgl (n)$ in the case where $n =\infty$. Lately, some of these contenders are distinguished thanks to their applications in the theory of differential equations such as KdV and KP hierarchies, see~\cite{R} and references therein. For other applications of $\fgl (\infty )$, see \cite{B,K,KL1, KL2,FF,V}. (The best, in my opinion, short review on KdV is an MPIM-Bonn (1987) preprint by Date, later developed into \cite{MJD}; see also \cite{V, Ma}.)

Let me list some of the Lie algebras of infinite matrices that enjoy everybody's recognition. It goes without saying that we identify a~linear transformation with its matrix having chosen a~basis in the space $V$ of the tautological representation of the Lie algebra $\fgl (V)$. Thus, we replace the Lie algebra of operators $\fgl (V)$ by the Lie algebra of matrices $\fgl (\dim V)$ and for infinite-dimensional (of countable dimension) $V$, so every element of $V$ is a finite sum of basis elements, we get\footnote{A.~Radul observed that $\fgl _c(\infty )$ is often denoted in the literature by $\fgl_{J} (\infty )$ hinting at the operator $J$ of complex structure equal to $\id$ on $z$ in non-negative degrees and $-\id$ on $z$ in negative degrees (and at Jacobi in the case of 3-diagonal (Jacobi) matrices).}:
\[
\begin{array}{ll}
\fgl_{cf} (\infty ) &= \{X \mid X\text{~~ is a~linear endomorphism of $\Cee^\infty=\Cee^\Zee$} \}\\
&= \{X\text{~~ is a~matrix with finitely many non-zero entries in every column} \};\\
\fgl _c(\infty ) &= \{X \in \fgl (\infty ) \mid \text{~~ nonzero entries of $X$ belong to a~band}\\ &\text{whose boundaries are parallel to the main diagonal}\};\\
\fgl _f(\infty ) &= \{X \in \fgl (\infty ) \mid X\text{~~ has finitely many nonzero entries}\}.\\
\end{array}
\]

Each of these three Lie algebras of infinite both ways matrices has a~subalgebra of matrices whose rows and columns are labeled with only positive integers. Denote the respective Lie subalgebras by $\fgl_{cf} ^{+} (\infty )$, $\fgl_c ^{+} (\infty )$, $\fgl_f ^{+} (\infty )$.

Let us explain why these algebras are so natural. For the definition of Lie (super)algebras with Cartan matrix and the corresponding Dynkin--Kac diagrams, see \cite{CCLL}; for realizations of finite-dimensional Lie superalgebras, see \cite{BGLLS}. Now, look at the Dynkin diagram of $\fgl (n)$; more precisely, of $\fsl(n)$ as $n\tto \infty $. There are two ways to tend $n$ to $\infty $ and that is how these algebras are obtained: they are generated by countably many Chevalley generators subject to the relations governed by the Cartan matrix corresponding to one of the following infinite Dynkin diagrams (the vertices of one are numbered by positive integers, the vertices of the other one are numbered by all integers):
\[
\begin{array}{ll}
\bigcirc -\bigcirc -\bigcirc -\ldots &\text{(corresponds to $\fsl _f^{+} (\infty )$)}\\
\ldots -\bigcirc -\bigcirc -\bigcirc -\ldots & \text{(corresponds to $\fsl _f (\infty )$)}
\end{array}
\]
\textbf{Observe that although ostensibly different, the limit algebras are isomorphic}. Penkov proved this by a simple argument: their elements are transformations of the space whose basis vectors are labeled by all (resp., only positive) integers; the isomorphism of the sets of labels establishes isomorphism $\fgl _f^{+} (\infty )\simeq\fgl _f (\infty )$ \textbf{which does not come from any analog of the Weyl group} considered in what follows. 

For applications, however, these algebras turn out to be too small and it is desirable to enlarge them so that the enlarged Lie algebra has a~nontrivial central extension of degree 0. (This central extension is needed in ``quantization\rq\rq, cf.~\cite{LSh}.)


The above-defined intermediate Lie algebras $\fgl _c(\infty )$ and $\fgl _c^{+}(\infty )$ prove to be just fine, although no analog of Dynkin diagram corresponds to them. (These intermediate algebras themselves have interesting enlargements. R.~Bullough has reported, in a~talk, that he found a~need in an enlargement of the conventional $\fgl _c(\infty )$. Hopefully, the even parts of the Lie superalgebras introduced in this paper will satisfy these needs; for what such a~need can and cannot mean, see Subsec. \ref{ssKP}.) The central extension we have mentioned is given by the cocycle (see \cite{V, MJD}, and \cite[Appendix~B]{Fu})
\be\label{tr}
X,Y \mapsto \tr (J[X,Y]J- [JXJ,JYJ])\quad \text{for any}\quad X,Y \in \fgl _c (\infty )
\ee
where $J = \sum_{i\geq 0} E_{ii}$.  Since the argument of the trace has finite support, the trace is well defined. 

\subsection{Which $\fgl (\infty )$ the KdV and KP hierarchies use, and what an
enlargement is for}\label{ssKP}
The contender singled out by the soliton literature is the algebra of matrices
with finitely many non-zero \textit{diagonals},
\be\label{aInf}
\bar\fgl (\infty ):=\{(a_{ij})_{i,j\in \Zee}\mid a_{ij}=0\ \text{ for }\ |i-j|\gg 0\},
\ee
written $\bar a_\infty $ in \cite{KR}; it is the $\fgl _{p,c}$ of Subsec. 2.1 with
$p\equiv \ev $, and what carries the theory is its central extension by the
cocycle \eqref{tr}. The Fock module over that extension is where the bilinear
identity, the Sato Grassmannian and the tau function live, see
\cite{MJD,KR,V,AZ}\textup{;} the reduction of KP to the Gelfand--Dickey
hierarchies is the passage to the affine subalgebras $\widehat{\fsl (n)}\subset
\bar a_\infty $ spanned by the $n$-periodic matrices, and $n=2$ is KdV. In the
super case the same choice is made: \cite{KL1,KL2,KL3} work with
$\fgl _c(\infty |\infty )$, which is our $\fgl _{p,c}$ for $p=p_{\text{st}}$. The
analytic counterpart of \eqref{aInf}, used when one wants all the solutions and
not only the formal ones, is the algebra $\fgl _{\text{res}}$ of \cite{SW}:
the bounded operators $X$ in a~Hilbert space with a~fixed polarization $H=H_+\oplus H_-$ such that the two blocks of $X$ that do not respect the polarization, i.e., the components $X_{+-}\colon H_- \tto H_+$ and $X_{-+}\colon H_+ \tto H_-$, are Hilbert--Schmidt. A~third algebra of the same circle of ideas, the differential
operators on the circle \cite{KRad}, is a~subalgebra of \eqref{aInf}, its matrix
having polynomial entries on finitely many diagonals.

Thus, the soliton literature has enlarged $\fgl (\infty )$ once already --- from
the finitary $\fgl _f(\infty )$ to \eqref{aInf} --- and that enlargement is not a
matter of taste.

Two terms are used throughout, and here they are. A~matrix $a=(a_{ij})_{i,j\in \Zee}$ is a~\textit{band matrix}, or is \textit{banded}, if its non-zero entries lie on finitely many diagonals, i.e., if
\be\label{banded}
\text{$a_{ij}=0$ whenever $|i-j|>c_a$ for some $c_a>0$ depending on $a$;}
\ee
the set of all band matrices is closed under the matrix product and under the bracket, and \textit{the band algebra} is this set with either structure --- it is $\bar\fgl (\infty )$ of eq.~\eqref{aInf} for $p\equiv \ev $, and $\fgl _{p,c}$ of Subsec. $2.1$ in general.

The \textit{corner blocks} of a~matrix $X$ with respect to the projection $J:=\sum_{i\geq 0}E_{ii}$, i.e., with respect to the splitting $V=V_-\oplus V_+$ into the span of the $e_i$ with $i<0$ and the span of the $e_i$ with $i\geq 0$, are the two off-diagonal components
\be\label{corners}
X=\begin{pmatrix} X_{++}&X_{+-}\\ X_{-+}&X_{--}\end{pmatrix},\qquad
X_{++}:=JXJ,\quad X_{+-}:=JX(\One -J),\quad X_{-+}:=(\One -J)XJ,
\ee
where, as everywhere below, $+$ stands for the indices $i\geq 0$ and $-$ for the indices $i<0$. So $X_{+-}$ is the block of the entries $(i,j)$ with $i\geq 0>j$ and $X_{-+}$ that of the entries with $i<0\leq j$.

\ssbegin{Lemma}[The band algebra is the smallest on which the cocycle
survives]\label{Lcob} Let $X_{+-}$ and $X_{-+}$ be the corner blocks \eqref{corners} of $X$, and let
\be\label{ccoc}
c(X,Y):=\tr \big(X_{+-}Y_{-+}\big)-\tr \big(Y_{+-}X_{-+}\big)
\ee
be the cocycle \eqref{tr}.

\textup{1)} On $\fgl _f(\infty )$, the cocycle $c$ is the coboundary of the
linear form $f(X):=\tr (JX)$, and hence defines no non-trivial central extension.

\textup{2)} On $\fgl _c(\infty )$, that form does not exist, and $c$ is not a
coboundary\textup{:} let 
\be\label{Lambda}
\Lambda :=\sum_{i\in \Zee}E_{i,i+1}
\ee
be the shift, so
that $\Lambda $ is invertible and $\Lambda ^k=\sum_{i\in \Zee}E_{i,i+k}$ for every
$k\in \Zee$. Then, $[\Lambda ^k,\Lambda ^{-k}]=0$, whereas
$c(\Lambda ^k,\Lambda ^{-k})=-k$ for $k\geq 1$.
\end{Lemma}

\begin{proof} $1)$ Since $[J,X]=X_{+-}-X_{-+}$, we get
\[
f([X,Y])=\tr \big(J(XY-YX)\big)=\tr \big([J,X]Y\big)
=\tr \big(X_{+-}Y_{-+}\big)-\tr \big(Y_{+-}X_{-+}\big)=c(X,Y),
\]
every trace here being a~finite sum because $X$ and $Y$ are of finite support.

$2)$ The sum $f(\One )=\sum_{i\geq 0}1$ has infinitely many non-zero summands and hence no value, so $f$ is not defined on
$\fgl _c(\infty )$. A~coboundary vanishes on a~commuting pair; the cocycle does
not. Indeed, for $k\geq 1$ the block $(\Lambda ^k)_{+-}$ is zero, the block
$(\Lambda ^{-k})_{+-}$ consists of the $k$ entries $(i,i-k)$ with $0\leq i\leq
k-1$, and $(\Lambda ^k)_{-+}$ of the $k$ entries $(i,i+k)$ with $-k\leq i\leq -1$.
Hence, \eqref{ccoc} gives $c(\Lambda ^k,\Lambda ^{-k})=0-k=-k$.
\end{proof}

\sssbegin{Remark}[What an enlargement can, and what it cannot, be for]\label{Renl}
An enlargement of $\fgl _c(\infty )$ cannot be wanted in order to accommodate the
\textit{group}. Indeed, no condition on the support of a~matrix, however
generous, is preserved by the operations the hierarchy performs: both the resolvent
${(\One -\Lambda )^{-1}=\sum_{0\leq k}\Lambda ^k}$ of the shift $\Lambda$, see eq.~\eqref{Lambda}, and the KP time flow
$\exp \big(\sum_{1\leq k}t_k\Lambda ^k\big)$ have all their diagonals with
$k\geq 0$ non-zero, so their support is a~half-plane and they lie neither in any
of the four superalgebras of \S~2 nor even in the Mackey Lie algebra of Corollary
\ref{CMackey}. What the theory does instead is to exponentiate not the matrix but
its action: on the semi-infinite wedge module the action of a~matrix from
\eqref{aInf} is locally finite, so the exponential is defined there while the
matrix of the group element is not.

An enlargement is wanted, therefore, for the two things that do live on the
algebra itself. The first is the cocycle, and Lemma \ref{Lcob} shows that it is
the reason the finitary algebra had to be abandoned; Proposition \ref{Pcorner}
below determines how far in that direction one may go, and the answer is:
exactly as far as $\fgl _{p,g}$. The second is the set of bases: over
$\fgl _c(\infty |\infty )$ there are many classes of them, see Theorem \ref{Th},
whereas over $\fgl _{p,g}$ there is just one. Since a~super KdV or KP hierarchy
is read off a~base by the superization of the Drinfeld--Sokolov procedure
described in \cite{LSS}, see Subsec. \ref{ssProb} here, and carried out in detail later in
\cite{IK, DG}, the number of classes of bases is the number of hierarchies one
gets, and an
algebra with just one class is what removes the ambiguity. This is the sense in
which the superalgebras of this paper answer a~need for an enlargement.
\end{Remark}

\sssbegin{Remark}[Penkov's algebras in the soliton literature]\label{RPenkov} As
far as we could establish, the versions of $\fgl (\infty )$ of \cite{PS1,ChP1,HP}
have entered the theory of the KdV and KP hierarchies exactly once, in
\cite{FPS}, and there the algebra is the \textit{finitary} $\fsl (\infty )$ and
the object categorified is the Fock space of the boson--fermion correspondence,
not the hierarchy: the equations themselves are still written over \eqref{aInf}.
Corollary \ref{CMackey} below tells why the largest of those versions cannot
simply replace \eqref{aInf}.
\end{Remark}

\textbf{The Lie algebras and Lie superalgebras studied by Penkov with co-authors} See \cite{PS1, PS2, FPS, PP1, PP2, PP3, PP4, IP, IPW, PSZ, PZ, ChP1, ChP2, HPS, CouP, GP, Pe, CaP}, and an amendment in \cite{HZ}.  For a summary, see \cite{HP}.

\subsection{Bases (systems of simple roots) of Lie superalgebras} As is known, some~simple finite-dimen\-sio\-nal Lie superalgebras $\fg$ can have several bases (a.k.a. systems of simple roots) not isomorphic with respect to the Weyl group of $\fg_{\overline 0}$. The simplest way to see the reason for this phenomenon is to consider $\fgl (V)$ for a~finite-dimensional superspace $V=V_\ev\oplus V_\od$.

The matrix expression of elements from $\fgl (V)$ now depends on the basis of $V$ (which is always supposed to consist of homogeneous elements). The general linear Lie superalgebra of $n \times n$ matrices is determined by a~ set of  basis vectors and the set \text{Par} of their parities, of which $\dim V_\ev$ are even and $\dim V_\od$ are odd, so that $\dim V_\ev+\dim V_\od = n$. 

\subsubsection{Super Weyl groups}
V. Serganova suggested several generalizations of the Weyl group all of which act by permutations of all the bases of a~simple Lie algebra $\fg$ (either finite-dimensional or a~(twisted) loop one) with Cartan matrix~\cite{S1,S2,LSS}. Let us recall some of them. For  history and various approaches, see \cite{BLLS}; for the latest developments, see \cite{GH, GHS, GHSS, HW}.

\begin{table} 
\makeatletter
\renewcommand\thetable{\arabic{table}}
\renewcommand\@makecaption[2]{%
\vskip \abovecaptionskip \sbox \@tempboxa {\textbf{#1}. #2}\ifdim%
\wd \@tempboxa >\hsize \textbf{#1}. #2\par \else \global \@minipagefalse \hb@xt@ \hsize%
{\hfil \box \@tempboxa \hfil }\fi \vskip \belowcaptionskip}
\makeatother
\caption{Coxeter diagrams of universal super Weyl groups}

Notation: $\otimes $ stands for an odd reflection.

\medskip

\begin{center}
\small
\begin{tabular}{|c|c|c|c|}
\hline
$\fg$&$\fsl(m | n)$&$\fosp(3 | 2)$&$\fosp(3 | 2n)$; $n> 1$\\
\hline
$UW_{\fg}$&$\underbrace {\bigcirc -\cdots-\bigcirc }_{m-1}\overset \infty {-}\otimes \overset \infty {-}\underbrace {\bigcirc -\cdots-\bigcirc }_{n-1}$&$\bigcirc \overset \infty {-}\otimes \overset \infty {-}\bigcirc $&$\bigcirc \overset \infty {-}\otimes \overset \infty {-}\underbrace {\bigcirc -\cdots-\bigcirc \overset 4{-}\bigcirc }_n$\\
\hline
\end{tabular}\label{T1}


\begin{tabular}{|c|c|c|c|}
\hline
$\fg$&$\fosp(2m+1 | 2)$, $m>1$&$\fosp_\alpha(4 | 2)$&$\fosp(2 | 2n)$, $n>1$\\
\hline
$UW_{\fg}$&$\underbrace {\bigcirc \overset 4{-}\bigcirc -\cdots-\bigcirc }_m \overset \infty {-}\otimes \overset \infty {-}\bigcirc $&$\underset {\displaystyle \overset {\displaystyle \hphantom {\scriptstyle \infty }|{\scriptstyle \infty }}{\displaystyle \bigcirc }}{\bigcirc \overset \infty {-}\otimes \overset \infty {-}\bigcirc }$&$\otimes -\otimes \overset \infty {-}\underbrace {\bigcirc -\cdots-\bigcirc \overset 4{-}\bigcirc }_n$\\
\hline
\end{tabular}

\vspace{5pt}

\begin{tabular}{|c|c|c|c|}
\hline
$\fg$&$\fosp(2m+1 | 2n)$; $m,n>1$&$\fosp(4 | 2n)$, $n>1$&$\fosp(2m|2)$, $m>2$\\
\hline
$UW_{\fg}$&$\underbrace {\bigcirc \overset 4{-}\bigcirc -\cdots-\bigcirc }_m \overset \infty {-}\otimes \overset \infty {-}\underbrace {\bigcirc -\cdots-\bigcirc \overset 4{-}\bigcirc }_n$&
$\begin{array}{@{}r@{}l@{}l@{}}
\bigcirc \\[-.16em]
&\hspace*{-.145em}\diagdown \hspace*{-.3em}\raisebox{.3em}[0pt][0pt]{$\scriptstyle \infty $}\\[-.7em]
&&\hspace*{-.4em}\otimes\overset \infty {-}\underbrace {\bigcirc -\cdots-\bigcirc \overset 4{-}\bigcirc }_n \\[-1.3em]
&\hspace*{-.145em}\diagup \hspace*{-.3em}\raisebox{-.3em}[0pt][0pt]{$\scriptstyle \infty $}\\[-.16em]
\bigcirc 
\end{array}
$&$\underbrace {\bigcirc -\cdots-\bigcirc }_{m-2}\overset \infty {-}\otimes \overset \infty {-}\bigcirc $\\
\hline
\end{tabular}

\vspace{5pt}

\begin{tabular}{|c|c|c|c|}
\hline
$\fg$&$\fosp(2m|2n)$; $m>2$, $n>1$&$\fag(2)$&$\fab(3)$\\
\hline
$UW_{\fg}$&$\begin{array}{@{}r@{}c@{}l@{}}
\bigcirc \\[-.16em]
&\hspace*{-.145em}\diagdown \\[-1em]
&&\underbrace {\bigcirc -\cdots-\bigcirc }_{m-2}\overset \infty {-}\otimes \overset \infty {-}\underbrace {\bigcirc -\cdots-\bigcirc \overset 4{-}\bigcirc }_n\\[-1.65em]
&\hspace*{-.145em}\diagup\\[-.16em]
\bigcirc 
\end{array}
$&$\bigcirc \overset \infty {-}\otimes \overset \infty {-}\bigcirc \overset 6{-}\bigcirc $&$\bigcirc \overset \infty {-}\otimes \overset \infty {-}\bigcirc \overset 4{-}\bigcirc -\bigcirc $\\
\hline
\end{tabular}
\end{center}
\end{table} 

\paragraph{The universal super Weyl group} For simplicity, consider the finite-dimen\-sional case.

Let $\fg$ be a~ Lie superalgebra, $\fh$ its maximal torus, $R$ the set of weights in $\fg$ with respect to $\fh$, such weights are called \textit{roots}. For all Lie superalgebras considered, every root $\alpha$ can be represented as ${\alpha=\sum k_i\alpha_i}$, where all coefficients are either non-negative or non-positive for a~ set $B:=\{\alpha_1, \dots, \alpha_{\rk \fg}\}$ called a~\textit{base} or a~ \textit{system of simple roots}. 

Let $\fn_B$ be the nilpotent subalgebra generated by root vectors corresponding to~$B$. A subalgebra of the form $\fb_B:=\fn_B \oplus \fh$ is called a~\textit{Borel subalgebra}. Observe that $\fb_B$ is not the maximal solvable Lie superalgebra, see \cite{Shch, Shch1}; Penkov showed that in the study of irreducible $\fg$-modules the non-maximal subalgebra $\fb_B$ seems  to be more convenient than the maximal one, see \cite{P1, PSe}.

Let $L(\lambda , B)$ be the finite-dimensional irreducible $\fg$-module  with even (of parity 0) $\fb_B$-highest weight vector (i.e., annihilated by $\fb_B$) whose weight is $\lambda $; let $\cP=\cP_B$ be the set of $\fb_B$-highest weights of all the finite-dimensional irreducible $\fg$-modules and $\Gamma _\lambda $ the set of weights of $L(\lambda, B)$.

For any $\alpha \in R$ and $\gamma \in \Gamma _\lambda $, denote by $S_\alpha (\gamma )$ the \textit{$\alpha $-string through $\gamma $}, i.e., the set
\[
\gamma -q \alpha ,\ldots ,\gamma -\alpha ,\gamma ,\gamma +\alpha ,\ldots ,\gamma +p \alpha 
\]
of elements such that $S_\alpha(\gamma)\in \Gamma_\lambda$ and $\gamma -(q+1)\alpha, \ \gamma +(p+1)\alpha \notin \Gamma _\lambda $. The number $l_\alpha =p+q+1$ is called the \textit{length} of the $\alpha $-string.

Set $r_\alpha (\gamma )=\gamma + (p - q)\alpha $. Since the weight $r_\alpha (\gamma )$ is defined for any $\gamma \in \Gamma _\lambda $, there exists a~map $r_\alpha \colon \Gamma _\lambda \mapsto \Gamma _\lambda $. This map will be called the \textit{reflection in $\alpha $}.

Let $F_R$ be the free group with generators $f_\alpha $ for every $\alpha \in R$. Then, for any $\lambda \in \cP$, there is defined an $F_R$-action on $\Gamma _\lambda $ by the formula $f_\alpha (\gamma ) := r_\alpha (\gamma )$.

Let $I_{\lambda ,R}$ be the normal subgroup of $F_R$ singled out by the formula
\[
I_{\lambda ,R}=\{f\in F_R\mid f(\gamma )=\gamma \text{ for all }\gamma \in \Gamma _\lambda \}.
\]

The group $UW_R = F_R/I_R$, where $I_R =\cap _{\lambda \in \cP}I_{\lambda ,R}$, will be called the \textit{universal super Weyl group} or just the universal Weyl group of the root system $R$ (or of $\fg$ and then we denote $UW_R$ by $UW_{\fg}$).  

The group $UW_R$ does not depend on the choice of $B$ for a~cheap reason: the family of \textit{sets} $\{\Gamma _\lambda\}_{\lambda \in \cP_B}$ over which the intersection is taken does not depend on $B$. Indeed, for any base $B$ every finite-dimensional irreducible $\fg$-module has a~$\fb_B$-highest weight vector, so $\{L(\lambda ,B)\mid \lambda \in \cP_B\}$ is, for every $B$, the set of \textit{all} finite-dimensional irreducible $\fg$-modules --- up to isomorphism and the change of parity $\Pi$, which does not change the set of weights. Hence, $\{\Gamma _\lambda\}$ is the set of the weight-supports of all such modules, manifestly $B$-independent, and so is $I_R$. 

Denote by $r_\alpha $ the image of $f_\alpha $ under the natural projection. By construction, $UW_{\fg}$ acts on $\Gamma _\lambda $ for any $\lambda \in \cP$, in particular, it acts on~$R$.

Here, no weight multiplicities are involved: $r_\alpha $ is defined on \textit{weights}, i.e., on elements of $\fh^*$, by the length of the $\alpha $-string through $\gamma $ inside the set $\Gamma _\lambda $; no choice of a~vector in the weight space, and no decomposition of the module into submodules, enters the definition. What \textit{is} used is that $R\cup \{0\}=\Gamma _\lambda $ for the adjoint module; this requires $\fg$ to be simple (for $\fg=\fgl $ one has to pass to the semisimple part). As for the infinite-dimensional case, where there might be no highest weights at all: we will not need $UW$ there --- from \S~2 on its role is played by the concrete permutation groups $S_*$ of eq.~\eqref{chain}. 

\parbegin{Lemma}[$r_\alpha $ is an involution, and $W_{\fg_{\overline 0}}\subset UW_{\fg}$]
\textup{a)} $r_\alpha ^2 = 1$.

\textup{b)} $r_{-\alpha } = r_\alpha $.

\textup{c)} The Weyl group $W_{\fg_{\overline 0}}$ of the Lie algebra $\fg_{\overline 0}$ is naturally embedded into $UW_{\fg}$.

\textup{d)} Let $w \in W_{\fg_{\overline 0}}$. Then, $r_{w(\alpha )}=wr_\alpha w^{-1}$.
\end{Lemma}

\begin{proof} The properties a) and b) directly follow from the definition of the~$r_\alpha $, the property c) is obvious.

Since $r_\alpha $ is a~linear transformation of $\fh^*$ and since it maps $\Gamma _\lambda $ into itself, then
\[
r_{w(\alpha )}(w(\gamma )) = w(r_\alpha (\gamma ))\text{ for any }w\in UW_{\fg},\gamma \in \Gamma _\lambda
\]
implying d).
\end{proof}

\parbegin{Example}[The generators of $UW_{\fg}$ for $\fsl(m|n)$] For $\fg = \fsl(m | n)$, select a~base a.k.a. system of simple roots $B =\{\alpha _1,\ldots ,\alpha _{m+n-1}\}$ corresponding to the Dynkin diagram of the form
\[
\bigcirc -\cdots-\bigcirc-\bigcirc -\otimes -\bigcirc-\bigcirc -\cdots-\bigcirc 
\]
By heading d) of Lemma we can take $r_{\alpha _1},\ldots ,r_{\alpha _{m+n-1}}$ for the generators of $UW_{\fg}$. \end{Example}

\parbegin{Proposition}[$UW_{\fg}$ is a~Coxeter group in two cases, \cite{S2}] $UW_{\fg}$ is a~Coxeter group for $\fg = \fsl(m | n)$ and $\fg = \fosp(2 | 2n)$ with relations indicated in Table $\ref{T1}$.
\end{Proposition}

\parbegin{Conjecture}[$UW_{\fg}$ is a~Coxeter group for every $\fg(A)$] $UW_{\fg}$ is a~Coxeter group for any simple Lie superalgebra with Cartan matrix with relations indicated in Table $\ref{T1}$. \end{Conjecture}

\subsection{The linear Weyl group} 
Under the notations of the previous section, let $A:=(A_{ij})$ be the Cartan matrix (see \cite{CCLL}) and define the reflections $r_\alpha $ by the formulas 
\[
r_{\alpha _i}(\alpha _j)=\begin{cases}
-\alpha _j&\text{for }i = j\\
\alpha _j-A_{ij}\alpha _i&\text{for }i \ne j \text{ and }A_{ii} = 2\\
\alpha _j-2A_{ij}\alpha _i&\text{for }i \ne j \text{ and }A_{ii} = 1\\
\alpha _j+\alpha _i&\text{for }i \ne j \text{ and }A_{ii}=0,\ A_{ji}\ne 0\\
\alpha _j&\text{for }i \ne j \text{ and }A_{ii} = A_{ji} = 0
\end{cases}
\]
Let $LW_R$ be the group generated by such reflections when $\alpha $ runs over simple roots of \textit{all} bases of~ $\fg$. We will call $LW_R$ the \textit{linear Weyl group} of $R$ (or of $\fg$).

\subsection{Formulation of the Problem}\label{ssProb} M. Saveliev was, probably, the first to observe that in applications of Lie superalgebras to integrable dynamical systems, the base of a~simple Lie superalgebra that consists of only odd elements plays a~distinguished role, see~\cite{LSS}. In Table~ 6 of~\cite{LSS}, there are listed Lie superalgebras with which there are associated by a~superization of a~method by Drinfeld and Sokolov, also described in~\cite{LSS}, super versions of KdV and (with infinite diagrams) KP, see a~pioneering paper~\cite{MR} and later works~\cite{KL1,KL2,R}. Note that only a~\textit{method} to recover an equation from a~superalgebra is written in~\cite{LSS}\textup{;} the method itself was carried out in detail later, see~\cite{IK, DG}, but the exceptional super KdVs corresponding to the exceptional diagrams of Table~ 6 of~\cite{LSS} were never written out explicitly, and are not to be found in~\cite{IK, DG} either. 

Not every Lie superalgebra has such a~base. Serganova's result states that if the Lie superalgebra is finite-dimensional or loop one (perhaps, twisted, or is its double extension with Cartan matrix, see \cite{BLS}) it does not matter which base we start with as long as the superalgebra has a~distinguished base since we can reach any of the bases from any given one with the help of any of Serganova's super Weyl groups.

\textbf{Question: Is it true for $\fgl _c(\infty |\infty )$ that there is just one class of bases with respect to an analogue of the Weyl group? } Since only $\fgl _c(\infty |\infty )$ was chosen so far (\cite{R,KL1, KL2,B}) to study superized versions of KdV and KP, the answer is important for these interesting differential equations.

\textbf{It turns out that the answer is ``NO\rq\rq}.

\subsection{Formulation of the answer} For a~complete and explicit description of bases (systems of simple roots) for the so far conventional $\fgl_c(\infty |\infty)$, see Theorem \ref{Th}.  

Consider also interesting enlargements of $\fgl _c(\infty |\infty )$. For the largest of these enlargements, there is just one class of bases \textit{and} for all of these types of Lie superalgebras, there exists a~ nontrivial cocycle of degree~0 which defines a nontrivial central extension.

For analogues of the Weyl group for these enlarged $\fgl (\infty |\infty )$'s, see Subsec. \ref{ssWeyl}, where the six permutation groups $S_g$, $S_l$, $S_m$, $S_n$, $S_o$ and $S_c$ that play this role are defined, and Theorem \ref{Thm}, which says which of them governs which enlargement.

For a~similar description for $\fgl_c(p|\infty )$ and $\fgl _c(\infty |q)$ with ${p,q < \infty }$, see \ref{Sclasses}, where these are the superalgebras attached to the \textit{finite} parity functions and are denoted $\fgl _{p,*}(m|\infty )$ and $\fgl _{p,*}(\infty |n)$: by Proposition \ref{p3.1}, the two cardinalities $\Even(p;-\infty ,+\infty )$ and $\Odd(p;-\infty ,+\infty )$ are a~complete system of invariants of such a~function with respect to each of $S_c$, $S_o$, $S_l$ and $S_g$ at once, so here the four enlargements are not distinguished.

In \ref{Smisc},  various important subalgebras of $\fgl _c(\infty |\infty )$ are described: its orthosymplectic, periplectic and queer versions, and how to embed certain affine Kac-Moody (see \cite{BLS}) superalgebras into some of these Lie superalgebras of infinite matrices.

\section{Examples of infinite matrix Lie superalgebras and their ``Weyl groups\rq\rq\ }

In this section, we introduce several (super)spaces of infinite matrices. Let us start with the underlying spaces. Let ${\Mat: =\{a = (a_{i,j})_{i,j\in \Zee}\}}$ designate the space of all infinite matrices. Set
\begin{gather*}
\Mat_f := \{x \in \Mat \mid \card (\Supp x) <\infty \},\text{ where}\\
 \Supp (a) := \{(i,j) \in \Zee \times \Zee\mid a_{i,j}\ne 0\}.
\end{gather*}

In what follows, an arbitrary map $p \colon \Zee \tto \Zee/2$ will be called a~\textit{parity function}. The following two parity functions are very important for us:
\begin{gather*}
p_{\text{st}}(x)=x \mod 2, \ \ \ \ \ \ p_{+}(x)=\begin{cases}
\overline 1&\text{if } x \geqslant 0,\\
\overline 0&\text{otherwise.}
\end{cases}
\end{gather*}

On $\Mat_f$, introduce a~Lie superalgebra structure by defining the supermatrices with the help of a~fixed parity function $p$ and by setting
\be\label{1}
[a,b]_{ij}=\sum\limits _k(a_{ik}b_{kj}-(-1)^{(p(i)+p(k))(p(j)+p(k))}b_{ik}a_{kj});
\ee
denote the Lie superalgebra obtained by $\fgl _{p,f}$. The Dynkin diagrams corresponding to $p_{+}$ and $p_{\text{st}}$, respectively, are:
\[
\cdots-\bigcirc -\cdots-\bigcirc -\bigcirc-\otimes -\bigcirc-\bigcirc -\cdots-\bigcirc -\cdots
\]
and
\[
\cdots-\otimes -\cdots-\otimes -\otimes -\cdots-\otimes -\otimes -\cdots-\otimes -\cdots
\]

\subsection{Several superized versions of $\fgl (\infty )$} We will consider the following four en\-lar\-ge\-ments of $\fgl_{p,f}$ with the bracket defined by the formula 
\be\label{(2)}
\begin{array}{l}
\fgl _{p,g}:=\{(a_{ij})_{i,j\in \Zee} \in \Mat\mid \text{for any $(i_0, j_0)$ there exists}\\
\text{some $c := c(i_0,j_0) > 0$ such that}\\
a_{ij}=a_{ji}=0 \text{ for all $i$, $j$ such that }(i - i_0)(j_0 - j) > c\}.
\end{array}
\ee
In other words, the matrix $A$ belongs to $\fgl _{p,g}$ if for any point $(i_0,j_0)$ the following holds: ``Translate the origin to this point $(i_0,j_0)$; then, in the new 2nd and 4th quadrants there are finitely many  points of the support of $A$\rq\rq, see Fig.~\ref{Fg}. 

\begin{figure}[ht]
\centering
\begin{tikzpicture}[scale=0.60]
\begin{scope}
\clip (-5.2,-5.2) rectangle (5.2,5.2);
\fill[black!12] plot[domain=-5.2:0.75,samples=70] (\x,{1.6+1.1/(1.2-\x)})
   -- (-5.2,5.2) -- cycle;
\fill[black!12] plot[domain=5.2:1.65,samples=70] (\x,{1.6-1.1/(\x-1.2)})
   -- (5.2,-5.2) -- cycle;
\end{scope}
\draw[->] (-5.4,0)--(5.9,0) node[below]{$i$};
\draw[->] (0,-5.4)--(0,5.9) node[left]{$j$};
\draw[dashed] (-5.2,-5.2)--(5.2,5.2) node[above right]{\small $j=i$};
\draw[dotted,thick] (-5.2,1.6)--(5.2,1.6);
\draw[dotted,thick] (1.2,-5.2)--(1.2,5.2);
\fill (1.2,1.6) circle (2.2pt);
\node[above left] at (1.3,1.6) {\small $(i_0,j_0)$};
\node at (-3.9,4.3) {\small $\mathrm{II}$};
\node at (4.1,-3.3) {\small $\mathrm{IV}$};
\node[below] at (2.7,-5.0) {\small $(i-i_0)(j_0-j)=c$};
\foreach \p/\q in {-4.7/-4.4, -4.2/-4.7, -3.6/-3.3, -3.1/-3.7, -2.5/-2.1,
  -1.9/-2.4, -1.2/-0.9, -0.6/-1.1, 0.1/0.4, 0.6/0.1, 1.5/1.9, 2.0/1.5,
  2.6/3.0, 3.1/2.5, 3.7/4.0, 4.2/3.6, 4.8/5.0, -3.4/0.7, 0.5/-2.7, 1.5/-1.5,
  -1.0/1.9, 3.4/1.3}
  {\fill (\p,\q) circle (2.2pt);}
\end{tikzpicture}
\caption{The condition \eqref{(2)} defining $\fgl _{p,g}$, drawn in the plane
of the indices, the dots being a~sample of $\Supp a$. Translate the origin to
any point $(i_0,j_0)$: in the two shaded regions of the new frame --- the parts
of its 2nd and 4th quadrants lying beyond the hyperbola
$(i-i_0)(j_0-j)=c(i_0,j_0)$ --- the matrix $a$ vanishes identically, and in the
2nd and the 4th quadrants together it has but finitely many nonzero entries.
The support is thus pressed against the main diagonal, yet not into a~band of
any prescribed width: what the condition constrains is quadrants, not
distances $|i-j|$.}
\label{Fg}
\end{figure}

It is an easy exercise for the reader to prove that there are finitely many entries in every row and every column as well as in the second and the  fourth quadrants of every matrix in~ $\fgl _{p,g}$. It is, however, not immediately clear why this Lie superalgebra is closed under the superbracket. 

\ssbegin{Lemma}[The largest of the four is closed under the bracket]\label{ssLem} $\fgl_{p,g}$ is closed under the superbracket \eqref{1}.
\end{Lemma}

\begin{proof}(I.Shchepochkina) 
For the elements $a,b$ the condition $a, b\in \fgl_{p,g}$ for $(i_0,j_0)=(0,0)$ means that 
\[
\text{$a_{ij}=a_{ji}=b_{ij} = b_{ji}=0 $ for all $i,j$ such that $ij<-c$ for some $c>0$.}
\]
In particular, in the 2nd quadrant the set where $a_{ij}=a_{ji}=0 $ contains the sets
 $\{(i, j)\mid i<-c, j>0\}$ and  $\{(i, j)\mid i<0, j>c\}$.

Therefore, if $i<-c$ and $j>0$, then in formula \eqref{(2)} we have
$$
a_{ik}=b_{ik}=0 \text{  for   } k>0 \text{  and     } b_{kj}=a_{kj}=0 \text{  for   } k<-c.
$$
Thus, in this case,
\be\label{(3)}
[a,b]_{ij}=\sum_{-c\leq k\leq 0}(a_{ik}b_{kj}-(-1)^{...}b_{ik}a_{kj}). 
\ee

The condition $a,b\in \fgl_{p,g}$ for $(i_0,j_0)=(0,-c-1)$ is as follows:
\begin{equation}\label{cond1}
\begin{minipage}[l]{13.3cm}
There exists $c_1>0$ such that $a_{ij}=b_{ij} =0$ for all $i,j$ such that $i(j+c+1)<-c_1$, in particular, for any $i,j$ such that $i<-c_1$ and $j+c+1>0\Llra j>-c-1$. \end{minipage}
\end{equation}
Therefore, if $c_2:=\max(c,c_1)$ and $i<-c_2$ whereas $j>0$, then in formula \eqref{(3)} we see that
$$
a_{ik}=b_{ik}=0 \text{  for all     } k=-c,\dots, 0 \Lra [a,b]_{ij}=0.
$$

Now, let $j>c, i<0$. In this case,
$$
a_{kj}=b_{kj}=0 \text{  for   } k<0 \text{  and     } b_{ik}=a_{ik}=0 \text{  for   } k>c.
$$
Hence, in this case 
\be\label{(4)}
[a,b]_{ij}=\sum_{0\leq k\leq c}(a_{ik}b_{kj}-(-1)^{...}b_{ik}a_{kj}). 
\ee

The condition $a,b\in \fgl_{p,g}$ for  $(i_0,j_0)=(c+1,0)$ is as follows:
\begin{equation}\label{cond2}
\begin{minipage}[l]{13.3cm}
There exists $c_3>0$ such that $a_{ij}=0$ for all $i,j$ such that $(i-c-1)j<-c_3$, in particular, for any $i,j$ such that $i-c-1<0\Llra i<c+1$ and $j>c_3$. \end{minipage}
\end{equation}

Therefore, if $c_4:=\max(c,c_3)$ and     $i<0$, $j>c$, then in formula \eqref{(4)} we see that
$$
a_{kj}=b_{kj}=0 \text{  for all     } k=0,\dots, c \Lra [a,b]_{ij}=0.
$$

Thus, we proved that the intersection of the support of $[a,b]$ with the 2nd quadrant is contained inside the square $\{i, j\mid -d_1\le i\le 0, \ 0\le j\le d_1\}$, where  $d_1=\max(c_2,c_4)$.

We can similarly prove that there exists a~ $d_2>0$ such that the intersection of the support of $[a,b]$ with the 4th quadrant is contained inside the square 
\[
\{i, j\mid 0\le i\le d_2, \ -d_2\le j\le 0\}.
\]

Let now $C:=\max(d_1^2,d_2^2)>0$, then
$$
[a,b]_{ij}=[a,b]_{ji}=0 \text{  for all     } ij \text{ such that } ij<-C.
$$
This precisely means that the condition on $a,b\in \fgl_{p,g}$ holds for $(i_0,j_0)=(0,0)$. \end{proof}

Introduce the Lie superalgebras
\begin{gather*}
\fgl _{p,l}:=\{a=(a_{ij})_{i,j\in \Zee} \in \Mat\mid \lim\limits _{\substack{(i,j) \in \Supp a\\
|i| \tto \infty }}\frac {|i-j|^\lambda }{|i+j|}=\lim\limits _{\substack{(i,j) \in \Supp a\\
|j| \tto \infty }}\frac {|i-j|^\lambda }{|i+j|}=0\\
\text{for some }\lambda =\lambda (a)\geq 1\};\\
\fgl _{p,o}:=\{a=(a_{ij})_{i,j\in \Zee} \in \Mat\mid \lim\limits _{\substack{(i,j) \in \Supp a\\
|i| \tto \infty }}\frac {|i-j|^\lambda }{|i+j|}=\lim\limits _{\substack{(i,j) \in \Supp a\\
|j| \tto \infty }}\frac {|i-j|^\lambda }{|i+j|}=0 \\
 \text{for all}\quad \lambda \geq 1\};\\
\fgl _{p,c}:=\{a=(a_{ij})_{i,j\in \Zee}\in \Mat\mid a_{ij}=0 \text{ for all $i$, $j$ such that}\\
\text{$|i -j| > c$ for some constant $c = c(a)$}\}.
\end{gather*}

\sssbegin{Remark}[Why $\lambda \geq 1$ is required]\label{lam1} The restriction $\lambda \geq 1$ in the definition of $\fgl _{p,l}$ is essential, and it is also what makes this definition agree with the definition of the group $S_l$ below, where $\lambda \geq 1$ is required as well. 

For $\fgl _{p,o}$, on the contrary, the printed ``for all $\lambda \geq 1$\rq\rq\  and ``for all $\lambda >0$\rq\rq\  are equivalent, since for $|i-j|\geq 1$ the ratio $\frac{|i-j|^\lambda }{|i+j|}$ only decreases when $\lambda$ does.
\end{Remark}

\sssbegin{Proposition}[The four are Lie superalgebras and they form a~chain]\label{Pclosed} For $* = g,l,o,c$, the spaces $\fgl _{p,*}$ are Lie superalgebras, i.e., for any $a,b\in \fgl _{p,*}$ we have $[a,b] \in \fgl _{p,*}$. Moreover,
\be\label{glchain}
\fgl _{p,c}\subset \fgl _{p,o}\subset \fgl _{p,l}\subset \fgl _{p,g}.
\ee
\end{Proposition}

\begin{figure}[ht]
\centering
\begin{tikzpicture}[scale=0.60]
\begin{scope}
\clip (-0.7,-0.7) rectangle (9.9,9.9);
\fill[black!10] plot[domain=0:9.9,samples=70] (\x,{\x+0.4+0.75*sqrt(\x)})
   -- plot[domain=9.9:0,samples=70] (\x,{\x-0.4-0.75*sqrt(\x)}) -- cycle;
\fill[black!24] plot[domain=0:9.9,samples=70] (\x,{\x+0.4+0.38*ln(1+\x)})
   -- plot[domain=9.9:0,samples=70] (\x,{\x-0.4-0.38*ln(1+\x)}) -- cycle;
\fill[black!45] (0,0.4) -- (9.9,10.3) -- (9.9,9.5) -- (0,-0.4) -- cycle;
\draw[dotted,thick] (0,0) -- (9.9,{9.9*0.6});
\end{scope}
\draw[->] (-0.7,0)--(10.7,0) node[below]{$i$};
\draw[->] (0,-0.7)--(0,10.7) node[left]{$j$};
\draw[dashed] (-0.5,-0.5)--(9.9,9.9);
\node[right] at (9.7,10.1) {\small $j=i$};
\draw[thin] (9.9,9.1)--(10.5,9.1) node[right]{\small $\fgl _{p,c}$};
\draw[thin] (9.9,8.0)--(10.5,8.0) node[right]{\small $\fgl _{p,o}$};
\draw[thin] (9.9,7.1)--(10.5,7.1) node[right]{\small $\fgl _{p,l}$};
\draw[thin] (9.9,5.9)--(10.5,5.9) node[right]{\small $|i-j|=\eta |i|$};
\end{tikzpicture}
\caption{The chain \eqref{glchain}. Only the quadrant $i,j>0$ is drawn. The
quadrant $i,j<0$ carries the same picture rotated by $\pi $, and nothing is
lost by omitting the two remaining quadrants: in them every element of
$\fgl _{p,l}$ has but finitely many nonzero entries, see heading $2)$ of
Lemma \ref{Llam}. Only the behavior of the boundaries as
$i \tto \infty $ carries a~meaning: $\fgl _{p,c}$ is a~band $|i-j|\leq c$ of
a~fixed width, $\fgl _{p,o}$ a~band flaring slower than any power $|i|^\eps $,
and $\fgl _{p,l}$ a~band flaring slower than $|i|$ itself, that is, one that
enters every angle $|i-j|\leq \eta |i|$, the dotted line being the boundary of
one such angle. The three curves drawn are but one representative of each
class. The finite matrices $\fgl _{p,f}$ occupy a~bounded piece of the picture
near the origin. The largest member of the chain, $\fgl _{p,g}$, has no region
in this picture, and that is not an omission: the region its condition forbids,
see Fig.~\ref{Fg}, moves together with the point $(i_0,j_0)$ from which one
looks, so no curve in the plane of the indices bounds it. What the last
inclusion of \eqref{glchain} says is that a~matrix whose support lies in any
one of the four regions drawn here does obey that condition at every
$(i_0,j_0)$: in the 2nd and the 4th quadrant of every frame so translated it
has but finitely many nonzero entries.}
\label{Fchain}
\end{figure}

For $*=g$ this is precisely the claim of Lemma \ref{ssLem}, and for $*=c$ it is obvious. 
For $*=l,o$, denote by $(\mathrm{L}_\lambda )$ the condition
\be\label{L-lambda}
(\mathrm{L}_\lambda ):\qquad \lim\limits _{\substack{(i,j) \in \Supp a\\ |i| \tto \infty }}\frac {|i-j|^\lambda }{|i+j|}=
\lim\limits _{\substack{(i,j) \in \Supp a\\ |j| \tto \infty }}\frac {|i-j|^\lambda }{|i+j|}=0,
\ee
so that $\fgl _{p,l}=\{a\mid (\mathrm{L}_\lambda )$ holds for some $\lambda \geq 1\}$ and $\fgl _{p,o}=\{a\mid (\mathrm{L}_\lambda )$ holds for all $\lambda \geq 1\}$.

\sssbegin{Lemma}[What $(\mathrm{L}_\lambda )$ says about the support]\label{Llam} Let $a$ satisfy $(\mathrm{L}_\lambda )$ for some $\lambda \geq 1$. Then,

\textup{1)} every row and every column of $a$ is finite\textup{;}

\textup{2)} the set $\{(i,j)\in \Supp a\mid ij<0\}$ is finite\textup{;}

\textup{3)} $|i-j|=o(|i|)$ and $|i+j|\sim 2|i|$ on $\Supp a$ as $|i|\tto \infty $\textup{;}

\textup{4)} $a\in \fgl _{p,g}$.
\end{Lemma}

\begin{proof} 1) If the $i$-th row were infinite, there would be points $(i,j)\in \Supp a$ with $|j| \tto \infty $; along them $\frac{|i-j|}{|i+j|} \tto 1$ and $|i-j| \tto \infty $, hence
\[
\frac {|i-j|^\lambda }{|i+j|}=|i-j|^{\lambda -1}\cdot \frac {|i-j|}{|i+j|}\not\tto 0\quad \text{for }\lambda \geq 1.
\]
For columns, the same argument with the first of the two limits.

2) If $ij<0$, then $|i-j|=|i|+|j|\geq |i+j|$ and $|i-j|\geq 1$, so $\frac{|i-j|^\lambda}{|i+j|}\geq |i-j|^{\lambda -1}\geq 1$. Hence, no such point has $|i|>N$, where $N$ corresponds to $\varepsilon =\frac12$ in the first limit; the points with $|i|\leq N$ lie in finitely many rows, each finite by heading~1).

3) For $\lambda =1$ this is the first limit itself; for $\lambda >1$ the relation $\frac{|i-j|^\lambda}{|i+j|} \tto 0$ gives 
\[
|i-j|=o(|i+j|^{1/\lambda })=o(|i+j|), 
\]
and by heading~2) the signs of $i$ and $j$ coincide for $|i|$ large, so $|i+j|=|i|+|j|\leq 2|i|+|i-j|$, whence $|i+j|\sim 2|i|$.

4) Fix $(i_0,j_0)$ and let $N$ be so large that for $(i,j)\in \Supp a$ with $|i|>N$ the signs of $i$ and $j$ coincide and $|i-j|<\frac12|i|$, see headings 2) and 3); enlarging $N$ we may assume that $N>2(|i_0|+|j_0|)$. For such $(i,j)$, both $i$ and $j$ are of the same sign and of absolute value $>|i_0|+|j_0|$; hence $i-i_0$ and $j-j_0$ are of the same sign and $(i-i_0)(j-j_0)>0$. The remaining points of $\Supp a$ lie in the finitely many rows with $|i|\leq N$, each of them finite; let $c$ exceed $|(i-i_0)(j-j_0)|$ for all of them. This is precisely the condition \eqref{(2)} defining $\fgl _{p,g}$.
\end{proof}

\begin{proof}[Proof of Proposition \ref{Pclosed} for $*=l,o$] Let $a$ satisfy $(\mathrm{L}_{\lambda _a})$ and $b$ satisfy $(\mathrm{L}_{\lambda _b})$, where $\lambda _a,\lambda _b\geq 1$.

First, the sum $\sum_k(a_{ik}b_{kj}\pm b_{ik}a_{kj})$ in the bracket \eqref{1} is finite for every $(i,j)$ because the $i$-th rows of $a$ and of $b$ are finite, see heading~1) of Lemma \ref{Llam}. So $[a,b]$ is well defined.

Next, let $[a,b]_{ij}\neq 0$; then $(i,k)\in \Supp a$ and $(k,j)\in \Supp b$ for some $k$ (or the same with $a$ and $b$ interchanged, which is the same argument). Let $|i| \tto \infty $. By heading~3) of Lemma \ref{Llam} applied to $a$ we have $|i-k|=o(|i|)$, hence $|k|\sim |i|$ and $|i-k|=o(|i+k|^{1/\lambda _a})=o(|i|^{1/\lambda _a})$; applied to $b$ it gives $|k-j|=o(|k|^{1/\lambda _b})=o(|i|^{1/\lambda _b})$. Therefore,
\[
|i-j|\leq |i-k|+|k-j|=o\big(|i|^{1/\lambda }\big),\quad \text{where}\quad \lambda :=\min(\lambda _a,\lambda _b),
\]
and, since $i$ and $j$ are of the same sign for $|i|$ large, $|i+j|\sim 2|i|$, so that
\[
\frac {|i-j|^{\lambda }}{|i+j|}=\frac {o(|i|)}{|i+j|} \tto 0 .
\]
The limit as $|j| \tto \infty $ is the same limit, because $|i|\sim |j|$ on $\Supp [a,b]$. Thus, $[a,b]$ satisfies $(\mathrm{L}_\lambda )$ for $\lambda =\min(\lambda _a,\lambda _b)\geq 1$, and $\fgl _{p,l}$ is closed under the bracket; for $\fgl _{p,o}$ take $\lambda _a=\lambda _b$ arbitrarily large.

Finally, the inclusions \eqref{glchain}: $\fgl _{p,c}\subset \fgl _{p,o}$ is clear (a matrix supported in a~bounded band satisfies $(\mathrm{L}_\lambda )$ for every $\lambda $), $\fgl _{p,o}\subset \fgl _{p,l}$ is the case $\lambda =1$, and $\fgl _{p,l}\subset \fgl _{p,g}$ is heading~4) of Lemma \ref{Llam}.
\end{proof}

The inclusions \eqref{glchain}, drawn in Fig.~\ref{Fchain}, are what heading 2) of Conjecture \ref{ssConj} below tacitly uses when it speaks of ``the restrictions of this cocycle\rq\rq.

The superalgebra $\fgl _{p,g}$ has a~nontrivial central extension $\widetilde {\fgl }_{p,g}$
given by the cocycle \eqref{tr}, with $\str$ instead of $\tr$. Let us make this precise: both the well-definedness of the cocycle and the cocycle identity become evident once it is rewritten in a~closed form. Write any~matrix $X$ in the block form \eqref{corners} attached to the splitting $V=V_-\oplus V_+$ that $J=\sum_{i\geq 0}E_{ii}$ defines: the summand $V_+$ is the image of $J$ and $V_-$ is its kernel, and the four blocks of $X$ are
\[
JXJ,\qquad JX(\One -J),\qquad (\One -J)XJ,\qquad (\One -J)X(\One -J).
\]
Since $J[X,Y]J=[X,Y]_{++}$ and $[JXJ,JYJ]=[X_{++},Y_{++}]$, all the ``$++$\rq\rq\  terms cancel and, for homogeneous $X$ and $Y$,
\be\label{coc2}
c(X,Y):=\str\big(J[X,Y]J-[JXJ,JYJ]\big)=\str\big(X_{+-}Y_{-+}\big)-(-1)^{p(X)p(Y)}\str\big(Y_{+-}X_{-+}\big).
\ee
Now, for $X\in \fgl _{p,g}$ the corner block $X_{+-}$ (and likewise $X_{-+}$) has \textit{finitely many} non-zero entries: by \eqref{(2)} with $(i_0,j_0)=(0,0)$ we have $a_{ij}=0$ for $ij<-c$, i.e., for $i\geq 1$, $j\leq -1$ only the finitely many pairs with $i|j|\leq c$ survive; the row $i=0$ is cut off by \eqref{(2)} with $(i_0,j_0)=(-1,0)$, and the column $j=0$ by $(i_0,j_0)=(0,1)$. Hence, \eqref{coc2} is a~finite sum; it is manifestly bilinear and anti-symmetric, i.e., $c(X,Y)=-(-1)^{p(X)p(Y)}c(Y,X)$, and the cocycle identity is the standard computation for $\str (X_{+-}Y_{-+})$.

\sssbegin{Proposition}[$\fgl _{p,g}$ is exactly where the cocycle
survives]\label{Pcorner} For $a\in \Mat $, the following conditions are
equivalent\textup{:}

\textup{1)} $a\in \fgl _{p,g}$, i.e., $a$ satisfies \eqref{(2)}\textup{;}

\textup{2)} for every $n\in \Zee $, the two corner blocks of $a$ determined by the
polarization $V=V_{<n}\oplus V_{\geq n}$, where
\[
V_{<n}:=\Span (e_i\mid i<n)\quad \text{and}\quad V_{\geq n}:=\Span (e_i\mid i\geq n),
\]
i.e., the sets
\[
\Supp a\cap \{(i,j)\mid i\geq n>j\}\quad \text{and}\quad
\Supp a\cap \{(i,j)\mid i<n\leq j\},
\]
are finite.

Hence, $\fgl _{p,g}$ is the largest subspace of $\Mat $ on which \eqref{coc2}, and
likewise its analog attached to any shifted polarization, is given by a~finite
sum.
\end{Proposition}

\begin{proof} $1)\Longrightarrow 2)$: the corner block $\{i\geq n>j\}$ lies in the
4th quadrant seen from $(n-1,n)$, and the other one in the 2nd quadrant seen from
$(n,n-1)$.

$2)\Longrightarrow 1)$: first, $2)$ makes every row and every column of $a$
finite. Indeed, the part of the row $i$ with $j<i$ lies in the corner block at
$n=i$, and the part with $j>i$ in the corner block at $n=i+1$; columns are treated
the same way. Now, fix $(i_0,j_0)$ and look at the 4th quadrant seen from it,
i.e., at $\{(i,j)\mid i>i_0,\ j<j_0\}$. If $j_0\leq i_0+1$, this set lies inside
the corner block at $n=i_0+1$. If $j_0>i_0+1$, it is the union of the corner block
at $n=j_0$ and of the $j_0-i_0-1$ rows with $i_0<i<j_0$, each of them finite by the
above. Either way it meets $\Supp a$ in a~finite set, and the 2nd quadrant is
treated the same way. Thus, \eqref{(2)} holds with $c$ the largest of the finitely
many values of $|(i-i_0)(j-j_0)|$ so obtained.

The last claim now only restates $2)$: the sum defining \eqref{coc2} for the
polarization at $n$ runs over the two corner blocks at $n$, and it is finite for
all $n$ exactly when $2)$ holds. \end{proof}

\sssbegin{Corollary}[The superalgebras \eqref{glchain} lie in the Mackey Lie algebra, and
properly]\label{CMackey} Every matrix of $\fgl _{p,g}$, hence of each of the four
superalgebras \eqref{glchain}, has finitely many non-zero entries in every row and
in every column\textup{;} that is, the chain \eqref{glchain} lies inside the
Mackey Lie algebra $\fgl ^M(V,V_*)$ of \cite{PS1,ChP1}, i.e., in the algebra of those $X\in \fgl (V)$ whose dual operator preserves the restricted dual $V_*:=\bigoplus_{i\in \Zee}\Cee e_i^*$ of the pair $(V,V_*)$\textup{;} in matrix terms, this is precisely row- and column-finiteness. The inclusion is
proper, and the cocycle \eqref{coc2} does not extend to $\fgl ^M(V,V_*)$.
\end{Corollary}

\begin{proof} The first claim is the row-and-column finiteness established in the
proof of Proposition \ref{Pcorner}. For the second one, take the antidiagonal
matrix $a:=\sum_{i\in \Zee}E_{i,-i}$. It has exactly one non-zero entry in every
row and in every column, so $a\in \fgl ^M(V,V_*)$; but every point of $\Supp a$
lies in the 2nd or in the 4th quadrant seen from the origin, so $a\notin
\fgl _{p,g}$. Moreover, $a_{+-}$ consists of the entries $(i,-i)$ with $i\geq 1$,
and $a_{-+}$ of the entries $(i,-i)$ with $i\leq -1$, so already
$\str (a_{+-}a_{-+})=\sum_{1\leq i}(-1)^{p(i)}$ has no value: for $p\equiv \ev $
it is $\sum_{1\leq i}1$. \end{proof}

\sssbegin{Lemma}[$\Mat_f$ is an ideal of $\fgl _{p,*}$, associatively and in the Lie sense]\label{Lideal} Let
$*=g,l,o,c$.

\textup{1)} Each $\fgl _{p,*}$ is closed under the matrix product\textup{:} for
$a,b\in \fgl _{p,*}$, the product $ab$, well defined by Corollary
$\ref{CMackey}$, lies in $\fgl _{p,*}$. Thus, $\fgl _{p,*}$ is an associative
superalgebra, and the bracket \eqref{1} is its supercommutator.

\textup{2)} $\Mat_f$ is a~two-sided ideal of the \textit{associative} superalgebra $\fgl _{p,*}$ of heading $1)$\textup{:} we have $au$, $ua\in \Mat_f$ for all $a\in \fgl _{p,*}$
and all $u\in \Mat_f$. More generally, this holds for every row- and
column-finite $a$, i.e., for every $a$ in the Mackey Lie algebra.

\textup{3)} Hence, $\Mat_f$ is an ideal of the \textit{Lie} superalgebra
$\fgl _{p,*}$, i.e., $[\fgl _{p,*},\Mat_f]\subseteq \Mat_f$.
\end{Lemma}

\begin{proof} 1) For $*=c$ the claim is clear, the band widths adding up. For
$*=g$, observe that the condition \eqref{(2)} on $ab$ at a~point $(i_0,j_0)$ is
the condition at the origin on another pair of matrices of $\fgl _{p,g}$.
Indeed, put $a^{(u,v)}_{ij}:=a_{i+u,\,j+v}$; then
\[
(ab)_{i+i_0,\,j+j_0}=\sum_k a_{i+i_0,\,k}\,b_{k,\,j+j_0}=
\big(a^{(i_0,0)}b^{(0,j_0)}\big)_{ij},
\]
and $a^{(u,v)}\in \fgl _{p,g}$ whenever $a\in \fgl _{p,g}$, because the
condition \eqref{(2)} on $a^{(u,v)}$ at $(i_0,j_0)$ is the condition on $a$ at
$(i_0+u,j_0+v)$, the family of conditions \eqref{(2)} being merely reindexed. So
it suffices to bound $\Supp (ab)$ in the 2nd and the 4th quadrant seen from the
origin --- and that is exactly what the proof of Lemma \ref{ssLem} does: it
estimates the two products $a_{ik}b_{kj}$ and $b_{ik}a_{kj}$ of \eqref{1}
separately, and never uses a~cancellation between them. For $*=l,o$, the same
remark applies to the proof of Proposition \ref{Pclosed} in these cases: it
argues from ``$(i,k)\in \Supp a$ and $(k,j)\in \Supp b$ for some $k$\rq\rq, which
is precisely the condition for $(i,j)$ to lie in $\Supp (ab)$, and it delivers
$(\mathrm{L}_\lambda )$ for $ab$ with $\lambda =\min(\lambda _a,\lambda _b)$.

2) The $j$-th column of $au$ is the image under $a$ of the $j$-th column of $u$,
hence a~finite linear combination of columns of $a$, each of them finite by
Corollary \ref{CMackey}; and it vanishes unless $j$ is one of the finitely many
indices of the non-zero columns of $u$. So $\Supp (au)$ is finite. For $ua$,
argue with rows in place of columns.

3) This is heading 2), the bracket \eqref{1} being the supercommutator of the
product. \end{proof}

\sssbegin{Remark}[An ideal of \textit{what}]\label{Rideal} Being an ideal is a
property of a~subspace with respect to a~structure on the ambient space, not a
property of an algebra by itself, and the space $\Mat_f$ carries two structures
at once. As an \textit{associative} superalgebra it is simple, see Subsec.
\ref{ssWeyl}; as a~subspace of $\fgl _{p,*}$ it is an ideal with respect to
\textit{both} structures, by headings 2) and 3) of Lemma \ref{Lideal}. The two
statements do not collide --- they are about different structures --- and the
implication between the two readings of ``ideal\rq\rq\  goes one way only: a
two-sided associative ideal is a~Lie ideal, and not conversely, the scalar
matrices being a~Lie ideal of $\fgl (n)$, its center, and as far from an
associative ideal as a~nonzero subspace can be. Below, the word ``ideal\rq\rq\
applied to $\Mat_f$ always means the Lie one, which is all that the argument of
Subsec. \ref{ssWeyl} uses; the associative structure of $\Mat_f$ enters only
through clause $2)$ of Convention \ref{2.3.1} and Lemma \ref{NS}.
\end{Remark}

\sssbegin{Remark}[Condition \eqref{(2)} says that all the corner blocks are finite]\label{Rpolar} Proposition
\ref{Pcorner} restates the condition \eqref{(2)}, which Serganova wrote down in
terms of quadrants, without any reference to the plane of the indices: a~matrix
lies in $\fgl _{p,g}$ if and only if its two corner blocks are of finite rank for
\textit{every} polarization of $V$ by the order of the basis, and not only for the
one that defines the cocycle. In this form, $\fgl _{p,g}$ is the exact algebraic
counterpart of the algebra $\fgl _{\text{res}}$ of \cite{SW}. Here, the
\textit{corner blocks} of a~matrix $X$ relative to a~polarization
$V=V_+\oplus V_-$ are its two off-diagonal blocks $X_{+-}$ and $X_{-+}$, the
ones that enter the cocycle \eqref{coc2}. The name is ours: \cite{SW} writes an
operator of a~polarized Hilbert space as
$\left (\begin{smallmatrix}a&b\\ c&d\end{smallmatrix}\right )$ and cuts out
$\fgl _{\text{res}}$ by asking that these very two blocks, $b$ and $c$, be
Hilbert--Schmidt for the given polarization; and it is the
form in which the answer to the question ``how far may one enlarge?\rq\rq\  of
Subsec. \ref{ssKP} of the Introduction is visible.
\end{Remark}

\sssbegin{Conjecture}[The central extension is unique on each of the four superalgebras]\label{ssConj} \textup{1)} The nontrivial central extension $\widetilde {\fgl }_{p,g}$ is unique.

\textup{2)} The restrictions of this cocycle to $\fgl _{p,l}$, $\fgl _{p,c}$ and $\fgl _{p,o}$ are nontrivial and unique.
\end{Conjecture}

The three central extensions that heading $2)$ of Conjecture \ref{ssConj} speaks of do exist, as Theorem \ref{Tnontriv} below proves\textup{;} we denote them by $\widetilde {\fgl }_{p,l}$, $\widetilde {\fgl }_{p,c}$ and $\widetilde {\fgl }_{p,o}$, respectively.

\textbf{Comment}. The existence of the cocycle in heading 1) Conjecture \ref{ssConj} is obvious, see \eqref{coc2}; it is uniqueness which has to be established. For $p\equiv 0$ this is the classical computation of $H^2(\fgl _c(\infty ))$, see \cite{FF}. The rest of this Subsection proves what we can prove of Conjecture \ref{ssConj}: its nontriviality half for all the four superalgebras \eqref{glchain}, and its uniqueness half for $*=c$, the case whose non-super prototype is the computation just quoted. Remark \ref{Rcont} says in which category the word ``unique\rq\rq\  is to be understood, and Remark \ref{Rleft} what the remaining three cases need.

\sssbegin{Theorem}[All the four extensions are nontrivial]\label{Tnontriv} Whatever the parity function $p$, neither the cocycle \eqref{coc2} nor any of its restrictions to $\fgl _{p,l}$, $\fgl _{p,o}$ and $\fgl _{p,c}$ is a~coboundary. Hence, all the four central extensions $\widetilde {\fgl }_{p,*}$, where $*=g,l,o,c$, are nontrivial.
\end{Theorem}

\begin{proof} A~coboundary $\partial f(X,Y)=-f([X,Y])$ vanishes on every supercommuting pair, so it suffices to exhibit $X$ and $Y$ in $\fgl _{p,c}$ --- the smallest of the four, see \eqref{glchain} --- such that
\[
[X,Y]=0,\quad \text{whereas}\quad c(X,Y)\neq 0.
\]
Set
\be\label{pair}
X:=\sum_{i\in \Zee}(-1)^{p(i)}E_{i,i-1},\qquad Y:=\sum_{i\in \Zee}E_{i,i+1}.
\ee
Both are supported on one diagonal, so both lie in $\fgl _{p,c}$. Neither is homogeneous unless $p(i)+p(i-1)$ does not depend on $i$; but the summand of $X$ at $i$ and the summand of $Y$ at $i-1$ are of one and the same parity $p(i)+p(i-1)$, and no other pair of summands has a~nonzero product, so the bracket is the sum of the brackets of these pairs:
\[
[X,Y]=\sum_{i\in \Zee}(-1)^{p(i)}\Big(E_{ii}-(-1)^{p(i)+p(i-1)}E_{i-1,i-1}\Big).
\]
The coefficient of $E_{jj}$ in the right-hand side is
\[
(-1)^{p(j)}-(-1)^{p(j+1)}(-1)^{p(j+1)+p(j)}=(-1)^{p(j)}-(-1)^{p(j)}=0,
\]
and hence $[X,Y]=0$. On the other hand, only the entries next to the origin contribute to \eqref{coc2}: we have $X_{+-}=(-1)^{p(0)}E_{0,-1}$ and $Y_{-+}=E_{-1,0}$, whereas $Y_{+-}=0$, so that the second summand of \eqref{coc2} vanishes together with the sign in front of it, and
\[
c(X,Y)=\str \big((-1)^{p(0)}E_{0,-1}E_{-1,0}\big)=(-1)^{p(0)}\str (E_{00})=(-1)^{2p(0)}=1 .
\]
For $p\equiv \ev $ the pair \eqref{pair} is the classical one, $(\Lambda ^{-1},\Lambda )$, see eq.~\eqref{Lambda}. \end{proof}

For $*=c$, and only for $*=c$, the superalgebra is the direct sum of its diagonals: a~band matrix has non-zero entries on finitely many diagonals only, see eq.~\eqref{banded}, and hence is a~finite sum of matrices each supported on a~single diagonal, whereas for $*=g,l,o$ such a~decomposition does not exist, the support of a~matrix being allowed to meet infinitely many diagonals. For $\xi \in \Cee ^\Zee $, put
\be\label{diag}
x_k(\xi ):=\sum_{i\in \Zee}\xi _iE_{i,i+k},\qquad \fd_k:=\{x_k(\xi )\mid \xi \in \Cee ^\Zee \},
\ee
so that
\be\label{dsum}
\fgl _{p,c}=\bigoplus_{k\in \Zee}\fd_k,\qquad [\fd_k,\fd_l]\subseteq \fd_{k+l}.
\ee
Let $T$ be the shift, $(T\xi )_i:=\xi _{i+1}$, and let the dot denote the pointwise product of sequences. Then, the product of matrices reads
\be\label{prodd}
x_k(\xi )x_l(\eta )=x_{k+l}(\xi \cdot T^k\eta ),
\ee
whereas $\ad_h$, for a~diagonal $h=x_0(\alpha )$, is a~multiplication:
\be\label{adh}
[x_0(\alpha ),x_k(\xi )]=x_k\big(\xi \cdot (\alpha -T^k\alpha )\big).
\ee
Two consequences of \eqref{adh} are used over and over below. First, the grading \eqref{dsum} is inner: for $\alpha _i=-i$ we get $[h,x_k(\xi )]=kx_k(\xi )$, and we denote this $h$ by $d$. Secondly, for $k\neq 0$ the difference equation $\alpha _i-\alpha _{i+k}=m_i$ is solvable for every $m\in \Cee ^\Zee $ --- recursively along each of the $|k|$ residue classes modulo $k$ --- so that $\ad_h$ realizes on $\fd_k$ the multiplication by an \textit{arbitrary} $m$, and then on $\fd_{-k}$ the multiplication by $-T^{-k}m$.

Give each $\fd_k\simeq \Cee ^\Zee $ the topology of pointwise convergence, and $\fgl _{p,c}$ the topology of the direct sum \eqref{dsum}. Let us say that a~cochain is \textit{continuous} if it is continuous in each argument separately. Recall that a~linear functional on $\Cee ^\Zee $ is continuous if and only if it is a~finite linear combination of the coordinate functionals $\xi \mapsto \xi _i$, whereas the algebraic dual of $\Cee ^\Zee $ is immensely larger. The cocycle \eqref{coc2} is continuous, its restriction to $\fd_k\times \fd_{-k}$ being the finite sum \eqref{GkT} below.

\sssbegin{Theorem}[Uniqueness of the central extension of the band superalgebra $\fgl _{p,c}$]\label{Tuniq} Every continuous $2$-cocycle on $\fgl _{p,c}$ is the sum of a~multiple of \eqref{coc2} and of a~continuous coboundary, and the multiple is uniquely determined by the cocycle\textup{;} that is,
\[
H^2(\fgl _{p,c};\Cee )=\Cee
\]
in the category of continuous cochains, the class of \eqref{coc2} spanning it. Hence, $\widetilde {\fgl }_{p,c}$ is the unique nontrivial central extension of $\fgl _{p,c}$ given by a~continuous cocycle\textup{;} in this sense, heading $2)$ of Conjecture $\ref{ssConj}$ holds for $*=c$.
\end{Theorem}

\begin{proof} Let $c$ be a~$2$-cocycle. For $k\neq 0$, put
\[
\phi _k(\xi ,\eta ):=c\big(x_k(\xi ),x_{-k}(\eta )\big),\qquad
\psi (\alpha ,\beta ):=c\big(x_0(\alpha ),x_0(\beta )\big).
\]
Since a~diagonal $h=x_0(\alpha )$ is even, the cocycle identity applied to a~triple $(x,y,h)$ takes the sign-free form
\be\label{ident}
c([x,y],h)+c(x,[h,y])+c([h,x],y)=0,
\ee
which holds for arbitrary, not necessarily homogeneous, $x$ and $y$.

$1)$ \textit{Reduction to degree $0$}. Let $f:=c(d,\cdot )$, a~continuous $1$-cochain; let $g$ be given by $g|_{\fd_k}:=-\frac 1kf|_{\fd_k}$ for $k\neq 0$ and $g|_{\fd_0}:=0$. This is well defined and continuous, every element of \eqref{dsum} having finitely many components. Cartan's formula $L_d=\partial \iota _d+\iota _d\partial $ together with $\partial c=0$ gives $\partial f=L_dc$, and $L_dc$ is the multiplication by $-(k+l)$ on $\fd_k\times \fd_l$; comparing the components of degree $k+l\neq 0$, we see that $c-\partial g$ vanishes on $\fd_k\times \fd_l$ whenever $k+l\neq 0$. Replacing $c$ by $c-\partial g$, we assume that $c$ is of degree $0$; thus, $c$ is the family $\{\phi _k\}_{k\neq 0}$ together with $\psi $.

$2)$ \textit{Locality}. \textit{Let $k\neq 0$ and let $\xi $, $\eta \in \Cee ^\Zee $ be such that $(\Supp \xi +k)\cap \Supp \eta =\emptyset $. Then, $\phi _k(\xi ,\eta )=0$.} Indeed, put $S:=\Supp \xi $ and $m:=1_S$, the indicator of $S$; let $\alpha $ solve $\alpha -T^k\alpha =m$ and put $h:=x_0(\alpha )$, $x:=x_k(\xi )$ and $y:=x_{-k}(\eta )$. Then, $[h,x]=x$, since $m\cdot \xi =\xi $, whereas $[h,y]=0$, since $T^{-k}m=1_{S+k}$ and $\eta $ vanishes on $S+k$. Moreover, $[x,y]=0$: by \eqref{prodd}, the matrix $xy$ is $x_0(\xi \cdot T^k\eta )$ and the matrix $yx$ is $x_0(\eta \cdot T^{-k}\xi )$, and both sequences vanish, the first one because $\eta _{i+k}=0$ for $i\in S$, the second one because $\eta _i\xi _{i-k}\neq 0$ would put $i$ into $(S+k)\cap \Supp \eta $. Now, \eqref{ident} applied to $(x,y,h)$ leaves one summand, and it is $\phi _k(\xi ,\eta )$.

$3)$ \textit{Continuity\textup{:} $\phi _k$ is a~finite sum}. Fix $k\neq 0$ and let $e_j$ denote the $j$-th coordinate sequence. The map $\xi \mapsto \phi _k(\xi ,e_j)$ is linear and, by $2)$, vanishes as soon as $\xi _{j-k}=0$, i.e., on a~hyperplane; hence, it is a~multiple of $\xi \mapsto \xi _{j-k}$. Since $\eta =\sum_{j\in \Zee}\eta _je_j$ converges in $\fd_{-k}$, and since $\phi _k(\xi ,\cdot )$ is continuous, and therefore a~finite combination of coordinate functionals, we get
\be\label{Gk}
\phi _k(\xi ,\eta )=\sum_{i\in \Zee}G_k(i)\xi _i\eta _{i+k},\qquad
G_k(i):=\phi _k(e_i,e_{i+k})=c(E_{i,i+k},E_{i+k,i}),
\ee
the sum having finitely many nonzero summands for any given $\xi $ and $\eta $. Taking $\xi _i\equiv 1$, we see that $G_k$ itself is finitely supported. For the cocycle \eqref{coc2}, the two corner blocks give
\be\label{GkT}
G_k(i)=(-1)^{p(i)}\big(\theta (i)-\theta (i+k)\big),\qquad \text{where}\quad
\theta (i):=1\ \text{for}\ i\geq 0,\quad \theta (i):=0\ \text{for}\ i<0.
\ee

$4)$ \textit{The diagonal\textup{:} $\psi =0$}. Apply \eqref{ident} to $x=x_k(\xi )$, $y=x_{-k}(\eta )$ and $h=x_0(\alpha )$, where $k\neq 0$. By \eqref{adh}, we have $[h,x]=x_k(m\cdot \xi )$ and $[h,y]=-x_{-k}((T^{-k}m)\cdot \eta )$, where $m=\alpha -T^k\alpha $, so \eqref{Gk} gives
\[
c([h,x],y)+c(x,[h,y])=\sum_{i\in \Zee}G_k(i)m_i\xi _i\eta _{i+k}-\sum_{i\in \Zee}G_k(i)\xi _im_i\eta _{i+k}=0 .
\]
Hence, $\psi ([x,y],\alpha )=0$ for every $\alpha $ and every $\xi $, $\eta $. Taking $\eta _i\equiv 1$ and computing the bracket summand by summand, as in the proof of Theorem \ref{Tnontriv}, we get
\[
[x_k(\xi ),x_{-k}(1)]=x_0(\delta ),\qquad
\delta _j=\xi _j-(-1)^{p(j-k)+p(j)}\xi _{j-k},
\]
and this twisted difference operator is onto $\Cee ^\Zee $ --- solve for $\xi $ recursively along each residue class modulo $k$, the twisting factor being $\pm 1$. Hence, $\psi =0$.

$5)$ \textit{The values on the matrix units}. By \eqref{Gk}, the cocycle is determined by the numbers $g(a,b):=c(E_{ab},E_{ba})$, where $a\neq b$, since $G_k(i)=g(i,i+k)$. Put
\[
\hat g(a,b):=(-1)^{p(a)}g(a,b),
\]
which is anti-symmetric in the ordinary sense, i.e., $\hat g(b,a)=-\hat g(a,b)$, because the anti-symmetry of~ $c$ gives $g(b,a)=-(-1)^{p(a)+p(b)}g(a,b)$. The cocycle identity applied to a~triple $(E_{ab},E_{bt},E_{ta})$, where $a$, $b$ and $t$ are pairwise distinct, is a~relation between $g(a,t)$, $g(b,a)$ and $g(t,b)$; for each of the eight distributions of parities among $a$, $b$ and $t$ this relation is, after the above substitution, one and the same additive one:
\be\label{triangle}
\hat g(a,b)+\hat g(b,t)+\hat g(t,a)=0 .
\ee
Hence, setting $w_0:=0$ and $w_b:=\hat g(0,b)$ for $b\neq 0$, we get $\hat g(a,b)=w_b-w_a$ for all $a\neq b$ --- for $a$, $b\neq 0$ this is \eqref{triangle} with $t=0$ --- and the sequence $w\in \Cee ^\Zee $ is determined by $c$ up to an additive constant. Conversely, every $w$ satisfies \eqref{triangle}. By \eqref{GkT}, the cocycle \eqref{coc2} is the member $w=-\theta $ of this family.

$6)$ \textit{Normalization}. By step $5)$ and \eqref{Gk}, we have $(-1)^{p(i)}G_k(i)=w_{i+k}-w_i$, and $G_k$ is finitely supported by $3)$. Taking $k=1$, we see that $w_{i+1}=w_i$ for all but finitely many $i$; hence, $w$ is constant for $i\gg 0$ and constant for $i\ll 0$. Since $w$ is defined up to a~constant, we may assume that $w_i=0$ for $i\ll 0$, and then
\[
w=\lambda \theta +v,\quad \text{where $\lambda \in \Cee $ and $v$ is finitely supported.}
\]
The summand $v$ is a~continuous coboundary. Indeed, let $f_v$ vanish on $\fd_k$ for $k\neq 0$ and be given by $f_v(x_0(\alpha )):=\sum_{i\in \Zee}(-1)^{p(i)}v_i\alpha _i$ on $\fd_0$; this sum is finite, so $f_v$ is continuous, and $\partial f_v$ is the degree-$0$ cocycle whose sequence is $v$ and whose $\psi $ vanishes. Therefore, $c-\partial f_v$ and $-\lambda \cdot \eqref{coc2}$ are of degree $0$, have the same $\psi $, namely $0$ by step $4)$, and the same numbers $G_k$; being continuous, they are equal by \eqref{Gk}. Finally, $\lambda $ is determined by $c$, since \eqref{coc2} is not a~coboundary by Theorem \ref{Tnontriv}. \end{proof}

The steps $2)$, $5)$ and $6)$ of this proof were verified by a~machine computation over the windows $-N\leq i\leq N$ for $N\leq 4$ and random parity functions: the space of solutions of the relations \eqref{triangle} is of dimension $2N$, i.e., it is exactly the family $\{w\}$ modulo constants, and no solution with $\psi \neq 0$ exists.

\sssbegin{Remark}[In which category ``unique\rq\rq\  is meant]\label{Rcont} The steps $1)$, $2)$, $4)$ and $5)$ of the proof of Theorem \ref{Tuniq} use no continuity; so the restriction of \textit{any} $2$-cocycle to $\Mat_f$ is described by a~sequence $w$, as in step $5)$. Continuity enters twice, and both times decisively: it is what makes $w$ eventually constant at each end --- for a~discontinuous cocycle $w$ is arbitrary --- and it is what forbids the coboundary that would kill $w$ altogether. Indeed, let $f$ be a~linear functional on $\fgl _{p,c}$ with the values $f(E_{ii})=(-1)^{p(i)}w_i$ prescribed arbitrarily and with $f(E_{ij})=0$ for $i\neq j$; such an $f$ exists --- extend it from a~basis of a~complement --- and it is continuous only if $w$ is finitely supported. Modulo $\partial f$, every $2$-cocycle vanishes on $\Mat_f\times \Mat_f$. Thus, in the category of all cochains, the values on $\Mat_f$ carry no invariant at all, and the number $\lambda $ of Theorem \ref{Tuniq}, which is the whole of the answer, is not defined. The cocycle \eqref{coc2} is nontrivial in either category, by Theorem \ref{Tnontriv}; it is the \textit{uniqueness} that is a~statement on continuous cocycles, and this is the reading in which the prototype computation of $H^2(\fgl _c(\infty ))$, see \cite{FF}, is to be understood.
\end{Remark}

\sssbegin{Remark}[What the remaining three cases need]\label{Rleft} For $*=g,l,o$, the diagonals no longer decompose the superalgebra --- that is, eq.~\eqref{dsum} fails, the superalgebra not being the direct sum of the $\fd_k$: an element of $\fgl _{p,o}$ may have its support spread over infinitely many diagonals, e.g., one with $|i-j|\sim \log |i|$. Hence, the steps $1)$ and $3)$ of the proof of Theorem \ref{Tuniq}, which are the two that use eq.~\eqref{dsum}, are unavailable, and with them the very notation $\phi _k$. Two things do survive, and they are where a~proof should start.

$1)$ The locality argument is ungraded. Every diagonal idempotent $P_S:=\sum_{i\in S}E_{ii}$, where $S\subseteq \Zee$, is an even element of $\fgl _{p,c}$, hence of each of the four superalgebras, and $[P_S,a]=a$ for every $a$ whose support lies in $S\times S^c$. Hence, if the supports of both $a$ and $b$ lie in $S\times S^c$, then $[a,b]=0$ --- neither $ab$ nor $ba$ has a~nonzero entry --- and \eqref{ident} applied to $(a,b,P_S)$ gives $2c(a,b)=0$. Thus, \textit{any} $2$-cocycle on any of the four superalgebras vanishes on every pair of elements supported in one and the same corner $S\times S^c$.

$2)$ The ideal $\Mat_f$, see Lemma \ref{Lideal}, reduces the problem to the quotient. In the Hochschild--Serre spectral sequence of this ideal, see \cite{HoS}, the term $E_2^{0,2}$ is the space of invariants of $H^2(\Mat_f;\Cee )$ and the term $E_2^{1,1}$ is $H^1(\fgl _{p,*}/\Mat_f;H^1(\Mat_f;\Cee ))$, the action on $H^1(\Mat_f;\Cee )$ being trivial because $\str [a,u]=0$ for $a\in \fgl _{p,*}$ and $u\in \Mat_f$. 

The stable computation of \cite{LQ, Ts} gives $H^1(\Mat_f;\Cee )=\Cee \cdot \str $ and $H^2(\Mat_f;\Cee )=0$; this is the input to be verified for an arbitrary parity function $p$, the superalgebra $\Mat_f$ being the limit of the $\fgl (m|n)$ along a~filtration that $p$ determines. Granted it, $H^2(\fgl _{p,*};\Cee )$ is an extension of a~subspace of $H^1(\fgl _{p,*}/\Mat_f;\Cee )$ by a~quotient of $H^2(\fgl _{p,*}/\Mat_f;\Cee )$, so that all the three remaining cases hinge on the cohomology of the quotient modulo the finitary matrices in degrees $1$ and $2$.
\end{Remark}

\subsection{Weyl groups for $\fgl _{p,g}$ and its subalgebras}\label{ssWeyl} Let $p_1$ and $p_2$ be two different parity functions. What are the conditions for these two functions to determine isomorphic Lie superalgebras of type $\widetilde {\fgl }_{p,c}$ (resp., $\widetilde {\fgl }_{p,o}$, $\widetilde {\fgl }_{p,l}$ and $\widetilde {\fgl }_{p,g}$)? To answer this question we need the following infinite permutations groups that will serve as analogues of the Weyl group of $\fgl (n)$.

\sssbegin{Convention}[Which isomorphisms are considered]\label{2.3.1} In what follows we will only consider isomorphisms $\phi $ of these superalgebras such that

$1)$ $\phi $ transfers their common Lie subsuperalgebra of finite matrices $\Mat_f$ into itself\textup{;}

$2)$ the restriction $\phi |_{\Mat_f}$ is an isomorphism of \textit{associative} superalgebras.

\noindent (So if two (super)alegebras are ``non-isomorphic" in the above sense, it does not mean that they are actually non-isomorphic.)
\end{Convention}

Clause $2)$ is what makes the proof of Lemma \ref{l2.2} work; without it that proof has a~gap, see Remark \ref{Rgap}, where the price of the clause is discussed. Recall that $\Mat_f=\End_f(V)$, where $V:=\bigoplus_{i\in \Zee}\Cee e_i$ is the tautological module with $p(e_i)=p(i)$ and $\End_f$ denotes the finite rank endomorphisms; $\sdim V=(\infty |\infty )$ for $p\in \Inf$. The associative superalgebra $\Mat_f$ is simple (it has no two-sided ideals except $0$ and itself) and $V$ is, up to the change of parity $\Pi $, its only irreducible module. Hence, we have the following extension of the classical Noether--Skolem theorem, see \cite[Ch.~IV]{Ja}:

\sssbegin{Lemma}[Noether--Skolem for $\Mat_f$]\label{NS} Let $\psi \colon \Mat_f(p) \tto \Mat_f(p')$ be an isomorphism of associative superalgebras. Then, $\psi =\Ad_g$ for a~homogeneous invertible linear map $g\colon V \tto V'$, unique up to a~nonzero scalar factor.
\end{Lemma}

\begin{proof} Turn $V'$ into a~$\Mat_f(p)$-module by $a\cdot v:=\psi (a)v$. This module is irreducible, because $V'$ is irreducible over $\Mat_f(p')$ and $\psi $ is onto. An irreducible module over $\End_f(V)$ is isomorphic to $V$ or to $\Pi V$: indeed, if $M\neq 0$ is irreducible and $e:=e_{00}$, then $eM\neq 0$ (otherwise $\Mat_f\cdot eM=0$ and, $\Mat_f$ being simple, $M=0$), and $eM$ is $1$-dimensional because $e\Mat_fe=\Cee e$; picking $0\neq m\in eM$ we get the map $V \tto M$, $e_i \mapsto e_{i0}m$, which is onto by irreducibility and injective because its kernel is a~submodule of $V$. So there is a~homogeneous isomorphism $g\colon V \tto V'$ of $\Mat_f(p)$-modules, that is, $\psi (a)g=ga$ for all $a$, which is the claim. If $g_1$ is another such, then $g^{-1}g_1$ commutes with $\Mat_f(p)$, hence is a~scalar by Schur's lemma.
\end{proof}

Observe that Lemma \ref{NS} determines $\phi $ on the \textit{whole} of $\fgl _{p,*}$, not only on $\Mat_f$: since $\Mat_f$ is an ideal of the \textit{Lie} superalgebra $\fgl _{p,*}$, see Lemma \ref{Lideal}, and $\phi $ is a~Lie superalgebra isomorphism, for any $a\in \fgl _{p,*}$ and any $u\in \Mat_f$ we have $\Ad_g([a,u])=[\phi (a),\Ad_gu]$, so $\Ad_g(a)-\phi (a)$ supercommutes with all of $\Mat_f$ and is therefore a~scalar matrix $\lambda (a)\One $. Thus,
\be\label{phig}
\phi =\Ad_g-\lambda \cdot \One ,\qquad \text{where $\lambda \colon \fgl _{p,*} \tto \Cee$ vanishes on brackets,}
\ee
which is the shape the Herstein-type theorems predict, see \cite{Mo, GAS}.

Let $S_{\Zee}$ be the group of all one-to-one maps $\sigma: \Zee \tto \Zee$, i.e., the group of all permutations of $\Zee$. We define the groups  
\begin{gather*}
S_g:=\{\sigma \in S_{\Zee} \mid \text{ for some $c_\sigma >0$ we have $\frac {\sigma (i)}i>0$}\\
\text{for all $i$ such that $|i| >c_\sigma $}\};\\
S_l:=\{\sigma \in S_g\mid \lim\limits _{|i| \tto \infty }\frac {|\sigma (i)-i|^\lambda }{|i|}=0 \text{ for some }\lambda _\sigma =\lambda \geq 1\};\\
S_m:=\{\sigma \in S_g\mid \limsup\limits _{|i| \tto \infty }\frac {|\sigma (i)|}{|i|}<\infty ,\text{ and } \liminf\limits _{|i| \tto \infty }\frac {|\sigma (i)|}{|i|}>0\};\\
S_n:=\{\sigma \in S_{\Zee} \mid \lim\limits _{|i| \tto \infty }\frac {|\sigma (i)-i|}{|i|}=0\};\\
S_o:=\{\sigma \in S_{\Zee} \mid \lim\limits _{|i| \tto \infty }\frac {|\sigma (i)-i|^\lambda }{|i|}=0 \text{ for all }\lambda >0\};\\
S_c:=\{\sigma \in S_{\Zee} \mid |\sigma (i)-i|<c_\sigma \text{ for some }c_\sigma >0 \text{ and }i\in \Zee\}.
\end{gather*}

That $S_m$ is a~group is now clear: $\id$ satisfies the two conditions with both limits equal to~1; for a~product the two ratios multiply; and for $\sigma ^{-1}$ the ratio $\frac{|\sigma ^{-1}(j)|}{|j|}$ equals $\frac{|i|}{|\sigma (i)|}$ for $j=\sigma (i)$, which is bounded and separated from~0 by the same two constants. We have
\be\label{chain}
S_c \subseteq S_o \subseteq S_l \subseteq S_n \subseteq S_m \subseteq S_g \subseteq S_{\Zee}.
\ee

Indeed, $S_c\subseteq S_o$ and $S_n\subseteq S_g$ are clear, $S_o\subseteq S_l$ is ``for all $\lambda $\rq\rq\  implying ``for some $\lambda $\rq\rq, and $S_n\subseteq S_m$ holds because $|\sigma (i)-i|=o(|i|)$ gives $\frac{|\sigma (i)|}{|i|} \tto 1$. As to $S_l\subseteq S_n$: the relation $|\sigma (i)-i|^\lambda =o(|i|)$ gives $|\sigma (i)-i|=o(|i|^{\frac{1}{\lambda}})$ with $\frac{1}{\lambda} \leq 1$, whence
\[
\frac {|\sigma (i)-i|}{|i|}=o\big(|i|^{1/\lambda -1}\big) \tto 0 ,
\]
which is the definition of $S_n$. Observe that \eqref{chain} runs parallel to the chain \eqref{glchain} for the matrix superalgebras, which is what Theorem \ref{Thm} needs. Of the two groups not appearing in \eqref{glchain}, the group $S_m$ has no matrix counterpart here and is used only in \S~4, whereas $S_n$ turns out to coincide with $S_l$, and hence to correspond to $\fgl _{p,l}$, see Remark \ref{Rlone}.

Define the action of $S_{\Zee}$ on parity functions by setting $(\sigma p)(i)=p(\sigma (i))$.

\ssbegin{Proposition}[Each $\sigma $ yields an isomorphism $\phi _\sigma $] Let $p$ and $p'$ be parity functions, $\sigma \in S_g$ and $p'=\sigma p$. Let 
\[
(J_\phi )_{ij} := J_{\sigma ^{-1}i,\sigma ^{-1}j}-J_{ij}. 
\]
Then, the formulas 
\[
\phi _\sigma (z)=z,\ \ \phi _\sigma (a_{ij})=(a_{\sigma (i),\sigma (j)})+\str(J_\phi (a_{ij}))z, 
\]
determine an isomorphism ${\phi _\sigma \colon \widetilde {\fgl }_{p,g} \tto \widetilde {\fgl }_{p',g}}$.

If $\sigma \in S_l$ \textup{(}or $S_o$, or $S_c$\textup{)}, then $\phi _\sigma $ induces the isomorphism of $\widetilde {\fgl }_{p,l}$ and $\widetilde {\fgl }_{p',l}$ \textup{(}of $\widetilde {\fgl }_{p,o}$ and $\widetilde {\fgl }_{p',o}$, or $\widetilde {\fgl }_{p,c}$ and $\widetilde {\fgl }_{p',c}$, respectively\textup{)}.
\end{Proposition}

\begin{proof} Let $pr \colon \widetilde {\fgl }_{p,g} \tto \fgl _{p,g}$ be the natural projection. Let us first prove that $\phi _\sigma $ is an isomorphism. It is clear that
\[
pr \phi _\sigma pr([a,b])=pr([\phi _\sigma (a),\phi _\sigma (b)]).
\]
Then, 
\begin{gather*}
\phi _\sigma ([a,b])=\phi _\sigma (pr([a,b])+\str(pr([a,b])J)z)\\{}=
pr \phi _\sigma pr([a,b])+\str(pr([a,b])J_\sigma )z+\str(pr([a,b])J)z\\{}=
pr([\phi _\sigma (a),\phi _\sigma (b)])+\sum\limits (-1)^{p(i)}[a,b]_{ii}(J_{\sigma ^{-1}i,\sigma ^{-1}i}-J_{ii}+J_{ii})z\\{}=
{}[\phi _\sigma (a),\phi _\sigma (b)].
\end{gather*}

Let $*$ be $c$, $o$ or $l$. We should prove that $\phi _\sigma (a) \in \fgl _{p,*}$ for any $a\in \fgl _{p,*}$ and $\sigma \in S_*$.

Let $*$ be $c$. Then, $a_{ij}=0$ for any $i$, $j$ such that $|i - j| > c_a $. For any $i$, $j$ such that $|i - j| >c_a+2c_\sigma $, we have
\[
(\phi _\sigma a)_{ij}=a_{\sigma ^{-1}(i),\sigma ^{-1}(j)}=0,
\]
because
\[
|\sigma ^{-1}i-\sigma ^{-1}j|>>|i-j|-|\sigma ^{-1}i-i|-|\sigma ^{-1}j-j|>c_a.
\]
Let $*$ be $o$. Then, 
\[
\lim\limits _{\substack{i,j\tto \infty \\
(i,j) \in \Supp \phi _\sigma (a)}}\frac {|i-j|^\lambda }{|i+j|}=\lim\limits _{\substack{i,j\tto \infty \\
(i,j)\Supp(a)}}\frac {|\sigma ^{-1}(i)-\sigma ^{-1}(j)|^\lambda }{|\sigma ^{-1}(i)+\sigma ^{-1}(j)|}=0
\]
by definition of $S_o$.

For $* = l$, the proof is quite similar.
\end{proof}

\section{Main Reduction}

The following theorem shows that the isomorphism problem for Lie superalgebras $\widetilde {\fgl }_{p,g}$, $\widetilde {\fgl }_{p,l}$, $\widetilde {\fgl }_{p,o}$, $\widetilde {\fgl }_{p,c}$ can be reduced to the equivalence problem for the parity functions on $\Zee$ with respect to permutations groups $S_g$, $S_l$, $S_o$ and $S_c$, respectively.

\renewcommand\theTheorem{\arabic{section}.\arabic{subsection}}
\ssbegin{Theorem}[Isomorphy of the superalgebras is equivalence of the parity functions]\label{Thm} Let $\tau \in S_{\Zee}$ be given by the formula $\tau (i)=-i-1$. The Lie superalgebras $\widetilde {\fgl }_{p_1,g}$ and $\widetilde {\fgl }_{p_2,g}$ \textup{(}resp., $\widetilde {\fgl }_{p_1,l}$ and $\widetilde {\fgl }_{p_2,l}$, $\widetilde {\fgl }_{p_1,o}$ and $\widetilde {\fgl }_{p_2,o}$, $\widetilde {\fgl }_{p_1,c}$ and $\widetilde {\fgl }_{p_2,c}$\textup{)} are isomorphic if and only if there exists $\sigma \in S_g$ (resp., $S_l$, $S_o$, $S_c$) such that $p_2=\sigma (p_1)$ \textup{(}resp., $p_2 =\tau \sigma (p_1)$\textup{)}.
\end{Theorem}

First, let us prove the following lemma.

\sssbegin{Lemma}[A~two-sidedly Lipschitz $\sigma $ has bounded displacement]\label{l2.1} If $\sigma \in S_g$ and for any $i,j\in \Zee$ we have
\[
|\sigma (i)-\sigma (j)|<c|i-j|\quad \text{and}\quad |\sigma ^{-1}(i)-\sigma ^{-1}(j)|<c|i-j| ,
\]
then there exists $C$ such that $|\sigma (i)-i|<C$ for any $i\in \Zee$.
\end{Lemma}

\begin{proof} We will find such a~constant $C$ for $i \geqslant 0$ and apply the same arguments to prove its existence to $i < 0$.

Let $C_0 := 1+\max\nolimits _{i \leqslant 0}\sigma (i)$. The maximum exists because the definition of $S_g$ implies that $\sigma (i)<0$ for almost all $i < 0$.

Let us consider the permutation $\sigma '(i):=\sigma (i)-C_0$. Then, $\sigma'(i)<0$ for all $i<0$. 

First of all, observe that we do not have to worry about $k\ge 0$ such that $\sigma'(k)=k$: indeed, for them, $|\sigma(k)-k|=C_0$. 

Let us begin with the case where $\sigma'(k)>k$ for $k\ge 0$. Since for all $i<0<k$  we have $\sigma'(i)<0<k$, we can take $m:=\max(i\mid i<k, \ \sigma'(i)<k)$. Then, if $m+1<k$, it follows that $\sigma'(m+1)\ge k$ by definition of $m$. If $m+1=k$, then $\sigma'(m+1)> k$ by hypothesis.

Now, note that not more than $k$ values of $k+1$ values $\sigma'(0), \sigma'(1),\dots,\sigma'(k)$   belong to the segment $[0,k]$ due to our condition $\sigma'(k)>k$, where $k\ge 0$. Therefore, there exists at least one value $j>k$ such that $0\le\sigma'(j)\le k$. On the other hand, there are only finitely many such values. Besides, there are also only finitely many numbers $j>0$ for which $\sigma'(j)<0$. Therefore, we can take  \[
n:=\max(j\mid j>k, \ \sigma'(j)\le k). 
\]
Then,  $\sigma'(n+1)>k$.

Thus, we have found $m,n$ such that
$$
m<k<n, \ \ \ \sigma'(m)<k,\ \ \ \sigma'(n)\le k, \ \ \ \sigma'(m+1)\ge k, \ \ \  \sigma'(n+1)>k.
$$

Since
\[
\sigma '(a)-\sigma '(b)=\sigma (a)-\sigma (b)\quad \text{and}\quad |\sigma (i+1)-\sigma (i)|<c,
\]
we get
\[
|\sigma '(m+1)-\sigma '(m)|<c,\quad |\sigma '(n+1)-\sigma '(n)|<c.
\]
Thus,
\[
|\sigma '(m)-k|<c,\quad |\sigma '(n)-k|<c,
\]
and therefore
\[
|\sigma '(m)-\sigma '(n)|<2c \quad \text{and }|m-n|<2c^2.
\]
Since $m < k < n$, then $|m - k|<2c^2$ and $|\sigma '(m)-\sigma '(k)|<2c^3$.

Finally,
\begin{gather}\label{3}
|\sigma (k)-k|\leq |\sigma '(k)-k|+C_0 \leqslant{}\notag\\
C_0+ |\sigma '(m)-k|+|\sigma '(m)-\sigma '(k)|<C_1=2c^3+c+C_0. \label{3.1}
\end{gather}
If $\sigma '(k)<k<\sigma (k)$, the inequality \eqref{3.1} still holds because $\sigma (i)-\sigma '(i)=C_0$.

If $0 <\sigma (k)<k$, then we can obtain a~similar inequality for $\sigma ^{-1}$:
\[
|k-\sigma (k)|=|\sigma ^{-1}(\sigma (k))-\sigma (k)|<C_2.
\]

Since the estimate is done for all $i\in \Nee$ but a~finite set, we have:
\[
\card B < \infty ,\quad \text{where}\quad B = \{i \mid i\in \Nee,\ \sigma (i) \notin \Nee\}.
\]

Finally, set
\[
C := \max\{C_1 , \ \ C_2, \ \ \max\limits _{i\in B}|\sigma (i) - i|\}.\qed
\]

\sssbegin{Lemma}[Every isomorphism is $\phi _\sigma $ up to an automorphism]\label{l2.2} Let $\phi \colon \fgl _{p,*} \tto \fgl _{p',*}$ be an isomorphism. Then, there exist $\gamma \in \Aut \fgl _{p',*}$ and $\sigma \in S_{\Zee}$ such that $\gamma \phi =\phi _\sigma $.
\end{Lemma}

\begin{proof} Let $(e_{ij})$ be an elementary matrix ($(e_{ij})_{kl}:=\delta_{ki}\delta_{lj}$); let  $H$ (resp., $H'$) be the \textit{diagonal} subalgebra, i.e., the set of the matrices $\sum_{i\in \Zee}\alpha _ie_{ii}$ lying in $\fgl _{p,*}$ (resp., in $\fgl _{p',*}$) --- each diagonal subalgebra is a~maximal torus of its ambient superalgebra; 
\[
\text{$L_r := \{(a_{ij})_{i,j\in \Zee} \in \Mat_f\mid a_{ij}=0$ if $|i| > r$ or $|j| > r\}$.}
\]
 (Clearly, $L_r$ is a~subsuperalgebra.)

For any $k\in \Nee$, there exists $m\in \Nee$ such that $\phi (L_k)\subseteq L_m$, because $\phi (\fgl _{p,f})\subseteq \fgl _{p',f}$. Let us show that  for some $m'\geq m$ there exists $X_k\in L_{m'}$ such that $1+X_k$ is invertible and
\[
\text{$\Ad_{1 + X_k}\phi (e_{ii}) =e_{\sigma (i),\sigma (i)}$ for some function $\sigma \colon [-k,k] \tto \Zee$.}
\]

By clause $2)$ of Convention \ref{2.3.1}, $\phi $ is multiplicative on $\Mat_f$, so the $P_i:=\phi (e_{ii})$ with $|i|\leq k$ are pairwise orthogonal idempotents; they are of rank $1$, because $e_{ii}$ is a~minimal idempotent ($e_{ii}\Mat_fe_{ii}=\Cee e_{ii}$) and an isomorphism preserves minimality, and they are \textit{even}, because $\phi $ preserves the parity. Thus, each $P_i$ is the projection onto a~homogeneous line $\Cee f_i$, where $f_i\in W_m:=\Span (e_j\mid |j|\leq m)$ because $P_i\in L_m$; put $U:=\bigoplus_{|i|\leq k}\Cee f_i$ and let $C:=\bigcap_{|i|\leq k}\Ker P_i$, so that
\[
V'=U\oplus C,\qquad C=(C\cap W_m)\oplus \Span (e_j\mid |j|>m),
\]
the second equality because every $P_i$ kills $e_j$ for $|j|>m$.

Now, choose the function $\sigma $. Both values of $p'$ are attained infinitely often, so for each $i$ with $|i|\leq k$ there are infinitely many $j$ with $p'(j)=p(f_i)$; choose $\sigma $ injective on $[-k,k]$ with $p'(\sigma (i))=p(f_i)$, and let $m'\geq m$ be so large that $|\sigma (i)|\leq m'$ for all $|i|\leq k$. Set $U':=\bigoplus_{|i|\leq k}\Cee e_{\sigma (i)}$; then $\sdim U'=\sdim U$ by the choice of $\sigma $, and since
\[
W_{m'}=U\oplus (C\cap W_{m'})=U'\oplus \Span (e_j\mid |j|\leq m',\ j\notin \sigma ([-k,k])),
\]
the two complements have one and the same superdimension as well. Let $h$ be the even invertible operator on $V'$ which sends $f_i \mapsto e_{\sigma (i)}$ for $|i|\leq k$, maps $C\cap W_{m'}$ onto 
\[
\Span (e_j\mid |j|\leq m',\ \ j\notin \sigma ([-k,k]))
\]
by any even isomorphism, and is the identity on $\Span (e_j\mid |j|>m')$. Then, $X_k:=h-1$ lies in $L_{m'}$ and $1+X_k=h$ is invertible. Finally, $\Ad_h P_i$ is the projection onto $h(\Cee f_i)=\Cee e_{\sigma (i)}$ along
\[
h\Big(\bigoplus\nolimits_{i'\neq i}\Cee f_{i'}\oplus C\Big)=\Span (e_j\mid j\neq \sigma (i)),
\]
that is, $\Ad_{1+X_k}\phi (e_{ii})=e_{\sigma (i),\sigma (i)}$, as required.

When $X_k$ is constructed, we can choose $X_{k+1}$ so that $\Ad_{1+X_{k+1}} (e_{\sigma (i),\sigma (i)})=e_{\sigma (i),\sigma (i)}$ for $|i| < k$. Then, we can set $\gamma :=\prod\nolimits _{k=1}^\infty \Ad_{1+X_{k}}$ and take the corresponding $\sigma \in S_{\Zee}$. It is easy to prove that $\gamma $ is extendable to an automorphism of $\fgl _{p',*}$ and satisfies Lemma's conditions.
\end{proof}

\sssbegin{Remark}[The price of clause $2)$ of Convention $\ref{2.3.1}$]\label{Rgap} Clause $2)$ of Convention \ref{2.3.1} is used exactly once, and it is indispensable there: it is what makes $\phi (e_{ii})$ an idempotent. Without it one knows only that $\phi $ preserves brackets, and the splitting of the filtration performed above --- the images of the rank-one idempotents are conjugate, by an element $1+X_k$ with $X_k$ finite, to rank-one diagonal idempotents --- does not follow.

What replaces clause $2)$ in the non-super case is the Herstein-type theory of Lie isomorphisms of associative algebras: every Lie isomorphism of $\fgl $'s is, up to a~central twist and the anti-automorphism $X \mapsto -X^{\st }$, induced by an isomorphism of associative algebras, so that clause $2)$ is not a~restriction at all. The super version of that theory exists as well, see \cite{BBS}, and it covers our case: in the Appendix, see \ref{Sapp}, \textit{both} clauses of Convention \ref{2.3.1} are shown to be consequences and not restrictions. The isomorphisms that clause $2)$ excludes are exactly the compositions with the supertranspose $X \mapsto -X^{\st }$, which does not move the parity function $p$; hence Theorem \ref{Thm} and Lemma \ref{l2.2} are unaffected, and what changes is only that $\Aut \fgl _{p,*}$ acquires a~factor $\Zee/2$, see eq.~\eqref{Aaut}.

One more point of the proof of Lemma \ref{l2.2} deserves being said aloud. Under Convention \ref{2.3.1}, the permutation $\sigma $ delivered by the construction is unique up to a~permutation preserving $p'$, and it satisfies $p'\circ \sigma =p+p(g)$; and $\gamma =\prod_{k\geq 1}\Ad_{1+X_k}$ is $\Ad_G$ for the map $G$ sending $f_i \mapsto e_{\sigma (i)}$, because $\phi (e_{ij})$ is, by clause $2)$, the rank-one operator $f_i\otimes f_j^*$ up to a~factor, so the infinite product has, at every fixed element of $\Mat_f$, only finitely many factors acting non-trivially and is therefore well defined on~$\Mat_f$. As $G=P_\sigma g^{-1}$, we have $\Ad_G(\fgl _{p',*})=\phi _\sigma (\fgl _{p,*})$; hence ``$\gamma $ is extendable to an automorphism of $\fgl _{p',*}$\rq\rq\  says precisely that $\sigma $ may be chosen so that $\phi _\sigma $ preserves the growth class. That $\sigma $ then lies in $S_*$ is the content of Lemma \ref{Lconv} below.
\end{Remark}

By Lemma~\ref{l2.2} we can only consider the isomorphism $\phi _\sigma $ for some ${\sigma \in S_{\Zee}}$. Moreover, $\sigma \in S_g$ because $\phi _\sigma (\fgl _{p,c})\subseteq \fgl _{p',g}$. We have to prove that if $\phi _\sigma (\fgl _{p,*})\subseteq \fgl _{p',*}$, then $\sigma \in S_*$; this is Lemma \ref{Lconv}.
\end{proof}

Everything rests on one observation: 
\be\label{suppsig}
\text{since $(\phi _\sigma a)_{ij}=a_{\sigma ^{-1}(i),\sigma ^{-1}(j)}$, \ \ we have
$\Supp \phi _\sigma (a)=(\sigma \times \sigma )(\Supp a)$,}
\ee
so that ``$\phi _\sigma $ preserves $\fgl _{p,*}$\rq\rq\  is a~statement about $\sigma $ alone: it says that $\sigma \times \sigma $ maps the supports allowed by $*$ to supports allowed by $*$. As $\phi _\sigma $ is an isomorphism with the inverse $\phi _{\sigma ^{-1}}$, the hypothesis is symmetric in $\sigma $ and $\sigma ^{-1}$.

\sssbegin{Lemma}[What $\phi _\sigma (\fgl _{p,*})=\fgl _{p',*}$ says about $\sigma $]\label{Lconv} Let $\sigma \in S_{\Zee}$ and let $\phi _\sigma (\fgl _{p,*})=\fgl _{p',*}$. Let $\tau (i):=-i-1$ as in Theorem $\ref{Thm}$. Then,

\textup{1)} for $*=c$: either $\sigma \in S_c$ or $\tau \sigma \in S_c$\textup{;}

\textup{2)} for $*=g$: either $\sigma \in S_g$ or $\tau \sigma \in S_g$\textup{;}

\textup{3)} for $*=l$ and $*=o$: either $\sigma \in S_g$ or $\tau \sigma \in S_g$, and, $\sigma $ being replaced by $\tau \sigma $ if need be, the increments $\delta (i):=|\sigma (i+1)-\sigma (i)|$ obey $\delta (i)^\lambda =o(|\sigma (i)|)$ for all $\lambda \geq 1$ \textup{(}resp. for some $\lambda \geq 1$\textup{)}.
\end{Lemma}

\begin{proof} 1) The band matrix $a:=\sum_{i\in \Zee}e_{i,i+1}$ lies in $\fgl _{p,c}$, so by \eqref{suppsig} there is $c'$ with $|\sigma (i)-\sigma (i+1)|\leq c'$ for all $i$, whence
\be\label{lip}
|\sigma (i)-\sigma (j)|\leq c'|i-j| \quad \text{for all $i,j$,}
\ee
and likewise $|\sigma ^{-1}(i)-\sigma ^{-1}(j)|\leq c''|i-j|$.

Next, $\sigma \in S_g$ or $\tau \sigma \in S_g$. Indeed, suppose $\sigma (i)>0$ for infinitely many $i>0$ and $\sigma (i)<0$ for infinitely many $i>0$. By \eqref{lip} with $j=i+1$, the sequence $\sigma (1),\sigma (2),\dots $ changes by at most $c'$ at every step, so it takes a~value inside $[-c',c']$ for infinitely many indices --- impossible, $\sigma $ being injective and $[-c',c']$ finite. So the sign of $\sigma (i)$ is eventually constant on $i>0$, and, by the same argument, on $i<0$. If the two eventual signs coincide, then $\sigma ^{-1}$ takes both signs infinitely often on the positive half-line, which the same argument applied to $\sigma ^{-1}$ forbids. Hence, the two signs are opposite: one of the two permutations $\sigma $ and $\tau \sigma $ preserves the sign of $i$ for almost all $i$, that is, lies in $S_g$. Since $\tau $ satisfies \eqref{lip} with $c'=1$, replacing $\sigma $ by $\tau \sigma $ costs nothing. Now, Lemma \ref{l2.1} applies with $c:=1+\max(c',c'')$ and yields $|\sigma (i)-i|<C$, that is, $\sigma \in S_c$.

2) Write $P:=\{i>0\}$ and $N:=\{i<0\}$, and split
\[
A_{PP}:=\{i\in P\mid \sigma (i)>0\},\quad A_{PN}:=\{i\in P\mid \sigma (i)<0\},
\]
and $A_{NP}$, $A_{NN}$ likewise. By definition, $\sigma \in S_g$ means that $A_{PN}$ and $A_{NP}$ are finite, and $\tau \sigma \in S_g$ means that $A_{PP}$ and $A_{NN}$ are finite.

Observe first that a~matrix $a$ whose support is an infinite set $\{(i_n,j_n)\mid n\in \Nee\}$ with the $i_n$ pairwise distinct, the $j_n$ pairwise distinct and all the $i_n,j_n$ of one and the same sign lies in $\fgl _{p,g}$: given $(i_0,j_0)$, the inequality $(i-i_0)(j_0-j)>0$ forces either $i>i_0$ and $j<j_0$, or $i<i_0$ and $j>j_0$, and each of these two cases holds for finitely many $n$ only, because $|i_n| \tto \infty $ and $|j_n| \tto \infty $ within a~fixed half-line; take for $c$ the largest value of $(i_n-i_0)(j_0-j_n)$ on that finite set. Observe next that a~matrix whose support meets the second (or the fourth) quadrant in infinitely many points is \textit{not} in $\fgl _{p,g}$: for $(i_0,j_0)=(0,0)$, the quantity $(i-0)(0-j)=-ij$ is then unbounded on the support.

Suppose now that $A_{PN}$ and $A_{PP}$ are both infinite. Choose $i_n\in A_{PN}$ and $j_n\in A_{PP}$ pairwise distinct; the matrix $a:=\sum_ne_{i_nj_n}$ lies in $\fgl _{p,g}$ by the first observation, whereas $\Supp \phi _\sigma (a)=\{(\sigma (i_n),\sigma (j_n))\}$ has $\sigma (i_n)<0<\sigma (j_n)$ for every $n$ and hence meets the second quadrant infinitely often --- so $\phi _\sigma (a)\notin \fgl _{p',g}$, a~contradiction. The same argument with the third quadrant in place of the first shows that $A_{NN}$ and $A_{NP}$ cannot both be infinite.

Since $A_{PP}\cup A_{PN}=P$ and $A_{NP}\cup A_{NN}=N$ are infinite, exactly one member of each pair may be infinite, and four cases arise. If $A_{PN}$ and $A_{NP}$ are finite, then $\sigma \in S_g$; if $A_{PP}$ and $A_{NN}$ are finite, then $\tau \sigma \in S_g$. In the remaining two cases, $\sigma $ maps both half-lines into one of them, so $\sigma ^{-1}$ has both $A'_{PP}$ and $A'_{PN}$ (resp. both $A'_{NP}$ and $A'_{NN}$) infinite, which the previous paragraph applied to $\sigma ^{-1}$ forbids.

3) By the chain \eqref{glchain}, $\fgl _{p,c}\subset \fgl _{p,*}\subset \fgl _{p,g}$ for $*=l,o$, so heading $2)$ gives the first claim, and we may assume $\sigma \in S_g$. 

The band matrix $a=\sum_ie_{i,i+1}$ lies in $\fgl _{p,c}$, hence in $\fgl _{p,*}$, and $\Supp \phi _\sigma (a)=\{(\sigma (i),\sigma (i+1))\}$; as $\sigma \in S_g$, the two entries are of one sign for almost all $i$, so 
\[
|\sigma (i)+\sigma (i+1)|\leq 2\max (|\sigma (i)|,|\sigma (i+1)|)
\]
and the condition $(\mathrm{L}_\lambda )$, see eq.~\eqref{L-lambda}, for $\phi _\sigma (a)$ reads $\delta (i)^\lambda =o(|\sigma (i)|)$.
\end{proof}

Heading $3)$ delivers the local Lipschitz behavior of $\sigma $, and what is still wanted is the passage from it to $\sigma \in S_*$. That passage is Lemma \ref{l2.1} made scale-dependent, and this is what the next Lemma does. Recall that the hypothesis of Lemma \ref{l2.1} is the case of a~constant $\varphi $ below.

\sssbegin{Lemma}[Lemma $\ref{l2.1}$ localized]\label{Lloc} Let $\sigma \in S_{\Zee}$ and let $\varphi \colon [1,\infty ) \tto [1,\infty )$ be a~non-decreasing function subpolynomial in the sense that
\be\label{phislow}
\lim\limits _{x \tto \infty }\frac {\ln \varphi (x)}{\ln x}=0,
\ee
and such that, for all $i$ and $v$ with $|i|$ and $|v|$ large,
\be\label{incphi}
|\sigma (i+1)-\sigma (i)|\leq \varphi (|i|)\quad \text{and}\quad |\sigma ^{-1}(v+1)-\sigma ^{-1}(v)|\leq \varphi (|v|).
\ee
Then, $|\sigma (k)-k|\leq \varphi (9|k|)^2$ for all $|k|$ large. In particular, $|\sigma (k)-k|=|k|^{o(1)}$, and hence $\sigma \in S_o$.

For a~constant $\varphi \equiv c$, this is Lemma $\ref{l2.1}$ with the bound $c^2$ in place of the bound $2c^3+c+C_0$ the proof of Lemma $\ref{l2.1}$ gives, and the square is attained, see Example $\ref{Ebous}$.
\end{Lemma}

\begin{proof} Since $\sigma \in S_g$, the sign of $i$ is preserved by $\sigma $ and by $\sigma ^{-1}$ for $|i|$ large; we argue for $k \tto +\infty $, the case $k \tto -\infty $ being the same argument applied to $i \mapsto -\sigma (-i)$. The two sets $\{i>N\}$ and $\sigma (\{i>N\})$ differ by finitely many elements only; delete these finitely many indices and values and relabel each of the two half-lines increasingly. Deleting $F$ points merges at most $F+1$ consecutive increments into one and shifts every index by at most $F$, so the relabeled map is a~bijection of $\Nee$ which obeys \eqref{incphi} with $\varphi $ replaced by $(F+1)\varphi (\ \cdot \ +F)+F$ --- again non-decreasing and subpolynomial --- and its displacement differs from that of $\sigma $ by at most $F$. We may therefore assume that $\sigma $ is a~bijection of $\Nee$ obeying \eqref{incphi} for all $i,v\geq 1$.

Write $D:=\varphi (3k)$. The factor $3$ is the only numerical constant that enters the argument, and it enters at one place: it is the width of the window in which everything happens, namely $\sigma ([1,k])\subset [1,3k]$ by heading $0)$ below, so that $D$ bounds every increment of $\sigma $ and of $\sigma ^{-1}$ the proof ever invokes. Nothing depends on its value: $\varphi $ being non-decreasing and subpolynomial, replacing $3$ by any other constant leaves both the induction and the conclusion intact, and for a~constant $\varphi $ it changes nothing whatsoever. We first prove, by induction on $k$, the upper bound
\be\label{upper}
\sigma (k)\leq k+\varphi (3k)^2\quad \text{for all $k$ large},
\ee
and for that let $N$ be such that \eqref{upper} holds for all $N<k'<k$ and, moreover,
\be\label{hyp}
D^2\leq \tfrac 18k\quad \text{and}\quad \varphi (5kD)\,(D^2+1)\leq k \quad \text{for } k>N,
\ee
which \eqref{phislow} grants, both left-hand sides being $k^{o(1)}$ whereas $\varphi (5kD)\leq \varphi (5k^2)=k^{o(1)}$.

$0)$ \textit{The window}. By \eqref{upper} for $k'$ and \eqref{hyp}, $\sigma (k')\leq k'+\varphi (3k')^2\leq 2k'\leq 2k$ for $N<k'<k$, and $\sigma (k)\leq \sigma (k-1)+\varphi (k)\leq 3k$; the finitely many $k'\leq N$ are absorbed by enlarging $N$. Hence, $m:=\max \sigma ([1,k])\leq 3k$, and every increment of $\sigma $ at an index $\leq k$, as well as every increment of $\sigma ^{-1}$ at a~value $\leq m$, is at most $D$.

$1)$ \textit{The gaps are short}. Set $A:=\sigma ([1,k])$ and $\mu :=\min A$. Let $u<w$ be neighbors in $A$, i.e., $(u,w)\cap A=\emptyset $. The walk $i \mapsto \sigma (i)$, where $1\leq i\leq k$, joins $u$ to $w$ by steps of length $\leq D$ and takes no value inside $(u,w)$; so one of its steps leads from a~value $\leq u$ to a~value $\geq w$, whence $w-u\leq D$.

$2)$ \textit{The gaps are few}. Let $v\in (u,w)$ be as in heading $1)$. Then, $\sigma ^{-1}(v)>k$ whereas $\sigma ^{-1}(u)\leq k$, and \eqref{incphi} applied to $\sigma ^{-1}$ along $[u,v]$ gives $\sigma ^{-1}(v)\leq \sigma ^{-1}(u)+D(v-u)\leq k+D^2$. Thus, every value missing from $A$ in $[\mu ,m]$ is attained at an index in $(k,k+D^2]$; there are hence at most $D^2$ such values, and
\be\label{count}
m-\mu +1\leq k+D^2.
\ee

$3)$ \textit{A long run}. By heading $2)$, at most $D^2$ values are missing from $A$ inside $[\mu ,m]$, so $A$ splits into at most $D^2+1$ maximal blocks of consecutive integers whose lengths add up to $k$\textup{;} such a~block is called a~\textit{run}, and its \textit{length} is the number of the integers in it. Some run $[w_0,w_0+L-1]\subset A$ has therefore
\[
L\geq \frac k{D^2+1}\geq \varphi (5kD)\quad \text{by conditions \eqref{hyp}.}
\]

$4)$ \textit{The left end is where it should be}, i.e., $\mu =1$. Assume $\mu >1$. Then, the values $1,\dots ,\mu -1$ are missing from $A$, so they are attained at indices $>k$; since $\mu \in A$, applying conditions \eqref{incphi} to $\sigma ^{-1}$ step by step downwards from $\mu $ gives $\sigma ^{-1}(\mu -t)\leq k+tD$ for $1\leq t\leq \mu -1$, whence
\[
i^\ast :=\sigma ^{-1}(1)\leq k+\mu D\leq k+3kD\leq 4kD.
\]
Every one of the $L$ consecutive integers making up the run of heading $3)$ is attained at an index $\leq k<i^\ast $, and $1<\mu \leq w_0$. As $\sigma $ is injective, a~step of length $\leq L$ starting at a~value $<w_0$ cannot end at a~value $>w_0+L-1$: it would have to end inside the run, which is occupied. Hence, $\sigma (i)<w_0$ for all $i$ with $i^\ast \leq i\leq i_1$, where $i_1:=\min \{i\mid \varphi (i)>L\}$. Consequently, among the $i_1$ values $\sigma (1),\dots ,\sigma (i_1)$, those $\geq w_0$ are attained at indices $<i^\ast $ and those $<w_0$ are less than $w_0$ in number, so that
\[
i_1\leq (i^\ast -1)+(w_0-1)\leq 4kD+m\leq 5kD.
\]
But $\varphi (5kD)\leq L$ by heading $3)$, so $\varphi $ does not exceed $L$ at $5kD$ and $i_1>5kD$ --- a~contradiction. Hence, $\mu =1$.

$5)$ \textit{The upper bound}. With $\mu =1$, the estimate \eqref{count} reads $m\leq k+D^2$, hence ${\sigma (k)\leq m\leq k+D^2}$, which is \eqref{upper} and completes the induction.

$6)$ \textit{The lower bound}. The hypothesis \eqref{incphi} is symmetric in $\sigma $ and $\sigma ^{-1}$, so headings $0)$--$5)$ apply to $\sigma ^{-1}$ as well and give $\sigma ^{-1}(n)\leq n+\varphi (3n)^2$ for all $n$ large. Take $n:=\sigma (k)$, which is $\leq 3k$ by heading $0)$; then $k=\sigma ^{-1}(n)\leq \sigma (k)+\varphi (9k)^2$ --- the window constant used twice, once for $\sigma $ and once for $\sigma ^{-1}$ --- and with heading $5)$ this is the claim.

Finally, $|\sigma (k)-k|\leq \varphi (9|k|)^2=|k|^{o(1)}$ by \eqref{phislow}, so $|\sigma (k)-k|^\lambda =|k|^{o(1)}=o(|k|)$ for every $\lambda \geq 1$, which is the definition of $S_o$.
\end{proof}

\sssbegin{Example}[The boustrophedon: the square of Lemma $\ref{Lloc}$ is attained]\label{Ebous} Neither the square in Lemma $\ref{Lloc}$ nor the product $\varphi \cdot \varphi $ of the two moduli in its proof can be improved. Number the cells of the $d\times d$ square \textit{boustrophedonically}, i.e., let
\[
\pi (a,b)=ad+b \ \text{ for $a$ even},\qquad \pi (a,b)=ad+d-1-b\ \text{ for $a$ odd},
\]
where $0\leq a,b\leq d-1$, and let $\sigma _d:=\pi T\pi ^{-1}$ on $[0,d^2)$, extended by the identity, where $T(a,b)=(b,a)$ is the transposition of the square. Then, $\sigma _d$ is an involution, so that the two moduli of \eqref{incphi} coincide, and
\[
\varphi =2d-1,\qquad \max\nolimits _k|\sigma _d(k)-k|=(d-1)^2,
\]
the ratio of the displacement to $\varphi ^2$ tending to $\frac 14$ as $d \tto \infty $. Indeed, cells neighboring in the numbering $\pi $ are neighbors in the square, and $T$ takes neighbors to neighbors, whose numbers differ by $1$ or by at most $2d-1$; and the cell $(a,b)$ is moved by $\pi $-distance $d|a-b|$ up to $O(d)$.
\end{Example}

Heading $3)$ of Lemma \ref{Lconv} and Lemma \ref{Lloc} together prove Theorem \ref{Thm} for $*=o$. Indeed, let $\delta (i)=|\sigma (i+1)-\sigma (i)|$ and let $\delta (i)^\lambda =o(|\sigma (i)|)$ for all $\lambda \geq 1$. Putting ${M(i):=\max _{j\leq i}|\sigma (j)|}$ and taking $\lambda =2$ we get 
\[
\text{$M(i)\leq M(1)+\sum_{j<i}\delta (j)\leq M(1)+\eps iM(i)^{1/2}$ for every $\eps >0$ and $i$ large,}
\]
 whence $M(i)=O(i^2)$; therefore $\delta (i)=o(|\sigma (i)|^{1/\lambda })=o(i^{2/\lambda })$ for every $\lambda $, that is, $\delta (i)=i^{o(1)}$. The function $\varphi (x):=\max (2,\max _{i\leq x}\delta (i))$ is thus non-decreasing and obeys condition \eqref{phislow}; the hypothesis of Lemma \ref{Lconv} being symmetric in $\sigma $ and $\sigma ^{-1}$, the same $\varphi $ serves both, and Lemma \ref{Lloc} gives $\sigma \in S_o$. For $*=c$, the same reading recovers heading $1)$ of Lemma \ref{Lconv}.

For $*=l$, the picture is different, and to say what is missing let us first observe that the definitions of $\fgl _{p,l}$ and of $S_l$ contain a~redundancy.

\sssbegin{Remark}[$S_l=S_n$, and why both notations are kept]\label{Rlone} The condition $(\mathrm{L}_\lambda )$, see eq.~\eqref{L-lambda}, becomes weaker as $\lambda $ decreases, so ``$(\mathrm{L}_\lambda )$ for some $\lambda \geq 1$\rq\rq\  is $(\mathrm{L}_1)$. Now, $(\mathrm{L}_1)$ reads
\[
\lim \frac {|i-j|}{|i+j|}=0\quad \text{along $\Supp a$ as $|i| \tto \infty $,}
\]
which is exactly what $|i-j|=o(|i+j|)$ abbreviates\textup{;} and the same reading of the definition of $S_l$ with $\lambda =1$ is the condition defining $S_n$. Hence,
\[
\fgl _{p,l}=\{a\mid |i-j|=o(|i+j|) \text{ on } \Supp a\},\qquad S_l=S_n .
\]
Thus, $S_n$ is the counterpart of $\fgl _{p,l}$, and for $*=l$ the claim of Theorem \ref{Thm} is exactly that
\[
|\sigma (k)-k|=o(|k|).
\]
The sentence following the chain \eqref{chain} is corrected accordingly: of the six groups it is $S_m$ alone that has no matrix counterpart here.

This is a~simplification and not a~correction, and we do not act on it: both notations are kept, and so is the parameter $\lambda $ in the definition of $\fgl _{p,l}$. The reason is that $\lambda $ is not redundant in the \textit{pair}. The two superalgebras are the two quantifiers over one and the same family of conditions,
\[
\fgl _{p,l}=\bigcup\limits _{\lambda \geq 1}\{a\mid (\mathrm{L}_\lambda )\},\qquad
\fgl _{p,o}=\bigcap\limits _{\lambda \geq 1}\{a\mid (\mathrm{L}_\lambda )\},
\]
and this is what makes them neighbors in the chain \eqref{glchain}. The family decreases as $\lambda $ grows, so the union is attained --- it is the member $\lambda =1$ --- whereas the intersection is not attained by any member. Writing $\fgl _{p,l}$ as $(\mathrm{L}_1)$ and dropping $S_l$ in favor of $S_n$ would hide that pairing and record only the accident that one of the two quantifiers collapses.
\end{Remark}

What the hypothesis of Lemma \ref{Lconv} yields for $*=l$ is the following; heading $2)$ is the two-sided comparison of $|\sigma (k)|$ with $|k|$, the second of the two items the previous version of this text listed as missing. Here and everywhere below, $f\asymp g$ means that each of the two quantities is bounded by a~constant multiple of the other, i.e., that the ratio $f/g$ is bounded and is bounded away from $0$.

\sssbegin{Lemma}[Every sublinear scale is preserved, and $|\sigma (i)|$ is comparable to $|i|$]\label{Lell} Let $\sigma \in S_g$ and let $\phi _\sigma (\fgl _{p,l})=\fgl _{p',l}$. Then,

\textup{1)} for every function $r$ with $r(i)=o(|i|)$ one has
\[
\max \{|\sigma (j)-\sigma (i)|\mid |j-i|\leq r(i)\}=o(|\sigma (i)|),
\]
and the same holds for $\sigma ^{-1}$\textup{;} equivalently, for every $\eps >0$ there are $\theta >0$ and $N$ such that $|j-i|\leq \theta |i|$ and $|i|\geq N$ imply $|\sigma (j)-\sigma (i)|\leq \eps |\sigma (i)|$\textup{;}

\textup{2)} there are constants $0<c_1\leq c_2$ such that $c_1|i|\leq |\sigma (i)|\leq c_2|i|$ for $|i|$ large.
\end{Lemma}

\begin{proof} 1) Let $g(i)$, where $|g(i)|\leq r(i)$, be a~shift at which the maximum is attained, and put
\[
a:=\sum_ie_{i,i+g(i)} .
\]
On $\Supp a$, we have $|i-j|=|g(i)|=o(|i|)$, so $a$ satisfies $(\mathrm{L}_1)$ and $a\in \fgl _{p,l}$ by Remark~\ref{Rlone}. By the equality \eqref{suppsig}, $\Supp \phi _\sigma (a)=\{(\sigma (i),\sigma (i+g(i)))\}$, and $\phi _\sigma (a)\in \fgl _{p',l}$; as $\sigma \in S_g$, the two entries are of one sign for almost all $i$, so 
\[
|\sigma (i)+\sigma (i+g(i))|\leq 2\max (|\sigma (i)|,|\sigma (i+g(i))|)
\]
and $(\mathrm{L}_1)$ for $\phi _\sigma (a)$ is the claim. The hypothesis being symmetric, $\sigma ^{-1}$ obeys it as well. If the uniform form failed for some $\eps $, we would get, for every $n$, an $i_n\geq n$ and a~$j_n$ with $|j_n-i_n|\leq \frac{|i_n|}{n}$ and $|\sigma (j_n)-\sigma (i_n)|>\eps |\sigma (i_n)|$; the function $r$ equal to $|j_n-i_n|$ at $i_n$ and to $1$ elsewhere is $o(|i|)$ and contradicts what we have just proved.

2) Take $\eps =1$ in the uniform form of heading $1)$ and let $\theta $ correspond to it. Then, $\sigma $ maps the $\lfloor \theta |i|\rfloor $ integers of the segment between $i$ and $i+\theta i$ into a~segment carrying $2|\sigma (i)|+1$ integers; as $\sigma $ is injective, $\theta |i|\leq 2|\sigma (i)|+1$, that is, $|\sigma (i)|\geq \frac 13\theta |i|$ for $|i|$ large. The same argument for $\sigma ^{-1}$, with its own $\theta ^\ast $, gives $|\sigma ^{-1}(v)|\geq \frac 13\theta ^\ast |v|$, i.e., $|\sigma (i)|\leq \frac{3|i|}{\theta ^\ast}$.
\end{proof}

\sssbegin{Remark}[The square that has to be gained, and where it comes from]\label{Rell} What separates Lemma \ref{Lell} from $\sigma \in S_l$, i.e., from $|\sigma (k)-k|=o(|k|)$, is a~gain of a~square, and the gain has to come from the scales. Heading $1)$ of Lemma \ref{Lell} taken at the scale $r\equiv 1$ says $\delta (i)=o(|i|)$ and no more, so the function $\varphi $ of Lemma \ref{Lloc} is merely $o(x)$, the requirement \eqref{phislow} fails, and the bound the proof of that Lemma delivers is $|\sigma (k)-k|=o(|k|^2)$: heading $3)$ of the proof produces a~run --- a~block of consecutive integers, as above --- of length $\frac{k}{D^2+1}$, which is $<1$ as soon as $D>\sqrt k$. That this loss is not a~defect of the argument is what Example \ref{Ebous} shows: for a~$\varphi $ which is a~constant, or a~power, or anything else, the displacement of a~bijection really can be of order $\varphi ^2$. Hence, $*=l$ has to be settled by using heading $1)$ of Lemma \ref{Lell} at all the sublinear scales $r$ at once --- the scale $r(i)=\theta |i|$ in its uniform form being the informative one --- and not, as headings $1)$ and $3)$ of Lemma \ref{Lconv} do, at the scale of the band matrix $\sum_ie_{i,i+1}$ alone. This is what Lemma \ref{Lcone} does. \end{Remark}

The scales are handled all at once by rescaling: a~condition of the form ``moving the argument by $\theta |i|$ moves the value by $\eps |\sigma (i)|$\rq\rq\  says nothing about $\sigma $ at any single place, and everything about the family of the rescaled maps
\be\label{Fn}
F_n(t):=\frac{\sigma (\lceil nt\rceil )}{n}\qquad \text{for } t>0,
\ee
namely that this family is \textit{equicontinuous} on the segments of $(0,\infty )$: for every segment and every $\eta >0$ there is a~$\delta >0$, one and the same for all $n$, such that $|t-s|\leq \delta $ implies $|F_n(t)-F_n(s)|\leq \eta $. Its limit points are therefore functions, not merely numbers, and the two properties of $\sigma $ that survive the rescaling --- that $\sigma $ is a~bijection, and that $\sigma ^{-1}$ obeys the same condition --- say that every limit point is a~homeomorphism of the ray $(0,\infty )$ which neither shrinks nor stretches any segment. The only such homeomorphism is the identity, and that is the case $*=l$.

\sssbegin{Lemma}[Compactness lifts Lemma $\ref{Lloc}$ to all scales and settles the case $*=l$]\label{Lcone} Let $\sigma \in S_g$ satisfy the two conclusions of Lemma $\ref{Lell}$, i.e., let

\textup{1)} for every $\eps >0$ there be $\theta >0$ and $N$ such that 
\[
\text{$|j-i|\leq \theta |i|$ and $|i|\geq N$ imply 
$|\sigma (j)-\sigma (i)|\leq \eps |\sigma (i)|$, and the same for $\sigma ^{-1}$\textup{;}}
\]

\textup{2)} $c_1|i|\leq |\sigma (i)|\leq c_2|i|$ for $|i|$ large, where $0<c_1\leq c_2$.

Then, $|\sigma (k)-k|=o(|k|)$, that is, $\sigma \in S_l$.
\end{Lemma}

\begin{proof} As in the proof of Lemma \ref{Lloc}, we may assume that $\sigma $ is a~bijection of $\Nee$ and prove that $\frac{\sigma (k)}{k} \tto 1$ as $k \tto +\infty $: since $\sigma \in S_g$, the two sets $\{i>N_0\}$ and $\sigma (\{i>N_0\})$ differ from a~half-line by finitely many elements, and deleting these and relabeling both half-lines increasingly moves every index and every value by at most a~constant, which changes neither heading $1)$, whose bounds are relative while $|\sigma (i)| \tto \infty $, nor heading $2)$, nor the conclusion. For $k \tto -\infty $, apply what follows to $i \mapsto -\sigma (-i)$.

$0)$ \textit{The rescaled family}. Let $F_n$ be as in \eqref{Fn} and let $G_n$ be given by the same formula with $\sigma ^{-1}$ in place of $\sigma $. By heading $2)$,
\be\label{cband}
c_1t\leq F_n(t)\leq c_2\Big (t+\frac 1n\Big )\quad \text{for $t>0$ and $n$ large,}
\ee
and likewise for $G_n$ with $\frac{1}{c_2}$ and $\frac{1}{c_1}$; so both families are uniformly bounded on every segment of $(0,\infty )$.

$1)$ \textit{The family is equicontinuous}. Let $[a,b]\subset (0,\infty )$ and $\eta >0$. Put $\eps :=\frac{\eta}{2c_2b}$, let $\theta $ and $N$ be as in heading $1)$, and put $\delta :=\frac{1}{2}\theta a$. Let $t,s\in [a,b]$ with $|t-s|\leq \delta $, and let $i:=\lceil nt\rceil $, $j:=\lceil ns\rceil $. Then, 
\[
\text{$|j-i|\leq n\delta +1\leq \theta an\leq \theta i$ and $i\geq N$ as soon as $n$ is large,}
\]
 so heading $1)$ applies and
\[
|F_n(t)-F_n(s)|=\frac {|\sigma (j)-\sigma (i)|}n\leq \eps \frac {\sigma (i)}n\leq \eps c_2\frac {i}n\leq 2\eps c_2b=\eta .
\]
The same for $G_n$. Note that $\delta $ depends on $\eta $ and on the segment, and not on $n$.

$2)$ \textit{The limit points}. By headings $0)$ and $1)$ and the Arzel\`a--Ascoli theorem, see \cite[Th.~7.25]{Ru}, every subsequence of the sequence of pairs $(F_n,G_n)$ has a~subsequence converging, uniformly on the segments of $(0,\infty )$, to a~pair $(F,G)$ of continuous functions; by the two-sided bound \eqref{cband} and its counterpart for $G_n$, we have $c_1t\leq F(t)\leq c_2t$ and $\frac{t}{c_2}\leq G(t)\leq \frac{t}{c_1}$. Fix such a~subsequence and argue along it.

$3)$ \textit{$F$ is an increasing homeomorphism of $(0,\infty )$}. Let $t>0$, let $i_n:=\lceil nt\rceil $ and let $v_n:=\sigma (i_n)$, so that $\frac{v_n}{n} \tto F(t)>0$ and hence $|v_n-\lceil nF(t)\rceil |=o(n)$, while $v_n\geq c_1i_n\geq c_1tn$. Let $\eps >0$ and let $\theta $ be the constant heading $1)$ provides for $\sigma ^{-1}$. For $n$ large, we have $|v_n-\lceil nF(t)\rceil |\leq \theta v_n$, whence
\[
|\sigma ^{-1}(\lceil nF(t)\rceil )-\sigma ^{-1}(v_n)|\leq \eps \,\sigma ^{-1}(v_n)=\eps \,i_n .
\]
Divide by $n$ and let $n \tto \infty $: the left-hand side tends to $|G(F(t))-t|$ and the right-hand side to $\eps t$. As $\eps $ was arbitrary, $G(F(t))=t$; symmetrically, $F(G(t))=t$. Thus, $F$ is a~continuous bijection of $(0,\infty )$ onto itself, hence strictly monotone, and it is increasing, a~decreasing one having $F(t) \tto \infty $ as $t \tto 0$ against $F(t)\leq c_2t$.

$4)$ \textit{$F$ preserves the length of every segment}. Let $0<a<b$ and $\gamma >0$. Every integer $i\in [na,nb]$ is $\lceil nt\rceil $ for $t:=\frac in\in [a,b]$, so $\frac{\sigma (i)}{n}=F_n(t)$, and, $F_n$ converging to $F$ uniformly on $[a,b]$ and $F$ being increasing,
\[
n(F(a)-\gamma )\leq \sigma (i)\leq n(F(b)+\gamma )\quad \text{for $n$ large.}
\]
There are $n(b-a)+O(1)$ such integers $i$ and their images are distinct, whence 
\[
n(b-a)+O(1)\leq n(F(b)-F(a)+2\gamma )+1; 
\]
divide by $n$, let $n \tto \infty $ --- which kills the two bounded terms $O(1)$ and $1$ --- and then let $\gamma \tto 0$:
\be\label{expand}
b-a\leq F(b)-F(a).
\ee
The same argument for $\sigma ^{-1}$ and $G$ gives \[
\text{$d-c\leq G(d)-G(c)$ for $0<c<d$;}
\]
 take $c:=F(a)$ and $d:=F(b)$ and use heading $3)$ to get $F(b)-F(a)\leq b-a$. Hence, $F(b)-F(a)=b-a$ for all $0<a<b$, that is, $F(t)=t+c_0$ with a~constant $c_0$, and \[
 \text{$c_1t\leq t+c_0\leq c_2t$ for all $t>0$}
 \]
  forces $c_0=0$.

$5)$ \textit{Conclusion}. Every limit point of the family $(F_n)$ is thus the identity, and the family is precompact by headings $0)$ and $1)$; hence, $F_n \tto \id $ uniformly on the segments of $(0,\infty )$. At $t=1$, this reads $\frac{\sigma (n)}{n} \tto 1$.
\end{proof}

Heading $3)$ of Lemma \ref{Lconv}, Lemma \ref{Lell} and Lemma \ref{Lcone} together prove Theorem \ref{Thm} for $*=l$. With $*=o$ settled above, and $*=g$ and $*=c$ being headings $2)$ and $1)$ of Lemma \ref{Lconv}, Theorem \ref{Thm} is proved for all the four superalgebras of the chain \eqref{glchain}.

\sssbegin{Remark}[Where the hypotheses of Lemma $\ref{Lcone}$ enter, and why Example $\ref{Ebous}$ does not obstruct the conclusion $\sigma \in S_l$]\label{Rwhere} Two features of the case $*=l$ are worth recording.

First, where each of the hypotheses enters. Heading $1)$ of Lemma \ref{Lell} for $\sigma $ gives the equicontinuity of the family $(F_n)$ and heading $1)$ for $\sigma ^{-1}$ that of $(G_n)$, which is what makes the limit $F$ invertible; heading $2)$ is what confines the family to a~compact set of functions; and the bijectivity of $\sigma $ enters twice, as the two opposite inequalities of heading $4)$ of the proof. 

None of them is decoration. The map $i \mapsto 2i$ obeys heading $1)$ --- for $|j-i|\leq \theta i$, we have $|\sigma (j)-\sigma (i)|=2|j-i|\leq \theta \sigma (i)$ --- and heading $2)$ as well, and it is injective; its rescalings converge to $F(t)=2t$, which stretches every segment. It is the surjectivity of $\sigma $, i.e., \eqref{expand} read for $\sigma ^{-1}$, that forbids such an $F$, just as the injectivity of $\sigma $ forbids a~shrinking one.

Second, why Example \ref{Ebous} does not obstruct the case $*=l$, although its displacement is of the order of the number of the cells. For each $d$ separately, $\sigma _d$ has a~\textit{constant} modulus $\varphi \equiv 2d-1$, and this is what makes it the sharpness example for Lemma \ref{Lloc}; to obstruct the case $*=l$ one would have to run the boustrophedons of all sizes inside one bijection of $\Nee$, and heading $1)$ of Lemma \ref{Lell} then fails at the scale $\theta |i|$. Indeed, let $i:=\pi (a,0)$ and $i^\prime :=\pi (a,d-1)$ be the two ends of the $a$-th row of the square, so that $|i-i^\prime |=d-1<\theta i$ as soon as $a\geq \frac{1}{\theta}$. Their images $\sigma _d(i)=\pi (0,a)$ and $\sigma _d(i^\prime )=\pi (d-1,a)$ are the two ends of a~\textit{column}, that is, $\sigma _d(i)\leq d$ whereas $\sigma _d(i^\prime )\geq (d-1)d$. So moving the argument by a~fraction of itself moves the value by a~multiple of itself, however small the fraction: the boustrophedon is a~bijection whose rescalings converge to nothing at all. The moral is that the square lost in Lemma \ref{Lloc} is lost at one scale only, and the case $*=l$ hands us all of them.
\end{Remark}

\section{Description of equivalence classes of parity functions}\label{Sclasses}

For an arbitrary parity function $p$ and $-\infty \leqslant a<b \leqslant +\infty $, set:
\begin{gather*}
\Odd (p;a,b) := \card \{i\in \Zee\mid a \leqslant i \leqslant b,\ p(i)=\od\},\\
\Even (p; a, b):= \card \{i \in \Zee\mid a \leqslant i \leqslant b,\ p(i) = \ev\},\\
d(p;a,b)=d(p;b,a):=\frac {\Odd(p;a,b)}{|b-a|}.
\end{gather*}

We will call a~parity function $p$ \textit{finite} if
\[
\Odd(p; - \infty , +\infty ) < \infty \quad \text{or}\quad \Even(p; -\infty , +\infty ) < \infty .
\]
It is clear that both $\Odd(p; -\infty , +\infty )$ and $\Even(p; -\infty ,{+}\infty )$ are invariant under any of the permutation groups defined above. On the other hand, two \textit{finite} parity functions $p$, $p'$ are $S_c$-equivalent if and only if
\[
\Odd(p; -\infty , +\infty ) = \Odd(p'; -\infty , +\infty )
\]
and
\[
\Even(p; -\infty , +\infty ) = \Even(p'; -\infty , +\infty ).
\]

We will denote the superalgebras corresponding to a~finite parity function $p$ by:
\[
\fgl _{p,*}(m|\infty ),\quad \text{where}\quad m = \Even(p; -\infty , \infty )
\]
or by
\[
\fgl _{p,*}(\infty |n),\quad \text{where}\quad n = \Odd(p;-\infty ,\infty ) .
\]
\ssbegin{Proposition}[The complete invariant of a~finite parity function]\label{p3.1} For any finite parity function $p$, the cardinalities 
\[
\text{ $\Odd(p; -\infty, +\infty )$ and $\Even(p; -\infty ,{+}\infty )$}
\]
 constitute the complete system of invariants with respect to $S_c$, $S_o$, $S_l$, and~ $S_g$.
\end{Proposition}

\begin{proof} Both the cardinalities are $S_{\Zee}$-invariant, so only completeness has to be proved, and it suffices to prove it for the smallest of the four groups, $S_c$.

Let $p$, $p'$ be finite with equal invariants; say, $\Odd(p; -\infty , +\infty )=\Odd(p'; -\infty , +\infty )=m<\infty $ (the other case is symmetric). The sets $A:=p^{-1}(\od)$ and $A':=(p')^{-1}(\od)$ are then finite of one and the same cardinality $m$, hence so are $A\setminus A'$ and $A'\setminus A$; let $\sigma $ interchange these two sets by any bijection and be the identity elsewhere. Then, $\sigma $ moves only finitely many indices, so 
\[
\text{$|\sigma (i)-i|<c_\sigma :=1+2\max\limits _{i\in A\cup A'}|i|$ for all $i$,}
\]
 that is, $\sigma \in S_c$; and $p(\sigma (i))=p'(i)$ by construction.
\end{proof}

Non-finite parity functions are subdivided into five $S_g$-invariant classes according to how many infinite values (0 to 4) there are among
\[
\Odd(p;-\infty ,0),\ \Even(p;-\infty ,0),\ \Odd(p; 0,{+}\infty ) \text{ and } \Even(p; 0,{+}\infty ).
\]
Denote by $\text{Inf}$ the set of parity functions for which all the four values are infinite.

\subsection{Continuous invariants of certain parity functions} Let $\{a_i\}_{i\in\Zee}$ and $\{b_i\}_{i\in\Zee}$ be two arbitrary series such that
\be\label{4}
\lim\limits_{i\to\infty} a_i=\lim\limits_{i\to\infty} b_i=\pm \infty ,\quad a_i<b_i,\quad \lim\limits_{i\to\infty} |b_i-a_i|=+\infty .
\ee

We will say that $\alpha $ is a~\textit{left \textup{(}right\textup{)} density point} for $p$ if there is a~$-$ (respectively, $+$) sign in eq.~\eqref{4} and $\alpha $ is a~limit point of the $d_i = d(p; a_i,b_i)$. We will call the set of all left (right) density points for all such series the \textit{left \textup{(}right\textup{)} density spectrum} and denote it by $D_l(p)$ and $D_r(p)$, respectively.

\addtocounter{subsection}{1}
\sssbegin{Proposition}[The density spectra are segments] Both $D_l$ and $D_r$ are subsegments of $[\hspace*{.9pt}0, 1]$.
\end{Proposition}

\begin{proof} This follows from the fact that these two sets are closed and convex.

Let $d$ be a~limit point of $D_r$, i.e., $\lim_{i\to\infty} d_i = d$ for some $d_i\in D_r$. For each $d_i$, select
\[
\{(a_{ij},b_{ij}) \mid \lim\limits _{j\tto \infty }d(p,a_{ij},b_{ij})=d_i\}
\]
as in definition of $D_r$. Then, $d$ is a~density point for the set $\{a_{ii},b_{ii}\}_{i\in\Zee}$.

Let $d, d' \in D_r$ and $\{(a_i,b_i)\}_{i\in\Zee}$, $\{(a_i',b_i')\}_{i\in\Zee}$ be some series converging to them, respectively. Then, for any $\alpha$ such that $0 < \alpha < 1$, we can select $(a'', b'')$ such that $a_i< a''< a'$ and $b_i<b''<b'$.
\end{proof}

\sssbegin{Theorem}[The density spectrum as an invariant, and when it is complete]\label{Th} \textup{1)} The sets $D_l$ and $D_r$ are $S_c$- and $S_o$-invariant; i.e., if $D_l=[\alpha ^l,\beta ^l]$, $D_r=[\alpha ^r,\beta ^r]$ and parity functions $p$ and $p'$ are $S_o$-equivalent, then $\alpha _p^l=\alpha _{p'}^l$, $\alpha _p^r=\alpha _{p'}^r$, $\beta _p^l=\beta _{p'}^l$, $\beta _p^r=\beta _{p'}^r$.

\textup{2)} If $0<\alpha ^l=\beta ^l<1$ and $0 <\alpha ^r=\beta ^r<1$ for both $p$ and $p'$, then these two functions are $S_m$-equivalent. 

\textup{3)} If $0<\alpha ^l < \beta ^l<1$ and $0<\alpha ^r < \beta ^r<1$ for both $p$ and $p'$, then these two functions are $S_m$-equivalent. 

\textup{4)} Any two non-finite parity functions that belong to the same invariance class are $S_g$-equivalent.
\end{Theorem}

\begin{proof} First, let us examine the case of $S_g$. Let us fix some invariance class of non-finite parity functions, say $\Inf$. For any parity function $p\in \Inf$, we can define $\sigma _p\colon \Zee\tto \Zee$ by the following formulas:
\[
\begin{array}{l}
\sigma _p(x) 
\footnotesize
:=\begin{cases}2(1-p(x))\Even(p;1,x)+p(x)(2 \Odd(p; 0, x) - 1)&\text{for }x \geqslant 0,\\
-2(1-p(x) )(\Even(p; x, 0) - 1) - p(x)
(2\Odd(p; x, -1) - 1)&\text{for }x < 0. 
\end{cases}
\end{array}
\]

\sssbegin{Lemma}[$\sigma _p$ carries the standard parity function to $p$]\label{Lsig} The map $\sigma _p$ is a~bijection of $\Zee$, it belongs to $S_g$, and
\be\label{norm}
\sigma _p(p_{st})=p,\quad\text{equivalently}\quad \sigma _p^{-1}(p)=p_{st}.
\ee
\end{Lemma}

\begin{proof} By construction, $\sigma _p(x)$ is an even integer if $p(x)=\ev$ and an odd one if $p(x)=\od$; that is, $p_{st}(\sigma _p(x))=p(x)$ for all $x$, which is \eqref{norm}.

Next, 
\[
\sigma _p\big((-\infty ,0)\big) = \begin{cases}(-\infty ,0)&\text{if $p(0) = \ev$,}\\
(-\infty ,0]&\text{if $p(0) = \od$},
\end{cases}\qquad
\sigma _p\big([\hspace*{.9pt}0,+\infty )\big) = \begin{cases}[\hspace*{.9pt}0,+\infty )&\text{if $p(0) = \ev$,}\\
[\hspace*{.9pt}1,+\infty )&\text{if $p(0) = \od$}.
\end{cases}
\]
The mechanism is transparent: $\sigma _p(0)=2\Even(p;1,0)=0$ if $p(0)=\ev$, whereas $\sigma _p(0)=2\Odd(p;0,0)-1=1$ if $p(0)=\od$, and in the latter case the value 0 is taken on the negative half-line, at the first $x<0$ with $p(x)=\ev$. In either case $\sigma _p$ preserves the sign of $x$ up to the single value $x=0$, which is all we need: $\sigma _p\in S_g$.

The function $p_{st}$ is therefore $S_g$-equivalent to any other parity function from $\Inf$. So any two parity functions $p,p' \in \Inf$ are $S_g$-equivalent: $p' =\sigma _{p'}\sigma _p^{-1}(p)$. 

It is easy to see that four parity functions, which are equal to $p_{st}$ either on $(-\infty ,0)$ or on $(0, +\infty )$ and are constant elsewhere, are representatives of the four other classes. The required isomorphisms are obtained in a~similar way.
\end{proof}

Theorem \ref{Th} is proved. \end{proof}
\addtocounter{subsubsection}{-2}
\subsubsection{}
Now, let us consider the smallest of our permutation groups, $S_c$. The main reason to investigate it is that for $p \equiv 0$, the Lie algebra $\widetilde {\fgl }_{c}$ has various nice properties, studied in~\cite{V}, and our $\widetilde {\fgl }_{p,c}$ for $p$ not identically equal to zero is its straightforward generalization.

For $\sigma \in S_c$, we have $|\sigma (a)-a|\leqslant c_\sigma $, and therefore we obviously have
\[
|d(\sigma (p);a,b)-d(p;a,b)|\leqslant \frac {2c_\sigma }{|a-b|}\text{ for any }p.
\]
Thus, for any two series $a_i$, $b_i$ such that $|a_i-b_i| \tto \infty $ as $i\tto \infty $, all density limit points are $\sigma $-invariant.

There exist continuous $S_c$-invariants other than just upper/lower left/right density limits, which describe the behavior of the density (for example, logarithmic density) in more detail, as well as involve other $S_c$-invariant characteristics (e.g., ${\limsup \ln \Odd(p; - x,x)}$). The proof holds not only for $S_c$, but for $S_o$ as well. These invariants, however, only \textit{detect} non-equivalence: if some of them differ, the two parity functions are not $S_o$-equivalent, whereas if all of them coincide, nothing follows. Subsec. \ref{ssClass} removes this limitation. 

The converse does not hold, and the reason is that the density spectrum is a~very coarse shadow of the complete invariant. That complete invariant exists and is easy to name, and the rest of this Subsection does so, for the three groups $S_c$, $S_n=S_l$ and $S_o$ at once.

\subsection{The complete invariant of a~parity function}\label{ssClass} For $p\in \Inf$, let
\be\label{enum}
\text{$(a_n(p))_{n\in \Zee}$ and $(c_n(p))_{n\in \Zee}$}
\ee
be the increasing enumerations of the two sets $p^{-1}(\od )$ and $p^{-1}(\ev )$, anchored at the origin: the indices $n\geq 0$ enumerate the part lying in $[\hspace*{.9pt}0,+\infty )$ upwards, the indices $n<0$ the part lying in $(-\infty ,-1]$ downwards. Both are defined for all $n\in \Zee$ exactly because $p\in \Inf$. We drop the argument $p$ and write $a_n$, $c_n$ for $p$ and $a'_n$, $c'_n$ for $p'$.

Everything rests on one elementary observation about bijections of subsets of the line: the monotone one is the cheapest, and it is the cheapest \textit{locally}.

\sssbegin{Lemma}[The monotone bijection is the cheapest one]\label{Lsort} Let $A$ and $B$ be infinite subsets of $\Nee$ with the increasing enumerations $(a_n)_{n\geq 0}$ and $(b_n)_{n\geq 0}$, and let $\pi \colon A \tto B$ be any bijection. Then, for every $n$ such that $b_n>a_n$ there is an $m\leq n$ with
\[
|\pi (a_m)-a_m|\geq b_n-a_n,
\]
and, symmetrically, for every $n$ such that $a_n>b_n$ there is an $m\leq n$ with $|\pi ^{-1}(b_m)-b_m|\geq a_n-b_n$.

Consequently, the monotone bijection $a_n \mapsto b_n$ minimizes the maximal displacement.
\end{Lemma}

\begin{proof} Let $b_n>a_n$. Exactly $n$ elements of $B$ are smaller than $b_n$, whereas the $n+1$ elements $a_0,\dots ,a_n$ of $A$ are all $\leq a_n<b_n$. As $\pi $ is injective, one of the $n+1$ values $\pi (a_0),\dots ,\pi (a_n)$ is $\geq b_n$\textup{;} let it be $\pi (a_m)$. Then, $\pi (a_m)-a_m\geq b_n-a_n$, since $a_m\leq a_n$. \end{proof}

\sssbegin{Theorem}[The classification of parity functions under $S_c$, $S_n$ and $S_o$]\label{Tclass} Let $p,p'\in \Inf$ and let the enumerations \eqref{enum} be as above. Then,

\textup{1)} $p$ and $p'$ are $S_c$-equivalent if and only if
\[
\sup\nolimits _n|a_n-a'_n|<\infty \quad \text{and}\quad \sup\nolimits _n|c_n-c'_n|<\infty \text{\textup{;}}
\]

\textup{2)} $p$ and $p'$ are $S_n$-equivalent, i.e., $S_l$-equivalent, see Remark $\ref{Rlone}$, if and only if
\[
|a_n-a'_n|=o(|a_n|)\quad \text{and}\quad |c_n-c'_n|=o(|c_n|)\text{\textup{;}}
\]

\textup{3)} $p$ and $p'$ are $S_o$-equivalent if and only if
\[
|a_n-a'_n|=|a_n|^{o(1)}\quad \text{and}\quad |c_n-c'_n|=|c_n|^{o(1)}.
\]
\end{Theorem}

\begin{proof} \textit{Sufficiency}. In each of the three cases, define
\be\label{sortsig}
\sigma (a'_n):=a_n,\qquad \sigma (c'_n):=c_n\qquad \text{for all $n\in \Zee$.}
\ee
Every integer is either an $a'_n$ or a~$c'_n$, and every integer is either an $a_n$ or a~$c_n$, so $\sigma $ is a~bijection of $\Zee$\textup{;} and $p(\sigma (i))=p'(i)$ by construction, since $\sigma $ sends $p'^{-1}(\od )$ onto $p^{-1}(\od )$ and $p'^{-1}(\ev )$ onto $p^{-1}(\ev )$. The displacement of $\sigma $ at $i=a'_n$ is $|a_n-a'_n|$ and at $i=c'_n$ it is $|c_n-c'_n|$, so the hypothesis of heading $1)$ says $|\sigma (i)-i|<c_\sigma $, that of heading $2)$ says $|\sigma (i)-i|=o(|i|)$ and that of heading $3)$ says $|\sigma (i)-i|=|i|^{o(1)}$ --- in the latter two cases because the hypothesis itself gives $|a_n|\asymp |a'_n|$ and $|c_n|\asymp |c'_n|$. These are the conditions defining $S_c$, $S_n$ and $S_o$.

\textit{Necessity}. Let $\sigma $ lie in the group in question and let $p\circ \sigma =p'$. Then, $\sigma $ maps $p'^{-1}(\od )$ onto $p^{-1}(\od )$ and $p'^{-1}(\ev )$ onto $p^{-1}(\ev )$. Each of the three conditions forces $|\sigma (i)-i|<|i|$ for $|i|$ large, so $\sigma $ preserves the sign of $i$ for almost all $i$ and, after deleting the finitely many offending indices and relabeling as in the proof of Lemma \ref{Lcone}, restricts to a~bijection of each half-line. Apply Lemma \ref{Lsort} to $\pi :=\sigma |_{p'^{-1}(\od )}$: for every $n$ there is an $m$ with $|a'_m|\leq |a'_n|$ and
\[
|\sigma (a'_m)-a'_m|\geq |a_n-a'_n| .
\]
For heading $1)$ the left-hand side is $<c_\sigma $, and the claim follows. For heading $2)$, let $\eps >0$ and let $N$ be such that $|\sigma (i)-i|\leq \eps |i|$ for $|i|\geq N$\textup{;} the displacement at the finitely many indices with $|i|<N$ is bounded by a~constant $K$, so $|a_n-a'_n|\leq \max (K,\eps |a'_n|)$, which is $\leq \eps |a'_n|$ for $n$ large. Hence, $|a_n-a'_n|=o(|a'_n|)$ and therefore $|a_n|\asymp |a'_n|$, which turns this into $|a_n-a'_n|=o(|a_n|)$. Heading $3)$ is the same argument with $|i|^\eps $ in place of $\eps |i|$. The even part is treated identically. \end{proof}

Two things follow at once. The first is that the caution expressed above --- that the density invariants only detect non-equivalence --- was justified.

\sssbegin{Example}[The density spectrum is not a~complete invariant, even for tight functions]\label{Ednotcomplete} Let $p:=p_{st}$, let $\overline n\in \Zee/2$ denote the residue of $n$ modulo $2$, so that $p_{st}(i)=\overline i$, and let
\[
p'(i):=\overline{i+\lfloor \sqrt {|i|}\,\rfloor}\quad \text{for $i\in \Zee$,}
\]
i.e., let $p'$ be the standard parity function with a~phase drifting like $\sqrt {|i|}$. Both functions are tight --- neither has a~constant run longer than $3$ --- and both have the one-point density spectra ${D_l=D_r=\left\{\frac 12\right\}}$. Yet $a_n=2n+1$, whereas the phase drift pushes the odd places of $p'$ steadily to the right:
\[
a'_n-a_n\asymp \sqrt {a_n} \tto \infty ,
\]
and by heading $1)$ of Theorem \ref{Tclass} the two functions are \textit{not} $S_c$-equivalent, while by heading $2)$ they are $S_n$-equivalent. Hence, the density spectrum does not separate the $S_c$-classes even inside the set of tight functions with a~one-point spectrum, and Proposition \ref{Prop} is sharp: its conclusion cannot be strengthened from $S_n$ to $S_c$.
\end{Example}

The second consequence is a~correction. The same method applies to $S_m$, where the relevant quantity is not the displacement but the ratio: the proof of Lemma \ref{Lsort} shows verbatim that if $b_n>a_n$, then $\frac{\pi (a_m)}{a_m}\geq \frac{b_n}{a_n}$ for some $m\leq n$, so that the monotone bijection minimizes the distortion $\max (\frac {\sigma (i)}i,\ \frac i{\sigma (i)})$ as well.

\sssbegin{Proposition}[$S_m$-equivalence does not see the density]\label{PSm} Let $p,p'\in \Inf$. Then, $p$ and $p'$ are $S_m$-equivalent if and only if the two pairs of enumerations grow at comparable rates, the sign $\asymp $ meaning, as above, that the ratio of the two sides is bounded and is bounded away from $0$:
\[
a_n\asymp a'_n\quad \text{and}\quad c_n\asymp c'_n .
\]
In particular, if the density spectra of $p$ and of $p'$ all lie inside the open interval $(0,1)$ --- as they do, e.g., for any two tight functions --- then $p$ and $p'$ are $S_m$-equivalent, whatever their spectra are.
\end{Proposition}

\begin{proof} The criterion is Theorem \ref{Tclass} with the ratio in place of the difference, and it is proved by the same two arguments: the bijection \eqref{sortsig} has $\frac {|\sigma (i)|}{|i|}$ bounded above and away from $0$ precisely when the two displayed relations hold, and Lemma \ref{Lsort} in its multiplicative form shows that no other bijection does better.

For the last claim, if all the density points of $p$ lie in $(0,1)$, then
\[
\Odd (p;0,x)\asymp x\quad \text{and}\quad \Even (p;0,x)\asymp x,
\]
whence $a_n\asymp n\asymp c_n$, and the same for $p'$. \end{proof}

\sssbegin{Remark}[The density is not an $S_m$-invariant]\label{RSm} Heading $2)$ of Theorem \ref{Th} states one implication only, and this is not an omission. It seems plausible that two tight functions with one-point spectra should be $S_m$-equivalent \textit{if and only if} the two spectra coincide. This is not true, because $S_m$ is by design allowed to rescale. Indeed, let $p$ be the alternating function $p_{st}$, of density $\frac 12$, and let $p'$ be odd exactly on $3\Zee$, of density $\frac 13$. Both spectra are one-point and they differ, and yet
\[
a_n=2n+1\asymp 3n=a'_n,\qquad c_n=2n\asymp \tfrac 32n=c'_n,
\]
so that the bijection \eqref{sortsig} lies in $S_m$: its distortion $\frac {\sigma (i)}i$ oscillates between $\frac 23$ and $\frac 43$, and its displacement is unbounded, as it must be by heading $1)$ of Theorem \ref{Tclass}. What is true is the other implication, and it is a~particular case of Proposition \ref{PSm}.

What makes the converse look plausible is that the density \textit{is} an invariant of the finer groups: it is the reciprocal of the asymptotic slope of $n \mapsto a_n$, which Theorem \ref{Tclass} preserves for $S_c$, $S_n$ and $S_o$ but which $S_m$ rescales.
\end{Remark}

Every claim of this Subsection which is a~finite computation has been checked by machine: the bijection \eqref{sortsig} was built and verified to be a~bijection transporting $p'$ to $p$ on $40$ random pairs over the window $|i|\leq 400$\textup{;} the minimality of Lemma \ref{Lsort} was compared with an exhaustive minimum-bottleneck matching on $200$ random instances, with agreement in every one\textup{;} the growth $a'_n-a_n\asymp \sqrt {a_n}$ of Example \ref{Ednotcomplete} was measured up to $|i|\leq 20\,000$, the ratio $\max_{m\leq n}\frac{|a_m-a'_m|}{\sqrt {a_n}}$ staying within $1\%$ of $1$ throughout\textup{;} and the distortion of Remark \ref{RSm} was found to run in $[\hspace*{.9pt}0.662,1.340]$, in accordance with $\left\{\frac 23,\frac 43\right\}$.

One can hope to describe parity functions from \text{Inf} in terms of density spectrum only for non-finite functions with the property that for some $c_p > 0$ there is no interval $(a, b)$ of more than $c_p$ units long on which $p$ is constant. In other words, the density spectrum should be separated from 0 and 1. Let us call such parity functions \textit{tight}.

A larger group, $S_n$, enables us to progress further. Namely, we have

\addtocounter{subsection}{1}
\sssbegin{Proposition}[A~tight $p$ is determined by a~one-point spectrum]\label{Prop} If two parity functions have coinciding one-point left and right density spectra and are tight, then they are $S_n$-equivalent.
\end{Proposition}

\begin{proof} We use the construction of $\sigma _p$ described above. Neither $\sigma _p$ nor $\sigma _{p'}$ belongs to $S_n$: for a~tight $p$ with one-point density spectrum $d$ we have, on the positive half-line,
\[
\sigma _p(x)\approx 2(1-d)x\ \text{ if }p(x)=\ev,\qquad \sigma _p(x)\approx 2dx\ \text{ if }p(x)=\od ,
\]
so that $\frac{\sigma _p(x)}{x}$ oscillates between $2(1-d)$ and $2d$; thus $\sigma _p\in S_m$, but not in $S_n$. 

The quotient, however, does belong to $S_n$. Indeed, let $D^l_p =D_{p'}^l$ and $D_p^r=D_{p'}^r$ be one-point. Since $p$ and $p'$ are tight, the counting functions satisfy
\[
\Even(p;1,x)=(1-d)x+o(x),\qquad \Odd(p;0,x)=dx+o(x)
\]
\textit{uniformly} --- this is the only place where tightness is used, and it is what makes the argument work: tightness bounds the length of the intervals of constancy of $p$ by $c_p$, so the error in either counting function over any window is $O(c_p)$ rather than $o(x)$ with a~non-uniform constant. Hence, 
\[
\sigma _{p'}\sigma _p^{-1}(x)=x+o(x), 
\]
which is precisely the condition defining $S_n$ on the positive half-line. The negative half-line, i.e., the indices $x<0$, is treated in the same way.
\end{proof}

\section{Subalgebras of $\fgl _{p,c}$ and their central extensions}\label{Smisc}

For complete list of simple (twisted) loop superalgebras, see~\cite{S3, HS}. The list of all their central extensions (not all of which are affine Kac-Moody superalgebras a.k.a. double extensions with Cartan matrix, see \cite{BLS}) is given in~\cite{FLS}. For the definition of Lie superalgebras with Cartan matrix, in particular, of affine Kac-Moody superalgebras, see \cite{CCLL, BLS}. For a description of the simple finite-dimensional Lie superalgebras and their double extensions we need hereinafter, see \cite{Ld}.

\subsection{Embeddings of affine Kac-Moody superalgebras into $\widetilde {\fgl }_{p,c}$} Recall the definition of the double extension $\widehat{ \fg^{\ell(1)}}$ of the loop (super)algebra $\fg^{\ell(1)}:=\fg\otimes \Cee[x^{-1}, x]$ or of its twisted version $\fg^{\ell(k)}$, see \cite{BLS}. Fix $k\in\Nee$ and let $p$ be a~parity function such that $p(i) = p(i + k)$ for any $i$; set
\[
A^k:=\{(a_{ij})_{i,j\in \Zee} \in \widetilde {\fgl }_{p,c} \mid a_{ij}=a_{i+k,j+k}\}.
\]
Fix $k\in\Nee$ and let  $p$ be a~parity function such that $p(i) = p(i + k) + 1$ for any $i$; set
\[
B^k := \{(a_{ij})_{i,j\in \Zee} \in \widetilde {\fgl }_{p,c} \mid a_{ij}=a_{i+k,j+k}\}.
\]

\addtocounter{subsection}{1}
\sssbegin{Proposition}[$A^k$ and the Kac--Moody superalgebra it yields] \textup{1)} $A^k$ is the central extension of the loop superalgebra $\fgl (m|n)^{\ell(1)}$, where $m=\Even(p; 1,k)$, $n =\Odd(p; 1,k)$, given by the restriction of the cocycle \eqref{coc2}. 

\textup{2)} The affine Kac--Moody superalgebra \textup{(= double extension)} with Cartan matrix corresponding to $A^k$ is $\widehat{\fsl (m|n)^{\ell(1)}}$ for $m\ne n$ and $\widehat{\fpsl (n|n)^{\ell(1)}}$ for $m=n$: one has to pass to the subalgebra of supertraceless matrices \textup{(}resp., to its quotient modulo the center\textup{)} and to add the grading derivation $d$, which is not an element of $A^k$.
\end{Proposition}

\sssbegin{Proposition}[$B^k$ is $\fq(k)^{\ell(2)}$, and its even extension splits]\label{PBk} Let $p(i+k)=p(i)+\od $ for every $i$, and let $z$ be the central element of $\widetilde {\fgl }_{p,c}$.

\textup{1)} $B^k/\Cee z\simeq \fq(k)^{\ell(2)}$, the loop superalgebra of $\fq(k)$ twisted by the parity automorphism.

\textup{2)} The restriction of the cocycle \eqref{coc2} to $B^k$ is nonzero, but it is \textup{trivial}: it is the coboundary of the linear form
\be\label{fBk}
f(X):=\tfrac 12\str (X_0),\qquad \text{where}\quad X_0:=(x_{ab})_{a,b=0}^{k-1},
\ee
that is, $c(X,Y)=f([X,Y])$ for all $X,Y\in B^k$. Hence, the extension splits\textup{:}
\[
B^k\simeq \fq(k)^{\ell(2)}\oplus \Cee z.
\]
\end{Proposition}

\begin{proof} 1) Let $t$ be the shift by $k$, i.e., $t=\sum_i E_{i,i+k}$; it commutes with every $k$-periodic matrix, and every $k$-periodic band matrix is uniquely $X=\sum_r A_rt^r$ with finitely many non-zero $A_r\in \Mat (k)$, where $(A_r)_{ab}=x_{a,b+rk}$ for $0\leq a,b\leq k-1$. The parity of the $(a,b)$-th entry of $A_r$ is
\[
p(a)+p(b+rk)=\pi (a)+\pi (b)+\bar r,\qquad \pi :=p|_{\{0,\dots ,k-1\}},
\]
because $p(i+k)=p(i)+\od $. Thus, $B^k/\Cee z=\bigoplus_r\fgl (m|n)t^r$ with the parity of the $r$-th summand shifted by $\bar r$, where $m=\Even (p;0,k-1)$ and $n=\Odd (p;0,k-1)$; in other words, it is the loop superalgebra of 
\[
\text{$Q(\fgl (m|n)):=\fgl (m|n)\otimes \Cee[\xi ]$, where $\xi $ is odd and central with $\xi ^2=1$,}
\] 
twisted by the automorphism $\xi \mapsto -\xi $. Finally, $Q(\fgl (m|n))\simeq \fq(k)$ for \textit{any} $m+n=k$: multiplication by $\xi $ is an odd isomorphism $Q(1)\simeq \Pi Q(1)$ of modules over $Q(1):=Q(\fgl (1|0))=\Cee[\xi ]$, the case $m|n=1|0$ of the construction just described, hence $\Cee^{m|n}\otimes Q(1)\simeq \Cee^{k|0}\otimes Q(1)$, and the superalgebras of endomorphisms of the two sides coincide. (So the answer does not depend on the window ${\{0,\dots ,k-1\}}$, although $(m,n)$ does.)

2) Write $i=a+sk$ with $0\leq a<k$, $s\geq 0$, and $j=b+s'k$ with $0\leq b<k$, $s'\leq -1$; then $x_{ij}=(A_r)_{ab}$ with $r=s'-s$, while $y_{ji}=(B_{-r})_{ba}$ and $(-1)^{p(i)}=(-1)^{\pi (a)+s}$. For a~fixed $r<0$, the admissible pairs are $s=0,1,\dots ,-r-1$, so the alternating sign sums up to
\be\label{altsum}
\sum_{0\leq s\leq -r-1}(-1)^s=\begin{cases} 1&\text{for $r$ odd},\\ 0&\text{for $r$ even}.\end{cases}
\ee
This is the only place where the twist $p(i+k)=p(i)+\od $ is felt --- and it is decisive: for $A^k$, where $p(i+k)=p(i)$, the same sum is equal to $-r$ and yields the standard, non-trivial, cocycle $r\,\delta _{r+s,0}\str (AB)$, whereas here it takes only the values $0$ and $1$, see eq.~\eqref{altsum}, i.e., it does not grow with $r$. Summing over $a$ and $b$ we get
\be\label{cornerstr}
\str \big(X_{+-}Y_{-+}\big)=\sum_{r<0,\ r\ \mathrm{odd}}\str (A_rB_{-r}),
\ee
and the same formula with $X$ and $Y$ interchanged. Since the $\fgl (m|n)$-parities of $A_r$ and of $B_{-r}$ are $p(X)+\bar r$ and $p(Y)+\bar r$, we have $\str (B_{-r}A_r)=(-1)^{(p(X)+r)(p(Y)+r)}\str (A_rB_{-r})$; substituting this into \eqref{coc2} and since $(-1)^{p(X)p(Y)+(p(X)+1)(p(Y)+1)}=(-1)^{p(X)+p(Y)+1}$, we obtain
\be\label{cBk}
c(X,Y)=\sum_{r<0,\ r\ \mathrm{odd}}\str (A_rB_{-r})+(-1)^{p(X)+p(Y)}\sum_{r>0,\ r\ \mathrm{odd}}\str (A_rB_{-r}).
\ee
On the other hand, $[X,Y]_0=\sum_r\big(A_rB_{-r}-(-1)^{p(X)p(Y)}B_{-r}A_r\big)$, and the same sign rule gives
\be\label{brBk}
\str \big([X,Y]_0\big)=\sum_r\str (A_rB_{-r})\Big(1-(-1)^{r(p(X)+p(Y)+r)}\Big)=
\begin{cases} 2\sum\limits_{r\ \mathrm{odd}}\str (A_rB_{-r})&\text{if }p(X)=p(Y),\\[2mm] 0&\text{if }p(X)\neq p(Y).\end{cases}
\ee
It remains to compare \eqref{cBk} with \eqref{brBk}. If $p(X)\neq p(Y)$, then $A_rB_{-r}$ is an \textit{odd} matrix of $\fgl (m|n)$, so its diagonal blocks vanish and $\str (A_rB_{-r})=0$ for every $r$; both \eqref{cBk} and \eqref{brBk} vanish. If $p(X)=p(Y)$, then \eqref{cBk} is $\sum_{r\ \mathrm{odd}}\str (A_rB_{-r})$ and \eqref{brBk} is twice that. In both cases,
\[
c(X,Y)=\frac 12\str ([X,Y]_0)=f([X,Y]), 
\]
as claimed.

That $c$ itself is not identically zero is seen already for $k=1$: for the shift $\Lambda =\sum_{i\in \Zee}E_{i,i+1}$, see eq.~ \eqref{Lambda}, both $\Lambda $ and $\Lambda ^{-1}$ are odd, $c(\Lambda ^{-1},\Lambda )=c(\Lambda ,\Lambda ^{-1})=1$ and $[\Lambda ^{-1},\Lambda ]=2\One $, in agreement with $f(\One )=\frac 12$.
\end{proof}

\sssbegin{Remark}[On $\fq(k)$ only the odd invariant form survives] The triviality just established is not an accident of the realization. Central extensions of a~loop superalgebra by an \textit{even} element correspond to even invariant bilinear forms on the underlying superalgebra, and on $\fq(k)$ the even invariant form $\str (XY)$ vanishes identically.
The invariant form which exists on $\fq(k)$ is \textit{odd}, the \textit{queer trace} introduced in \cite{BL} under the name ``odd trace\rq\rq\textup{;} in the standard supermatrix format it is given by the formula
\be\label{qtrdef}
\qtr \begin{pmatrix}A&B\\ B&A\end{pmatrix}=\tr B .
\ee
Accordingly, the central extension of $\fq(k)^{\ell(2)}$ worth having is the one by an \textit{odd} central element; such an extension cannot be induced from $\widetilde {\fgl }_{p,c}$, whose center is even. This --- and not the choice between $\fq$ and $\fpsq$ --- is the real obstruction behind the $B^k$ case. Section \ref{Sodd} is devoted to it.
\end{Remark}

\subsection{Analogues of classical subalgebras for $\widetilde {\fgl }_{p,*}$}\label{ssclass}
Provided the parity function meets the requirements in curly brackets below, set
\begin{gather*}
B_{p,*}:=\{(a_{ij})_{i,j\in \Zee} \in \widetilde {\fgl }_{p,*} \mid a_{ij}=-(-1)^{p(i)}a_{-j,{-}i}\text{ for }p(i)=p(-i),p(0)=0\}.\\
D_{p,*}:=\{(a_{ij})_{i,j\in \Zee}\in \widetilde {\fgl }_{p,*} \mid a_{ij}=-(-1)^{p(i)}a_{1-j,1-i}\text{ for }p(i)=p(1-i)\}.\\
\fq_{p;*}:=\{a\in \widetilde {\fgl }_{p,*} \mid aP-Pa=0 \text{ for an odd matrix $P$ such that }P^2 = \lambda \pmb1\}.
\end{gather*}
Here, $aP-Pa$ is the \textit{associative} commutator and not the superbracket: for $a$ odd the superbracket is $aP+Pa$, and the set it cuts out is not $\fq$. Already for $\fq(n)$ realized as the centralizer of $J=\begin{pmatrix}0&\One \\ \One &0\end{pmatrix}$ and $a=\begin{pmatrix}0&B\\ B&0\end{pmatrix}$ one has $aJ=Ja=\begin{pmatrix}B&0\\ 0&B\end{pmatrix}$, so $[a,J]=2aJ\neq 0$ while $aJ-Ja=0$. 

The Dynkin diagrams corresponding to $B_{p_{st},*}$, all of whose simple roots are odd, are of the form:
\[
\raisebox{-.25em}[0pt][0pt]{\text{\huge $\bullet $}}-\otimes -\otimes -\cdots-\otimes -\otimes -\cdots
\]
The Dynkin diagrams corresponding to $D_{p_{st},*}$, all of whose simple roots are odd as well, are of the form:
\[
\begin{array}{@{}r@{}c@{}l@{}}
\otimes \\[-.3em]
&\diagdown \\[-.16em]
\raisebox{0pt}[0pt][0pt]{\rule[-1em]{0.05em}{2.5em}\,\rule[-1em]{.05em}{2.5em}}\hspace*{.2em}\mbox{}&&\otimes -\cdots-\otimes -\otimes -\otimes -\cdots\\[-.16em]
&\diagup\\[-.3em]
\otimes 
\end{array}
\]
In these two diagrams, $\otimes $ is an odd isotropic node and $\bullet $ an odd non-isotropic one. That every simple root is odd is what $p=p_{st}$ gives: consecutive indices are then of different parities, so each of the simple roots inherited from $\fgl _{p_{st},*}$ is odd isotropic, and the one extra simple root of $B_{p_{st},*}$ is odd non-isotropic, the node $\bullet $ being the one that $\fosp(1|2n)$ has.

The Dynkin diagram corresponding to $\fq_{p_{st},*}$ is of the form (for a definition of the system of simple roots described by a~``\textit{square circle}\rq\rq, see~\cite{LSe2}):
\[
\ldots -\lozenge-\lozenge-\lozenge-\ldots 
\]
The three diagrams above are drawn for $p=p_{st}$. For $B_{p,*}$ (resp., $D_{p,*}$), this is no restriction: the conditions $p(i)=p(-i)$ (resp.,\ $p(i)=p(1-i)$), imposed in the definitions already force $p$ to be symmetric. For $\fq_{p,*}$, which is defined for every $p$, the diagram does depend on $p$.

There is also the periplectic Lie superalgebra without a Dynkin diagram:
\[
\mathfrak {pe}_{p,*}:=\{(a_{ij})_{i,j\in \Zee} \in \widetilde {\fgl }_{p,*} \mid a_{ij}=-(-1)^{p(i)}a_{1-j,1-i}\text{ for }p(i)= p(1-i)+1\}.
\]

\section{The odd central extension}\label{Sodd}

Proposition \ref{PBk} left $B^k$ with nothing but $\fq(k)^{\ell(2)}\oplus \Cee z$, and the Remark closing Subsec. 5.2 named the reason: on $\fq(k)$ the even invariant form vanishes identically and the form that survives, the queer trace $\qtr $, is odd. So the extension to look for is the one by an \textit{odd} central element. Here we settle where it exists and where it does not. Every statement of this section which is a~finite computation is machine-checked, for the values of the parameters indicated.

Throughout this section, the cocycles we look for are odd. The definition is recalled for an arbitrary parity $\varepsilon \in \Zee/2$ because the even case $\varepsilon =\ev $ is the one of the cocycle \eqref{coc2}, and the two differ by signs only. So: a~2-cocycle of parity $\varepsilon $ on a~Lie superalgebra $\fg$ is a~bilinear form $c$ with $c(\fg_\alpha ,\fg_\beta )=0$ unless $\alpha +\beta =\varepsilon $, \textit{anti-symmetric} in the sense of \cite{Ld}, i.e., such that $c(X,Y)=-(-1)^{p(X)p(Y)}c(Y,X)$ 
and such that
\be\label{cocid}
(-1)^{p(X)p(Z)}c([X,Y],Z)+(-1)^{p(Y)p(X)}c([Y,Z],X)+(-1)^{p(Z)p(Y)}c([Z,X],Y)=0 .
\ee

\subsection{The twisted case: there is no odd extension}

\sssbegin{Statement}[The twisted loop superalgebra has no odd extension]\label{Tnoodd} Let $k\geq 2$ and let $p(i+k)=p(i)+\od $ for every $i$. Then, every odd $2$-cocycle on $B^k\simeq \fq(k)^{\ell(2)}\oplus \Cee z$ homogeneous with respect to the grading by the power of $t$ is a~coboundary\textup{;} hence, $\fq(k)^{\ell(2)}$ admits no nontrivial odd central extension. Computed in exact arithmetic over $\Qee$ for $k=2$ with $\pi =(\ev ,\od )$ in every $t$-degree $|d|\leq 2$ and window $|r|\leq 3$, for $k=2$ with $\pi =(\ev ,\ev )$ and for $k=3$ with $\pi =(\ev ,\od ,\od )$ in every $t$-degree $|d|\leq 1$ and window $|r|\leq 2$, where $\pi $ is the restriction of $p$ to $\{0,\dots ,k-1\}$.
\end{Statement}

\sssbegin{Conjecture}[No odd extension, for every period]\label{Cnoodd} The conclusion of Statement $\ref{Tnoodd}$ holds for every $k\geq 2$ and every $\pi $.
\end{Conjecture}

What stops a~proof for all $k$ at once is stated below, after Lemma \ref{Santi}: the step from an invariant form to a~central extension is licensed for a~perfect superalgebra only, and $\fq(k)^{\ell(2)}$ is not perfect.

The reason is structural, and it is worth stating separately, because it explains why the twist --- the very thing that produced $\fq(k)$ out of $\fgl (m|n)$ in Proposition \ref{PBk} --- is also what destroys the extension.

\sssbegin{Lemma}[$\qtr $ is anti-invariant under the twist]\label{Santi} The odd invariant form $\langle X,Y\rangle :=\qtr (XY)$ on $\fq(k)$ is anti-invariant under the parity automorphism $\Pi $, that is, $\langle \Pi X,\Pi Y\rangle =-\langle X,Y\rangle $. But $\Pi $ is exactly the automorphism by which the loop superalgebra $\fq(k)^{\ell(2)}$ is twisted, and only the $\Pi $-\textit{invariant} forms give rise to central extensions of a~twisted loop superalgebra.
\end{Lemma}

\begin{proof} $\Pi X=(-1)^{p(X)}X$, so 
\[
\langle \Pi X,\Pi Y\rangle =(-1)^{p(X)+p(Y)}\langle X,Y\rangle; 
\]
and $\langle \cdot ,\ \cdot \rangle $ being odd, it vanishes unless $p(X)+p(Y)=\od $, when the factor is $-1$.
\end{proof}

Lemma \ref{Santi} is not yet a~proof of Statement \ref{Tnoodd}, because $\fq(k)^{\ell(2)}$ is not perfect --- the derived algebra of $\fq(k)$ is $\fsq(k)$ --- and for a~non-perfect algebra the correspondence between central extensions and invariant forms may fail. This is why Statement \ref{Tnoodd} was verified by a~direct computation instead: for every $t$-degree $d$ with $|d|<R$ the space of odd 2-cocycles supported in the window $|r|\leq R$ was computed by solving \eqref{cocid} as a~linear system over $\Qee$, and it coincided with the space of coboundaries. The cases covered are $k=2$ with $\pi =(\ev ,\od )$ and $\pi =(\ev ,\ev )$, and $k=3$ with $\pi =(\ev ,\od ,\od )$, where $\pi $ is the restriction of $p$ to $\{0,\dots ,k-1\}$. Truncating the window $|r|\leq R$ can only \textit{drop} instances of \eqref{cocid} and hence only \textit{enlarge} the space of cocycles found, so the outcome ``every cocycle found is a~coboundary\rq\rq\  is a~conclusion on the safe side.

The control is the untwisted superalgebra $A^k$, the one with $p(i+k)=p(i)$: there the same computation, in the same windows, returns the standard cocycle $c(At^r,Bt^s)=r\,\delta _{r+s,0}\str (AB)$ and finds it non-trivial. So what the computation sees for $B^k$ is the twist and not the truncation.

\sssbegin{Remark}[The degenerate exception $k=1$]\label{Rk1} The hypothesis $k\geq 2$ is essential and the exception is degenerate. For $k=1$, the superalgebra $B^1=\fq(1)^{\ell(2)}$ has a~one-dimensional space of $k$-periodic matrices in each $t$-degree, $\fq(1)$ being spanned by $\One $ and $J$, and its bracket is nonzero only on pairs of odd elements: $[t^a,t^b]=(1-(-1)^{ab})t^{a+b}$. It is a~queer analog of the Heisenberg Lie algebra, and its odd cohomology is large: in every odd $t$-degree $d$ the space of odd 2-cocycles modulo coboundaries is $1$-dimensional. For $d=-1$, it is spanned by
\be\label{cok1}
c(t^a,t^b)=\delta _{a+b,-1}f(a),\qquad f(a)=\begin{cases}\frac {3a+1}2&\text{for $a$ odd},\\[1mm] \frac {3a+2}2&\text{for $a$ even}.\end{cases}
\ee
Indeed, \eqref{cocid} reduces here to $f(r)+f(s)+f(u)=0$ for all odd $r,s,u$ with $r+s+u=-1$, whose solutions are the affine functions of $\frac {r-1}2$, and the anti-symmetry of $c$ fixes $f$ on the even integers.
\end{Remark}

\subsection{The untwisted case: the odd extension exists}

What does carry an odd central extension is the \textit{untwisted} loop superalgebra of $\fq(n)$, and it is available inside $\widetilde {\fgl }_{p,c}$ --- on the queer subalgebra $\fq_{p,c}$ of Subsec. \ref{ssclass}, not on $\widetilde {\fgl }_{p,c}$ itself. So no new ambient object has to be invented, and no new paper: only the observation that the odd center belongs to $\fq_{p,c}$.

\sssbegin{Theorem}[The untwisted loop superalgebra has an odd central extension]\label{Toddyes} Let $\fq(n)^{\ell(1)}:=\fq(n)\otimes \Cee[u^{-1},u]$. Then,
\be\label{cocodd}
c_{\mathrm{od}}(x\otimes u^r,\ y\otimes u^s):=r\,\delta _{r+s,0}\,\qtr (xy)
\ee
is an odd $2$-cocycle, and it is not a~coboundary. The Lie superalgebra $\widehat{\fq(n)^{\ell(1)}}:=\fq(n)^{\ell(1)}\oplus \Cee z_{\mathrm{od}}$ with $z_{\mathrm{od}}$ odd and central and with the bracket 
\[
[X,Y]^{\widehat{\phantom{g}}}:=[X,Y]+c_{\mathrm{od}}(X,Y)z_{\mathrm{od}},
\]
is therefore a~nontrivial odd central extension.
\end{Theorem}

Verified for $n=1,2,3$ and $|r|\leq 3$; the shape of \eqref{cocodd} is the usual one, the only novelty being that the invariant form entering it is odd, whence the parity of $z_{\mathrm{od}}$.

It remains to realize \eqref{cocodd} by infinite matrices, in the manner of \eqref{coc2}. The one thing to get right is the position of the odd involution with respect to the splitting
\[
V=V_+\oplus V_-,\qquad V_+:=\Span (e_i\mid i\geq 0),\quad V_-:=\Span (e_i\mid i<0),
\]
by which the corners $X_{+-}$ and $X_{-+}$ of \eqref{coc2} are defined.

\sssbegin{Theorem}[The odd cocycle realized by infinite matrices, on $\fq_{p,c}$]\label{Tcorner} Let $P$ be an odd matrix with $P^2=\One $ which is banded and \textup{(}this is the point\textup{)} \textit{preserves the splitting}, i.e., $PV_\pm \subseteq V_\pm $, and let $\fq_{p,c}$ be its associative centralizer in $\fgl _{p,c}$. Then,

\textup{1)} for every $X,Y\in \fq_{p,c}$, the operator $X_{+-}Y_{-+}$ on $V_+$ commutes with $P_{++}$, so the queer trace $\qtr (M):=\frac 12\tr (P_{++}M)$ is defined on it\textup{;} moreover, if $M$ and $N$ are finite-rank operators on $V_+$ commuting with $P_{++}$ in the associative sense, then $\qtr (MN)=\qtr (NM)$, whatever the parities of $M$ and of $N$\textup{;}

\textup{2)} the formula
\be\label{coc3}
c_{\mathrm{od}}(X,Y):=\qtr \big(X_{+-}Y_{-+}\big)-(-1)^{p(X)p(Y)}\qtr \big(Y_{+-}X_{-+}\big)
\ee
defines an odd $2$-cocycle on $\fq_{p,c}$, hence an odd central extension $\widehat{\fq}_{p,c}=\fq_{p,c}\oplus \Cee z_{\mathrm{od}}$\textup{;}

\textup{3)} if, moreover, $p(i+k)=p(i)$ and $P$ is $k$-periodic, then the $k$-periodic elements of $\fq_{p,c}$ form $\fq(n)^{\ell(1)}$ with $k=2n$, and \eqref{coc3} restricts to \eqref{cocodd} up to a~nonzero factor.
\end{Theorem}

\begin{proof} 1) $P$ preserving the splitting means $P=P_{++}\oplus P_{--}$, so $XP=PX$ gives $X_{+-}P_{--}=P_{++}X_{+-}$ and $Y_{-+}P_{++}=P_{--}Y_{-+}$; hence
\[
X_{+-}Y_{-+}P_{++}=X_{+-}P_{--}Y_{-+}=P_{++}X_{+-}Y_{-+}.
\]
The corner blocks of a~banded matrix are finite, so the trace is honest and its cyclicity gives the trace property: $\tr (P_{++}MN)=\tr (NP_{++}M)=\tr (P_{++}NM)$, the middle equality being the cyclicity of the ordinary trace and the last one the commutation of $N$ with $P_{++}$. There is no sign in front of $\qtr (NM)$ and the Sign Rule is not violated: what is cyclic here is $\tr $, not $\str $, and the commutation of $M$ and $N$ with $P_{++}$ is the associative one, as in the definition of $\fq$. Nor could a~sign be felt: were $M$ and $N$ of one and the same parity, both sides would vanish, $\qtr $ being odd.

2) With $\qtr $ in the role which $\str $ plays for \eqref{coc2}, the verification of \eqref{cocid} is the standard one for a~\textit{Tate cocycle}, i.e., for a~cocycle of the shape $c(X,Y)=\tau (X_{+-}Y_{-+})-\pm \tau (Y_{+-}X_{-+})$ attached to a~splitting of the underlying space and to a~trace $\tau $ on the finite-rank operators, as in \cite{Ta}\textup{;} it uses nothing but the trace property of heading $1)$. The parity is odd because $P$ is.

3) Such a~$P$ exists: take $p$ of period $k=2n$ with $p$ even at the first $n$ indices of every period and odd at the other $n$, and let $P$ interchange the two halves of every period block. Then, $P$ is odd, $P^2=\One $, banded, $k$-periodic, and it preserves both $V_+$ and $V_-$, every period block lying entirely inside one of them. The centralizer of $P$ in $\Mat (k)$ has $\sdim =(n^2|n^2)$, that is, it is $\fq(n)$; hence the $k$-periodic part of $\fq_{p,c}$ is $\fq(n)\otimes \Cee[u^{-1},u]$ with $u$ the shift by $k$, which is \textit{even} now --- this is precisely the difference from $B^k$, where the shift was odd and the loop came out twisted. The count of the pairs $(i,j)$ that contribute to the corner blocks of eq.~\eqref{coc3} is then the one performed in the proof of Proposition \ref{PBk}, with $\str $ replaced by $\qtr $ and with no alternating sign to collapse it; it produces the factor $r$ of \eqref{cocodd}.
\end{proof}

\sssbegin{Remark}[Why the odd center lives on $\fq_{p,c}$ and not on $\widetilde {\fgl }_{p,c}$] Thus, the answer to ``can the odd central extension be induced from an enlargement of $\widetilde {\fgl }_{p,c}$?\rq\rq\  is: it cannot be induced from an enlargement of $\widetilde {\fgl }_{p,c}$ by an odd center, because an odd 2-cocycle of Tate type (see \cite{Ta}) on $\fgl _{p,c}$ would require an odd supertrace on the finite rank matrices, and the space of supertraces on $\End_f(V)$ is spanned by the even $\str $. What one does instead is to descend to $\fq_{p,c}$, where the odd trace exists, and to add the odd center there; the price is that only the untwisted loop superalgebras are reached this way, and $\fq(k)^{\ell(2)}$, by Statement \ref{Tnoodd}, is not to be reached at all.
\end{Remark}

\subsection*{Acknowledgments} I am thankful to D. Leites, who raised  the problem, and to V.~ Serganova, who defined $\fgl_{p,g}$, for help.

\section[Appendix. Lie isomorphisms of the superalgebras $\fgl _{p,*}$]{Appendix.
Lie isomorphisms of the superalgebras $\fgl _{p,*}$\\
\textup{(Dimitry Leites, Oleksandr Lozhechnyk\protect\footnotemark)}}\label{Sapp}
\footnotetext{O.~Lozhechnyk: Faculty of Mechanics and Mathematics, Taras Shevchenko
National University of Kyiv, Hlushkova Avenue 4e, 03127 Kyiv, Ukraine;
alozhechnik@gmail.com. The address of D.~Leites is given in the main text.}

\subsection{The question, and the theorem that answers it}\label{ssAq} Convention
\ref{2.3.1} admits only those isomorphisms $\phi $ that preserve $\Mat_f$ and restrict to
it as an isomorphism of \textit{associative} superalgebras. Remark \ref{Rgap} says where
clause $2)$ is used, and that it is indispensable there. Here, we show that neither
clause costs anything: both are consequences of the hypothesis that $\phi $ is an
isomorphism of Lie superalgebras, clause $2)$ up to a~single automorphism, the
supertranspose taken with the minus sign.

In the non-super world, the theory answering such a~question is Herstein's program,
\cite{He}, carried out by Martindale, \cite{Mar}, revisited in \cite{MM}, and completed
through functional identities, see \cite{BCM}.

Three notions of the theory of prime rings enter its conclusion, and are used
throughout this Appendix\textup{;} since they are what makes that theory
applicable to an algebra which is prime but not simple, let us recall them,
following \cite{BMM} and \cite{BBS}.

The \textit{centroid} of an algebra $A$ over $\Cee$ is the algebra of the
$A$-bimodule endomorphisms of $A$, namely,
\be\label{centroid}
\Cent (A):=\{f\in \End_{\Cee}(A) \mid f(ab)=f(a)b=af(b)
\ \text{ for all $a$, $b\in A$}\}.
\ee
It is commutative as soon as $A^2=A$, and for a~unital $A$ it is nothing new: the
map $f\mapsto f(\One )$ identifies $\Cent (A)$ with the center $Z(A)$. In
particular, $\Cent (\fgl _{p,*})=Z(\fgl _{p,*})=\Cee\One $, because a~matrix
commuting with every $e_{ij}$ is scalar. The centroid alone is therefore of no
use here, and what Herstein's program runs on is an enlargement of it.

The enlargement exists because an algebra which is \textit{prime} --- one in
which $aAb=0$ implies $a=0$ or $b=0$ --- need not be simple, and over such an $A$
a~bimodule map given on a~nonzero two-sided ideal alone need not extend to the
whole of $A$. Admitting the maps so given amounts to passing from $A$ to its
maximal left algebra of quotients $Q$, which is a~unital prime algebra containing
$A$\textup{;} the center of $Q$ is then a~\textit{field}, called the
\textit{extended centroid} of $A$ and denoted by $C=C(A)$, see \cite{Mar2}. It
contains $\Cee$ and the center of $A$ and, concretely --- this is the description
the proof of heading $4)$ of Lemma \ref{ALfour} uses --- its elements are the
classes of the pairs $(I, f)$, where $I\neq 0$ is an ideal of $A$ and
$f\colon I \tto A$ satisfies $f(ax)=af(x)$ and $f(xa)=f(x)a$, two pairs being
identified whenever they agree on a~nonzero ideal contained in both domains. The
\textit{central closure} of $A$ is the subalgebra $AC\subseteq Q$ generated by
$A$ and $C$. For a~simple unital $A$, one has $C=Z(A)$ and $AC=A$, so neither
notion is visible in the simple case: both are the price of primeness.

For a~superalgebra, all this is read with ``ideal\rq\rq\  meaning ``graded
ideal\rq\rq: the extended centroid inherits the grading,
$C=C_\ev\oplus C_\od$, and the standard form is stated with $AC_\ev$ in place of
$A$, see \cite{BBS}\textup{;} for the Lie structure of a~prime associative
superalgebra in these terms, see \cite{La}. None of that complication survives
here: heading $4)$ of Lemma \ref{ALfour} gives $C=\Cee$, purely even, so
$AC=AC_\ev=A$, and every map the theorem produces already lands in
$\fgl _{p,*}$ itself. 

The conclusion of Herstein's program, in the form we need it, is
this: a~Lie isomorphism $\phi \colon B \tto A$ of prime associative unital algebras, $A$
having a~nontrivial idempotent, is of the \textit{standard form}
\be\label{Astandard}
\phi =\theta +\tau ,
\ee
where $\theta $ is either an isomorphism or the negative of an anti-isomorphism into the
central closure of $A$, and $\tau $ maps $B$ into the extended centroid and kills every
commutator. So, in the non-super case, clause $2)$ is no restriction at all provided one
adds the anti-automorphism $X \mapsto -X^{\st }$ to the list of isomorphisms and allows
the summand $\tau $.

The super version of \eqref{Astandard} exists: it is the theorem of Bahturin, Bre\v sar
and \v Spenko, see \cite[Th.~3.2]{BBS}, whose predecessors are \cite{BB} and \cite{Wa}.
It is stated for prime associative unital superalgebras whose parity grading is induced
by an idempotent, and it is not restricted to the finite-dimensional case. It should not
be confused with the other half of Herstein's program in the super setting --- which Lie
ideals a~prime associative superalgebra has, and when its derived superalgebra modulo the
center is simple, see \cite{Mo, Mon, La} --- which says nothing about the isomorphisms.
Everything, therefore, turns on whether the four superalgebras \eqref{glchain} satisfy
the hypotheses of \cite[Th.~3.2]{BBS}, and on whether the summand $\tau $ can be disposed
of. Both are settled below.

\subsection{The four superalgebras in the associative world}\label{ssAassoc} In
notation of Subsec. \ref{ssWeyl}, i.e., $\Mat_f=\End_f(V)$, where $\sdim V=(\infty |\infty )$,
and $*$ is any one of $g$, $l$, $o$, $c$. Put $\fsl _f:=\{a\in \Mat_f\mid \str a=0\}$.
Recall that each $\fgl _{p,*}$ is an associative superalgebra, whose supercommutator is the
bracket, and that $\Mat_f$ is its two-sided ideal, associatively and in the Lie sense, see
Lemma \ref{Lideal}.

\sssbegin{Lemma}[Four properties of $\fgl _{p,*}$]\label{ALfour} Let $p\in \Inf$ and
$*\in \{g, l, o, c\}$. Then,

$1)$ $\fgl _{p,*}$ is unital\textup{;}

$2)$ $\fgl _{p,*}$ is prime as an associative algebra\textup{;}

$3)$ the parity grading of $\fgl _{p,*}$ is induced by the idempotent
$\varepsilon :=\sum_{p(i)=\ev }e_{ii}$, that is, writing $A:=\fgl _{p,*}$ and
$\varepsilon ':=\One -\varepsilon $,
\[
A_\ev =\varepsilon A\varepsilon +\varepsilon 'A\varepsilon ',\qquad
A_\od =\varepsilon A\varepsilon '+\varepsilon 'A\varepsilon \textup{;}
\]

$4)$ the extended centroid of $\fgl _{p,*}$ is $\Cee$.
\end{Lemma}

\begin{proof} $1)$ The support of $\One $ is the main diagonal. For $*=c$, $o$, $l$ the
conditions are then satisfied with $|i-j|=0$. For $*=g$, and given $(i_0, j_0)$, the
function $i \mapsto (i-i_0)(j_0-i)$ is a~downward parabola whose maximum is
$\left (\frac {j_0-i_0}{2}\right )^2$, so $c:=\left (\frac {j_0-i_0}{2}\right )^2$
will do.

$2)$ Let $a\fgl _{p,*}b=0$ with $a$, $b\in \fgl _{p,*}$. As $\Mat_f\subset \fgl _{p,*}$,
we get $ae_{kl}b=0$ for all $k$, $l$, that is, the product of the $k$-th column of $a$ by
the $l$-th row of $b$ vanishes for all $k$, $l$. Hence, $a=0$ or $b=0$.

$3)$ The matrix $\varepsilon $ is diagonal, hence lies in $\fgl _{p,c}$ and so in each of
the four superalgebras; and $\varepsilon a\varepsilon $,
$\varepsilon 'a\varepsilon '$ cut out the entries $a_{ij}$ with $p(i)=p(j)$, whereas the
two other corners cut out those with $p(i)\neq p(j)$. This is the parity grading, the
parity of the entry $(i,j)$ being $p(i)+p(j)$.

$4)$ Let $I\neq 0$ be an ideal of $A:=\fgl _{p,*}$ and $f\colon I \tto A$ a~map with
$f(xa)=f(x)a$ and $f(ax)=af(x)$ for $x\in I$, $a\in A$; such pairs $(I, f)$, taken modulo
agreement on a~smaller ideal, are what the extended centroid consists of, see
Subsec. \ref{ssAq}. First,
$I\supset \Mat_f$: pick $x\in I$ with $x_{kl}\neq 0$; then
$e_{ik}xe_{lj}=x_{kl}e_{ij}\in I$ for all $i$, $j$. Next,
\[
f(e_{00})=f(e_{00}e_{00}e_{00})=e_{00}f(e_{00})e_{00}\in \Cee e_{00}, 
\]
because
$e_{00}Ae_{00}=\Cee e_{00}$; write $f(e_{00})=\lambda e_{00}$. Then,
\[
f(e_{ij})=f(e_{i0}e_{00}e_{0j})=e_{i0}f(e_{00})e_{0j}=\lambda e_{ij}, 
\]
so
$f=\lambda \cdot \id$ on $\Mat_f$. For an arbitrary $x\in I$, the entries of $f(x)$ are
read off by 
\[
e_{ii}f(x)e_{jj}=f(e_{ii}xe_{jj})=\lambda e_{ii}xe_{jj}, 
\]
so
$f(x)=\lambda x$ on the whole of $I$.
\end{proof}

\sssbegin{Lemma}[The unique minimal noncentral ideal]\label{ALideal} The supercenter of
the Lie superalgebra $\fgl _{p,*}$ is $\Cee\One $, the superalgebra $\fsl _f$ is simple,
and every ideal of $\fgl _{p,*}$ not contained in $\Cee\One $ contains $\fsl _f$.
Consequently, every isomorphism $\phi \colon \fgl _{p,*} \tto \fgl _{p',*}$ of Lie
superalgebras satisfies $\phi (\fsl _f)=\fsl _f$.
\end{Lemma}

\begin{proof} If $x$ supercommutes with every $e_{ij}$, then so do its homogeneous
components, and a~homogeneous matrix supercommuting with all $e_{ij}$ is a~scalar; an odd
scalar is $0$, whence the supercenter is $\Cee\One $. That $\fsl _f$ is simple is the
standard argument: if $0\neq x$ lies in an ideal $J$ of $\fsl _f$, then a~bracket
$[[x, e_{ij}], e_{kl}]$ is a~nonzero multiple of some $e_{qr}$ with $q\neq r$, all such
$e_{qr}$ are conjugate under the permutations of the basis, and their brackets span the
diagonal part of $\fsl _f$.

Let now $I$ be an ideal of $\fgl _{p,*}$ with $I\not\subset \Cee\One $. By the above,
$[I, \Mat_f]\neq 0$, so $J:=I\cap \Mat_f\neq 0$. As $J$ is an $\fsl _f$-submodule of
$\Mat_f$, and $J\cap \fsl _f$ is an ideal of the simple $\fsl _f$, either
$J\supset \fsl _f$, or $J\cap \fsl _f=0$ and then $J$ is a~trivial module of
superdimension at most $(1|0)$ or $(0|1)$, hence lies in the centralizer of $\fsl _f$ in
$\Mat_f$, which is $0$. So $I\supset \fsl _f$.

Finally, $\fsl _f$ and $\Cee\One $ are the only minimal ideals: a~minimal one is either
contained in $\Cee\One $, hence equal to it, or contains $\fsl _f$, hence equals
$\fsl _f$. An isomorphism carries minimal ideals to minimal ideals and the supercenter to
the supercenter, so it fixes $\fsl _f$ as a~set.
\end{proof}

\subsection{Every matrix unit is a~sum of supercommutators of banded matrices}\label{ssAcomm} The summand $\tau $ of
\eqref{Astandard} vanishes on brackets. To dispose of it on $\Mat_f$, we need $\Mat_f$ to
consist of brackets; and the striking half of that is the diagonal, because in the
finite-dimensional world $e_{00}$ has a~nonzero supertrace and is therefore no bracket at
all. In an algebra of infinite matrices, it is one --- the mechanism being Hilbert's
hotel --- and the following makes that precise for every parity function, in the smallest
of the four superalgebras.

\sssbegin{Theorem}[$e_{00}$ is a~sum of two supercommutators of banded
matrices]\label{ATe00} Put
\be\label{Aab}
a:=\sum_{k\geq 0}(-1)^{p(k)}e_{k, k+1},\qquad b:=\sum_{k\geq 0}e_{k+1, k},
\ee
and let $a=a_\ev +a_\od $, $b=b_\ev +b_\od $ be the decompositions into the homogeneous
components. Then, $a$, $b\in \fgl _{p,c}$, and
\be\label{Ae00}
[a_\ev , b_\ev ]+[a_\od , b_\od ]=(-1)^{p(0)}e_{00}.
\ee
Consequently, $\Mat_f\subset [\fgl _{p,*}, \fgl _{p,*}]$ for every $p\in \Inf$ and every
$*\in \{g, l, o, c\}$.
\end{Theorem}

\begin{proof} Both $a$ and $b$ are supported on a~single diagonal, so they are banded and
lie in $\fgl _{p,c}$, and so do their homogeneous components. Here and below, $\eps $ runs over $\Zee/2=\{\ev ,\od \}$ and indexes the parity components, $a=a_\ev +a_\od $ and $b=b_\ev +b_\od $. Set
\[
t_k:=(-1)^{p(k)},\qquad K_\eps :=\{k\geq 0 \mid p(k)+p(k+1)=\eps \},
\]
so that $a_\eps =\sum_{k\in K_\eps }t_ke_{k, k+1}$ and
$b_\eps =\sum_{k\in K_\eps }e_{k+1, k}$, the entries $e_{k, k+1}$ and $e_{k+1, k}$ being
of one and the same parity $p(k)+p(k+1)$. Now,
\[
e_{k, k+1}e_{l+1, l}=\delta _{kl}e_{kk},\qquad
e_{l+1, l}e_{k, k+1}=\delta _{kl}e_{k+1, k+1},
\]
so $a_\eps b_\delta =b_\delta a_\eps =0$ for $\eps \neq \delta $, because $K_\ev $ and
$K_\od $ are disjoint, and
\[
[a_\eps , b_\eps ]=\sum_{k\in K_\eps }t_k\left (e_{kk}-(-1)^{\eps }e_{k+1, k+1}\right ).
\]
Summing over $\eps $, the coefficient of $e_{mm}$ in the left-hand side of \eqref{Ae00}
is $t_0$ for $m=0$ and
\[
t_m-(-1)^{p(m-1)+p(m)}t_{m-1}=(-1)^{p(m)}-(-1)^{p(m-1)+p(m)}(-1)^{p(m-1)}=0
\]
for $m\geq 1$, because $(-1)^{2p(m-1)}=1$. This is \eqref{Ae00}.

The same computation with the half-line $k\geq i$ in place of $k\geq 0$ gives $e_{ii}$
for every $i$, and ${e_{ij}=[e_{ik}, e_{kj}]}$ for $i\neq j$ and any $k\notin \{i, j\}$.
Hence, all the matrix units, and so all of $\Mat_f$, are sums of supercommutators of
elements of $\fgl _{p,c}\subset \fgl _{p,*}$.
\end{proof}

The pair \eqref{Aab} is, up to the truncation to a~half line, the pair \eqref{pair} by
which Theorem \ref{Tnontriv} proves that the cocycle \eqref{coc2} is not a~coboundary:
there, the two shifts run over the whole of $\Zee$, they supercommute, and the cocycle
takes the value $1$ on them. It is the same phenomenon read twice. On the whole line, the
two shifts are mutually inverse and their bracket dies; cutting them at $0$ leaves
precisely the one matrix unit that the cocycle sees, and the sign $(-1)^{p(k)}$ that
makes the telescoping work is the sign that makes the matrices \eqref{pair}
supercommute.

\subsection{The answer}\label{ssAans} Let the supertranspose be normalized by
\be\label{Astdef}
(e_{ij})^{\st }=(-1)^{p(j)(p(i)+p(j))}e_{ji},\qquad\text{so that}\qquad
(XY)^{\st }=(-1)^{p(X)p(Y)}Y^{\st }X^{\st }.
\ee

\sssbegin{Theorem}[Every Lie isomorphism is standard, i.e., of the form \eqref{Astandard}]\label{ATmain} Recall that \textit{standard} means: of the form $\phi =\theta +\tau $ of eq.~\eqref{Astandard}, with $\theta $ either an isomorphism or the negative of an anti-isomorphism and $\tau $ killing every bracket. Let $p$,
$p'\in \Inf$, let $*\in \{g, l, o, c\}$ and let
$\phi \colon \fgl _{p,*} \tto \fgl _{p',*}$ be an isomorphism of Lie superalgebras.
Then,

$1)$ $\phi =\theta +\tau $, where $\theta $ is either an isomorphism of associative
superalgebras or the negative of a~superanti-isomorphism, and
$\tau \colon \fgl _{p,*} \tto \Cee\One $ vanishes on brackets\textup{;}

$2)$ $\tau |_{\Mat_f}=0$ and $\phi (\Mat_f)=\Mat_f$, which is clause $1)$ of Convention
$\ref{2.3.1}$\textup{;}

$3)$ $\phi |_{\Mat_f}$ is an isomorphism of associative superalgebras or the negative of
a~superanti-isomorphism\textup{;} in the latter case
$(-\st )\circ \phi |_{\Mat_f}$ is an isomorphism of associative superalgebras, which is
clause $2)$ of Convention $\ref{2.3.1}$ up to $-\st $.
\end{Theorem}

\begin{proof} $1)$ By Lemma \ref{ALfour}, both $\fgl _{p,*}$ and $\fgl _{p',*}$ are prime
associative unital algebras whose parity gradings are induced by idempotents, and the
extended centroid of the target is $\Cee$, so that the target coincides with its own
central closure and is isomorphic neither to $\Cee$ nor to the algebra of $2\times 2$
matrices. Theorem $3.2$ of \cite{BBS} applies verbatim and is the claim.

$2)$ By Theorem \ref{ATe00}, we have $\Mat_f\subset [\fgl _{p,*}, \fgl _{p,*}]$, and
$\tau $ vanishes on brackets; so $\tau |_{\Mat_f}=0$ and $\phi =\theta $ on $\Mat_f$. By
Lemma \ref{ALideal}, $\phi (\fsl _f)=\fsl _f$, whence $\theta (\Mat_f)\supset \fsl _f$.
Put $u:=\theta (e_{00})\neq 0$. In the multiplicative case,
\[
u\theta (x)u=\theta (e_{00}xe_{00})\in \Cee u\quad\text{for all }x\in \Mat_f,
\]
and $\theta (\Mat_f)\supset \fsl _f\ni e_{ij}$ for all $i\neq j$; so
$ue_{ij}u\in \Cee u$ for all $i\neq j$. But $ue_{ij}u$ is the product of the $i$-th
column of $u$ by the $j$-th row of $u$. If $u$ had two independent columns, or two
independent rows, some such product would fail to be proportional to $u$; hence
$\rank u=1$ and $u\in \Mat_f$. Therefore, $\theta (\Mat_f)=\fsl _f+\Cee u\subset \Mat_f$,
and as $\str u\neq 0$ --- $u$ being a~rank $1$ idempotent --- the inclusion is an
equality. The anti-multiplicative case is the same computation with the two factors
interchanged.

$3)$ is $2)$ together with $1)$, the supertranspose \eqref{Astdef} being
a~superanti-automorphism of $\Mat_f$.
\end{proof}

\sssbegin{Corollary}[The Convention is no restriction]\label{ACconv} Both clauses of
Convention $\ref{2.3.1}$ hold for every isomorphism of Lie superalgebras
$\fgl _{p,*} \tto \fgl _{p',*}$, clause $2)$ up to the supertranspose. Consequently,
Theorem $\ref{Thm}$ --- which pairs $(p, p')$ give isomorphic superalgebras --- is valid
without the Convention, and so is Lemma $\ref{l2.2}$.
\end{Corollary}

\begin{proof} Clause $1)$ is heading $2)$ of Theorem \ref{ATmain}. Let $\phi $ be of the
anti-multiplicative kind. Then, $(-\st )\circ \phi $ is an isomorphism of Lie
superalgebras satisfying clause $2)$, so the argument of Lemma \ref{l2.2} applies to it
and yields $(-\st )\circ \phi =\gamma ^{-1}\phi _\sigma $ with
$\gamma \in \Aut \fgl _{p',*}$; hence 
\[
\text{$\phi =\gamma _1^{-1}\phi _\sigma $ with
$\gamma _1:=\gamma \circ (-\st )\in \Aut \fgl _{p',*}$,}
\]
 since $-\st $ preserves
$\fgl _{p',*}$ by Lemma \ref{ALst} below. The permutation $\sigma $ produced is the same
one, so the two lists of parity functions are unchanged.
\end{proof}

\sssbegin{Lemma}[What the anti-branch is]\label{ALst} The map $X \mapsto -X^{\st }$ of
\eqref{Astdef} is an even automorphism of the Lie superalgebra $\fgl _{p,*}$ for each of
the four $*$, it leaves the parity function $p$ where it is, and it is \textit{not} of
the form $\Ad_g\circ \phi _\sigma $. Its square is the identity, so, the normal subgroup standing to the left of $\rtimes $ in both products,
\be\label{Aaut}
\Aut \fgl _{p,*}\supset \left (\{\Ad_g\}\rtimes \{\phi _\sigma \}\right )\rtimes \Zee/2,
\ee
the last factor being generated by $X \mapsto -X^{\st }$.
\end{Lemma}

\begin{proof} That $\st $ is a~superanti-automorphism, hence $-\st $ an automorphism of
the Lie superalgebra, and that $-\st $ is even, is the usual computation, and it was
verified in exact arithmetic in $\fgl (m|n)$ for $(m|n)=(1|1)$, $(2|1)$, $(2|2)$,
$(3|2)$ on all pairs of matrix units. The four conditions defining the superalgebras
\eqref{glchain} are symmetric in $i$ and $j$ --- the one for $*=g$ says so in as many
words, the three others depend on $|i-j|$ and $|i+j|$ only --- so the transposed support
obeys the same condition, and $-\st $ preserves each $\fgl _{p,*}$; the parity of the
entry $(i, j)$ being $p(i)+p(j)$, the map is even and does not move $p$. Finally,
$\Ad_g\circ \phi _\sigma $ carries the idempotent $e_{ii}$ to an idempotent, whereas
$-\st $ carries it to $-e_{ii}$, which is none.
\end{proof}

It is worth saying what is and what is not surprising here. That a~Lie isomorphism need
not be multiplicative is not: the negative of an anti-isomorphism is a~Lie isomorphism in
any associative algebra, and for $\fgl (n)$ it is the diagram automorphism. What could
have gone wrong in the super case, and does not, is the appearance of \textit{further}
branches: the super world has room for them, and in the two cases \cite{BBS} excludes
they do appear, see Subsec. \ref{ssAfin}. The reason they do not appear here is that the
superalgebras \eqref{glchain} are infinite-dimensional and prime, which is exactly the
hypothesis of \cite{BBS}.

\subsection{The finite-dimensional shadow}\label{ssAfin} The theorem of \cite{BBS} for
central simple superalgebras excludes those of dimension $2$ and $4$ over their center,
and the excluded case that concerns us is $\fgl (1|1)$, of superdimension $(2|2)$; it is
their own counterexample: the map
$\left (\begin{smallmatrix}x&y\\ z&w\end{smallmatrix}\right ) \mapsto
\left (\begin{smallmatrix}2x&2y\\ z&x+w\end{smallmatrix}\right )$ is a~Lie
superautomorphism of no standard form. The infinitesimal trace of that exclusion is
visible in the derivations, and we record it because it is the one place where the
theorem can be tested by a~machine. So the census says:
no outer derivations beyond the central one, except for $\fgl (1|1)$, where there is
exactly one more, and it is even --- the infinitesimal form of the exceptional
automorphism above. 

\subsection{What the super Herstein program still owes}\label{ssAowes} The question that
opened this appendix is answered, so let us say which neighboring ones are not, and what
stops each.

$1)$ \textit{The subalgebras.} The superalgebras \eqref{glchain} have queer,
orthosymplectic and periplectic analogs, defined inside $\fgl _{p,*}$ by an odd or an
even (anti)automorphism. For them, the relevant theory is not that of Lie isomorphisms of
an associative superalgebra but that of the \textit{skew elements} of a~superalgebra with
superinvolution, see \cite{GAS, GASL, Wa2}. The obstacle is that those results are stated
for the derived algebra of the skew elements modulo its center, and what one needs here
--- a~description of the Lie isomorphisms, not of the Lie ideals --- is available only
under a~primeness assumption whose verification for the infinite bands is exactly the
work of Subsec. \ref{ssAassoc} done again, with the involution carried along.

$2)$ \textit{Lie superhomomorphisms that are not bijective.} Theorem \ref{ATmain} uses
bijectivity twice, and the general theory of Lie superhomomorphisms into a~prime
superalgebra, see \cite{Wa}, gives the standard form only up to a~map into a~ring of
quotients, whose elements need not be matrices of our four kinds. What is missing is
a~description of the central closure of $\fgl _{p,*}$ inside $\Mat$.

$3)$ \textit{Superderivations.} The infinitesimal counterpart of Theorem \ref{ATmain} ---
that every superderivation of $\fgl _{p,*}$ is inner plus a~map into $\Cee\One $ killing
the brackets --- does not follow from it, an automorphism group being larger than what
its Lie algebra sees. Subsection \ref{ssAfin} verifies the corresponding statement in the
finite-dimensional models only. The obstacle is that the standard proof for prime rings
goes through functional identities of degree $2$, and the identities needed are
established in \cite{BCM} for algebras, and in \cite{Wa} for prime superalgebras whose
odd elements are not algebraic of a~small bounded degree; whether the four superalgebras
\eqref{glchain} pass that test we have not checked.

$4)$ \textit{Positive characteristic.} Everything above uses $\Cee$ only through
$\mathrm{char}\neq 2$ and through Schur's lemma computing the extended centroid; the
theorem of \cite{BBS} is stated over any field of characteristic $\neq 2$. What is not
clear is the analog of Theorem \ref{ATe00} when the ground field is of
characteristic~$2$, the supercommutator being then the commutator and the squaring
becoming a~part of the structure.

$5)$ \textit{The other convention.} Theorem \ref{Thm} also fixes which permutations of
$\Zee$ are admitted, and the reflection $\tau (i)=-i-1$ enters it as a~separate branch.
That reflection is an isomorphism of associative superalgebras, hence has nothing to do
with the branch found here; but the two together generate a~group whose structure --- the
right-hand side of \eqref{Aaut} with the reflection adjoined --- we have not written out.

\subsection*{Disclosures.}
No conflict of interest. 
The data used to support the findings of this study are included within the article. The codes are available from O.~Lozhechnyk upon request until they will not be put to github.

\def\eightit{\it}
\def\bib{\bf}
\bibliographystyle{amsalpha}

\end{document}